\documentclass[runningheads]{llncs}
\usepackage{epsfig}
\usepackage{amsmath}
\usepackage{graphicx}
\usepackage{amsfonts}

\def\bt{\begin{theorem}\rm}
\def\et{\end{theorem}}
\def\bc{\begin{corollary}\rm}
\def\ec{\end{corollary}}
\def\bproof{\begin{proof}}
\def\endproof{\end{proof}}
\def\bx{\begin{example}}
\def\ex{\end{example}}
\def\bl{\begin{lemma}}
\def\el{\end{lemma}}
\def\bp{\begin{proposition}}
\def\ep{\end{proposition}}

\def\bdf{\begin{definition}\rm}
\def\edf{\end{definition}}
\def\ba{\begin{array}}
\def\ea{\end{array}}
\def\be{\begin{equation}}
\def\ee{\end{equation}}
\def\bd{\begin{description}}
\def\ed{\end{description}}
\def\bu{\begin{enumerate}}
\def\eu{\end{enumerate}}
\def\bi{\begin{itemize}}
\def\ei{\end{itemize}}

\newbox\bigstrutbox
\setbox\bigstrutbox=\hbox{\vrule height12.5pt depth5pt width0pt}
\def\bigstrut{\relax\ifmmode\copy\bigstrutbox\else\unhcopy\bigstrutbox\fi}

\newbox\Bigstrutbox
\setbox\Bigstrutbox=\hbox{\vrule height17.5pt depth5pt width0pt}
\def\Bigstrut{\relax\ifmmode\copy\Bigstrutbox\else\unhcopy\Bigstrutbox\fi}

\def\endproof{{\hfill{\phantom{123}}\hfill\Large $\bf\bigstrut \Box$}\vskip .15cm}
\def\endproofs{{\vskip -.8cm \endproof}}
\def\seq{{\,\stackrel{*}{=}\,}}

\def\ds{\displaystyle}

\def\A{{\bf A}}
\def\B{{\bf B}}
\def\C{{\bf C}}
\def\D{{\bf D}}
\def\E{{\bf E}}
\def\F{{\bf F}}
\def\G{{\bf G}}
\def\H{{\bf H}}
\def\I{{\bf I}}
\def\J{{\bf J}}
\def\K{{\bf K}}
\def\L{{\bf L}}
\def\M{{\bf M}}
\def\N{{\bf N}}

\def\P{{\bf P}}
\def\Q{{\bf Q}}
\def\R{{\bf R}}
\def\S{{\bf S}}
\def\T{{\bf T}}
\def\U{{\bf U}}
\def\V{{\bf V}}
\def\W{{\bf W}}
\def\X{{\bf X}}
\def\Y{{\bf Y}}

\def\a{{\bf a}}
\def\b{{\bf b}}
\def\c{{\bf c}}
\def\d{{\bf d}}
\def\e{{\bf e}}
\def\f{{\bf f}}

\def\l{{\bf l}}
\def\m{{\bf m}}
\def\n{{\bf n}}

\def\p{{\bf p}}
\def\q{{\bf q}}

\def\s{{\bf s}}
\def\t{{\bf t}}
\def\u{{\bf u}}
\def\v{{\bf v}}
\def\w{{\bf w}}
\def\x{{\bf x}}
\def\y{{\bf y}}
\def\z{{\bf z}}

\def\0{{\bf 0}}
\def\1{{\bf 1}}
\def\2{{\bf 2}}
\def\3{{\bf 3}}
\def\4{{\bf 4}}
\def\5{{\bf 5}}
\def\6{{\bf 6}}
\def\7{{\bf 7}}
\def\8{{\bf 8}}
\def\9{{\bf 9}}
\def\tr{{\rm tr}}

\def\ce{{\e}^*}
\def\cn{{\n}^*}
\def\ca{{\a}^*}
\def\RS{{\mathbb R}^{3,3}}
\def\RF{{\mathbb R}^{4}}
\def\RT{{\mathbb R}^{3}}

\def\ds{\displaystyle}
\def\diag{{\rm diag}}
\def\CL{{{\cal C}\hskip -.02cm l}}

\begin{document}

\pagestyle{headings}

\mainmatter

\title{Three-Dimensional Projective Geometry with Geometric Algebra}

\titlerunning{Three-Dimensional Projective Geometry with Geometric Algebra}

\author{Hongbo Li, Lei Huang, Changpeng Shao, Lei Dong}

\authorrunning{Hongbo Li, Lei Huang, Changpeng Shao, Lei Dong}

\institute{Academy of Mathematics and Systems Science\\
           Chinese Academy of Sciences\\
           Beijing 100190, China\\
           hli@mmrc.iss.ac.cn}

\maketitle

\begin{abstract}

The line geometric model of 3-D projective geometry has the nice property that the Lie algebra $sl(4)$ 
of 3-D projective transformations is isomorphic to the bivector algebra of $\CL(3,3)$, and line geometry is closely
related to the classical screw theory for 3-D rigid-body motions. The canonical homomorphism from $SL(4)$ to
$Spin(3,3)$ is not satisfying because it is not surjective, and the projective transformations
of negative determinant do not induce orthogonal transformations in the Pl\"ucker coordinate space of lines.

\vskip .2cm
This paper presents our contributions
in developing a rigorous and convenient algebraic framework for the study of 3-D projective geometry with Clifford algebra. To
overcome the unsatisfying defects of the Pl\"ucker correspondence, we propose a group $Pin^{sp}(3,3)$
with $Pin(3,3)$ as its normal subgroup, to quadruple-cover the group of 
3-D projective transformations and polarities. We construct spinors in factored form that generate 3-D reflections and
rigid-body motions, and extend screw algebra from the Lie algebra of rigid-body motions to
other 6-D Lie subalgebras of $sl(4)$, and construct the corresponding cross products and virtual works.

\end{abstract}

\section{Introduction}
\setcounter{equation}{0}

The study of the geometry of lines in space was invented by Pl\"ucker with his
introduction of the now so-called Pl\"ucker coordinates of lines. It became an active
research topic with the establishment of screw theory by Balls \cite{balls}, where the 6-D
Pl\"ucker coordinates of a line are decomposed into a pair of 3-D vectors,
called the {\it screw form} of the line, and the inner product and
cross product of vector algebra are extended to screw forms. 

A pair of force and torque, called a {\it wrench}, 
are naturally represented by a pair of 3-D vectors, and are geometrically interpreted as a line in space along which the force acts,
together with a line at infinity about which the torque acts. On the other hand, 
a pair of infinitesimal rotation and translation, called an infinitesimal {\it screw motion} or {\it rigid-body motion} 
or {\it twist}, are represented by
the rotation axis and the translation vector, and are again naturally represented by a pair of 3-D vectors. Geometrically
the translation is a special ``rotation" about an axis that is at infinity, so the translation vector represents a line at infinity.
For a wrench $(\f, \q)$ and an infinitesimal screw motion $(\v, \t)$, where $\f$ is the force direction multiplied with
the magnitude of force, $\q$ is the composed torque,
$\v$ is the rotation axis direction multiplied with the angle of rotation, and $\t$ is the moment of the 
screw motion, the {\it virtual work} of the wrench along the infinitesimal screw motion is the ``crossed" inner product
\be
\left(\ba{c}
\f\\
\q
\ea\right)\cdot 
\left(\ba{c}
\v\\
\t
\ea\right)
:=
\f\cdot \t+\q\cdot \v.
\label{def:innerproduct6}
\ee

The inner product (\ref{def:innerproduct6}) gives the 6-D space of wrenches a signature $\RS$, where a pure force has zero inner product with
itself, called a {\it null vector}. A {\it positive} (or {\it negative}) vector of $\RS$ is one having positive (or negative) inner product
with itself. A positive vector is interpreted as a pure force together with an extraneous torque so that the pair follow the right-hand rule,
while for the negative vector, the force and torque follow the left hand rule. The group $SL(4)$ which
acts in the 4-space of homogeneous coordinates of points, can be lifted to a group action in the 6-D space of wrenches by acting upon
the Pl\"ucker coordinates of the lines representing the wrenches. The image of the lift is $SO_0(3,3)$, the connected component of $SO(3,3)$
containing the identity \cite{beffa}. As $PSL(4)=SL(4)/Z_2$ is the group of orientation-preserving projective transformations,
the crossed inner product provides a 6-D orthogonal geometric model of wrenches to study 3-D projective geometry of points.

The same inner product (\ref{def:innerproduct6}) also gives the 6-D space of twists the same signature
$\RS$. The interpretation of a null vector of $\RS$ in the setting of twists, is that it represents an infinitesimal pure 
rotation or pure translation. A positive vector represents an infinitesimal screw motion where the translation along the screw axis
follow the right-hand rule with the orientation of the rotation, while a negative vector represents a left-handed infinitesimal
screw motion. The lift of group $SL(4)$ to $SO_0(3,3)$ then makes 
the space of infinitesimal {\bf rigid-body motions} a 6-D orthogonal geometric model to 
study the orientation-preserving {\bf projective geometry} of points.

The bold-faced words clearly reveal a conflict. The group of rigid-body motions is 6-D, while the group $SL(4)$ is 15-D; 
the former is much smaller. What sense does it make to investigate projective transformations via rigid-body motions? 
Furthermore, the inner product (\ref{def:innerproduct6}) is between the space of wrenches and the space of twists, indicating
that the two spaces need to be identified, yet they have to be different spaces by nature. Understanding (\ref{def:innerproduct6})
as a pairing between a linear space and its dual space does not make much difference, as the same
inner product exists in either space. 

Line geometry and screw theory are closely related to each other.
In history, the screw forms were first used by Clifford in the name of biquaternions, also known as
dual quaternions, in describing 3-D Euclidean transformations. Later on, Balls \cite{balls}, Study
\cite{study}, Blaschke\cite{blaschke} established screw theory and developed {\it dual vector algebra} out of
Clifford's dual quaternions, also known as {\it screw algebra}. Nowadays
line geometry together with screw theory have important applications in mechanism
analysis, robotics, computer vision and computational geometry
\cite{chevallier},
\cite{bayro},
\cite{selig},
\cite{daniilidis},
\cite{merlet},
\cite{pottmann},
\cite{dai}.

For two vectors $(\x_1, \y_1)^T\in \RS$ and $(\x_2, \y_2)^T\in \RS$, where $\x_i, \y_j\in \RT$, their {\it cross product}, also called
{\it dual vector product}, is defined as follows:
\be
\left(\ba{c}
\x_1\\
\y_1
\ea\right)\times
\left(\ba{c}
\x_2\\
\y_2
\ea\right)
:=\left(\ba{c}
\x_1\times \x_2\\
\x_1\times \y_2+\y_1\times \x_2
\ea\right).
\ee
This product is covariant under the subgroup of $SO_0(3,3)$ that is the lift of the group of rigid-body motions $SE(3)$, but not so
under the whole group $SO_0(3,3)$. In other words, it is not a valid operator in 3-D projective geometry; it is valid only for
Euclidean geometry. 

In dual vector algebra, the {\it dual inner product} of two vectors of $\RT$ is defined to be a {\it dual number}. 
A dual number is of the form $\lambda +\epsilon \mu$ where $\lambda, \mu\in \mathbb R$,\, $\epsilon^2=0$ and $\epsilon$ commutes with everything.
This numbers system is a ring instead of a field, and the corresponding polynomials and modules are drastically different from the usual ones.
The dual inner product is invariant under the lift of $SE(3)$ to $SO_0(3,3)$, so it is suitable for Euclidean geometry only.

In \cite{li01}, it was pointed out that 
dual vector algebra and dual quaternions can be realized in the conformal geometric algebra $\CL(4,1)$,
and can be extended to arbitrary dimensions. The Euclidean geometric parts of the wrench model and the twist model 
have no conflict, and their identification is natural. The twist model should not have anything beyond
Euclidean geometry, otherwise it would be absurd. The wrench model, or more generally the model of lines in space, 
deserves further attention.

For invariant computing in projective geometry, the traditional algebraic tool is Grassmann-Cayley algebra and bracket algebra
\cite{li08}. The study of projective geometry by Clifford algebra was initiated by Hestenes and Ziegler \cite{hestenes91}, and
Stolfi \cite{stolfi}. 
The representation of projective transformations by spinors was initiated by Doran {\it et al.} \cite{doran}, where a homomorphism of 
the Lie algebra $gl(n)$ into $so(n,n)$ was proposed, making it possible to construct projective transformations by elements of
$Pin(n,n)$. Following this line, Goldman and Mann \cite{goldman14} discovered for many 3-D projective transformations their bivector
generators in $\CL(4,4)$. 
Considering that the dimension of $so(4,4)$ is $C_8^2=28$, while the dimension of $sl(4)$ is 15, the embedding space of $sl(4)$ seems
too high \cite{dorst15}.

A classical result \cite{beffa} states that the group $SL(4)$, which
acts upon the 4-space of homogeneous coordinates of points, is in fact isomorphic to the group $Spin_0(3,3)$, the connect component
of $Spin(3,3)$ containing the identity, and the isomorphism is realized via the Pl\"ucker coordinates of lines and the adjoint action of
$Spin(3,3)$ upon $\RS$. This canonical isomorphism indicates the possibility of using the wrench model, the model of spatial lines, to study
projective geometry with $\CL(3,3)$. 

In AGACSE 2009, Li and Zhang \cite{li09} proposed a new model of 3-D projective geometry
by taking the null vectors of ${\mathbb R}^{3,3}$ as algebraic generators, and defining 
points and planes as the two connected
components of the set of null 3-spaces of ${\mathbb R}^{3,3}$ respectively.
Whenever an element of $Spin_0(3,3)$ acts upon $\RS$, it induces a projective transformation via the outermorphism
of the action upon the null 3-vectors representing 3-D points and planes.
This approach was later followed by Klawitter \cite{klawitter14}, who proposed an explicit expression of the spinor 
inducing a projective transformation in $4\times 4$ matrix form, and recently by 
Dorst \cite{dorst15}, who constructed bivector generators for many 3-D projective transformations.

When viewed from the homogeneous coordinates model $\RF$ of 3-D projective geometry, the ${\mathbb R}^{3,3}$ model seems to have too many defects.
The map from $SL(4)$ to $SO(3,3)$ is not surjective, nor injective. The 
projective transformations of negative determinant cannot be lifted to $O(3,3)$, and conversely, the elements of
$O(3,3)$ with negative determinant do not correspond to any projective transformation, but represent projective polarities
where points are all mapped to planes. In the ${\mathbb R}^{3,3}$ model, while lines are represented by vectors,
the 3-D points and planes are represented by null 3-vectors, whose embedding vector space has dimension $C_6^3=20$. 
To make things worse,
the mapping from $\RF$ to the null 3-vectors is quadratic, and defining the subset of null 3-vectors in the 20-D vector space
they span is difficult.

From the mathematical viewpoint, establishing the space $\RF$ of homogeneous coordinates from the 6-space
${\mathbb R}^{3,3}$ spanned by lines requires rigorous mathematical argument. It is the converse procedure of Pl\"ucker's
construction of line coordinates from point coordinates. The well-definedness of the points and planes, and the covariance of the
construction under suitable transformations of ${\mathbb R}^{3,3}$ need to be established. The benefits of using the null 3-vectors
instead of the linear space $\RF$ to represent points need to be discovered. The groups $SO_0(3,3)$
and $Spin_0(3,3)$ are too small to cover the whole group of all 3-D projective transformations and polarities, and finding 
suitable covering groups to provide spin representations for all 3-D projective transformations and polarities is indispensable.

So compared with other models of Geometric Algebra for classical geometries 
\cite{hestenes84},
\cite{hestenes98},
\cite{lasenby},
\cite{dorst07},
\cite{li08},
\cite{hildenbrand},
\cite{fleuystad},
\cite{kanatani}, the line geometric model of 3-D projective geometry is much less developed. 
When every problem raised above is solved, then for the group of 3-D Euclidean transformations, a highly mature subject 
of study in Geometric Algebra, a one-to-one correspondence among the representations in the line geometric model and in other
models need to be set up. 

As mentioned before, 
the screw algebra is valid only for Euclidean motions, and the corresponding group $SE(3)$ is only a subgroup of $SL(4)$.
When $SE(3)$ is replaced by another 6-D Lie subgroup of $SL(4)$, then the Lie algebra $se(3)$ of Euclidean motions is replaced
by another 6-D Lie subalgebra of $sl(4)$. Correspondingly, we can introduce new screw forms for the 6-D Lie subalgebra,
together with the new ``virtual work" of a wrench, which is still a vector of $\RS$, along an infinitesimal
``projective motion" represented by a screw form of the Lie subalgebra. The 6-D Lie subalgebra has its own Lie bracket,
so the corresponding screw forms should have a different cross product. The new ``virtual work" should be related to the new
cross product, or even be completely determined by it.

We can go one step further by
decomposing the 15-D algebra $sl(4)$ into the direct sum of five 3-D vector spaces, so that instead of using only pairs
of 3-D vectors as in classical screw theory for the screw forms of $se(3)$, we can use 5-tuples of 3-D vectors to represent
screw forms of $sl(4)$, and develop a ``super-screw theory", equipped with ``super-cross product" and ``super-virtual work".

For the purpose of developing a mathematically rigorous model out of the peculiar, unfamiliar and seemingly ineffective
line geometric model, for more effectively describing and manipulating 3-D projective transformations with Geometric Algebra, 
and with the ambition to further extend screw theory to projective geometry, we
picked up the research subject again in 2014, and after one-year's hard work, we are confident to announce that
from the algebraic viewpoint, this model is sufficiently mature now. The main contributions are summarized
as follows:

1. Rigorous establishment of the $\RS$ model for 3-D projective geometry.

While in the classical model of
projective line geometry only $SL(4)$ has spin representation, and all the spinors are
in $Spin_0(3,3)$, a rather unsatisfying limitation, the new model completely
overcomes the limitation by providing pin group representations for all 3-D projective
transformations and polarities, thus enlarging the transformation group four times. 

The group of {\it linear regularities} of $\RS$ is defined by
\be
RL(3,3):=\{\B\in GL(3,3)\,|\, \B^T\J\B=\pm \J\}, 
\ee
where $\J$ is the matrix form of the metric of $\RS$.
Only when we computed the group acting upon the null 3-vectors induced by $RL(3,3)$ did we find the complete version of 
the line geometric model. The group $RL(3,3)$ double covers the whole group of 3-D projective transformations and
polarities in this manner, and the group $Pin^{sp}(3,3)$ quadruple-covers the latter, hence it
can be used to construct versors for all kinds of 3-D projective transformations and polarities.

The well-definedness of points and planes in the $\RS$ model, and the covariance of the representations 
are established. Some nice properties of reflections in $\RS$ are found, together with the classification of 3-D projective transformations
induced by two reflections in $\RS$.

2. Construction of spinors in factored form inducing 3-D reflections and rigid-body motions, and discovery of the relation between the 
cross product of the screw forms of $se(3)$ and the virtual work.

For 3-D reflections and rigid-body motions, the spinors inducing them in $Pin^{sp}(3,3)$ in factored form
are discovered. Since the bivector Lie algebra of $\CL(3,3)$ is isomorphic to $sl(4)$, any element of $se(3)$ has a bivector form,
and the cross product of the bivectors equals the cross product of their screw forms as vectors of $\RS$. 
On the other hand, a wrench is only a vector of $\RS$. To make pairing with a bivector, a vector needs to be first upgraded
to a bivector of $\Lambda^2(\RS)$ by making inner product (tensor contraction) with a trivector. 
We show that this trivector is exactly the one complementary to the trivector defining the cross product of the screw forms of $se(3)$,
and the latter trivector is exactly the lift of the quadratic form of ${\mathbb R}^{3,0,1}$.
This correspondence shows the intrinsic the connection between the virtual work and the cross product of $se(3)$. 
The connection between the wrench interpretation and the twist interpretation of the line geometric model is now clarified.

3. Extension of the cross product and virtual work of the screw forms of $se(3)$ to other 
6-D Lie subalgebras of $sl(4)$.

For many 6-D Lie subalgebras of $sl(4)$, we have developed the corresponding screw forms together with the cross product and virtual work
that are completely determined by the Lie bracket of the subalgebra. In particular, for $so(\K)$ where $\K$ is a quadratic form
of $\RF$ with rank $\geq 3$, we have established the corresponding screw forms, cross products and virtual works, 
and discovered a striking fact: the trivectors for constructing
new cross products and new virtual works are exactly the lifts of the quadratic form $\K$ by the Pl\"ucker correspondence and the
dual Pl\"ucker correspondence to the trivector space. This result demonstrates that there is no intrinsic connection between the 
$se(3)$-interpretation and the wrench interpretation of line geometry.

This paper is organized as follows. Section 2 is on the Pl\"ucker model of 3-D projective geometry, the invariant group $PR(3)$
and its quadruple-covering group $Pin^{sp}(3,3)$. Section 3 is on the covariance of the Pl\"ucker transform and the dual
Pl\"ucker transform, and the explicit expression of the induced group element of $PR(3)$ by an element of $RL(3,3)$. Section 4 is on
properties of projective transformations induced by one or two reflections of $\RS$, and the construction of spinors in factored form inducing
rigid-body motions. Section 5 is on the bivector representation of $se(3)$, the screw forms and the cross products of screw forms.
For reflections in space, the generating spinors are also constructed. Section 6 is on extension of classical screw theory to
6-D Lie subalgebras of $sl(4)$ other than $se(3)$. A framework of super-screw theory is also presented.

\section{The Pl\"ucker model of 3-D projective geometry}
\setcounter{equation}{0}

\bdf
A real $n$-space is said to have {\it signature ${\mathbb R}^{p,q,r}$}, 
if it has a basis with respect to which
the metric of the $n$-space is diagonal, where the diagonal elements are composed of 1 of multiplicity $p$, and $-1$ of multiplicity $q$, and 
0 of multiplicity $r=n-p-q$.
${\mathbb R}^{p,q}$ stands for ${\mathbb R}^{p,q,0}$, called a {\it non-degenerate} inner-product space;
${\mathbb R}^{p}$ stands for ${\mathbb R}^{p,0,0}$, called a {\it Euclidean} inner-product space;
${\mathbb R}^{0,0,r}$ is called a {\it null} inner-product space, and
${\mathbb R}^{p,1}$ is called a {\it Minkowski} inner-product space.
\edf

\bdf
A {\it Witt decomposition} of $\RS$ refers to a decomposition of $\RS$ into the direct sum of two null 3-spaces, say $\I_3, \J_3$.
A {\it Witt basis} of the decomposition is composed of a basis of $\I_3$, say $\E_1, \E_2, \E_3$, and the corresponding
{\it Witt-dual basis} of $\J_3$, denoted by $\E_1', \E_2', \E_3'$, such that $\E_i\cdot \E_j'=\delta_{ij}$.
\edf

\bl \cite{crumeyrolle}
$\RS$ has infinitely many Witt decompositions.
Let $\RS=\I_3\oplus \J_3$ be a fixed Witt decomposition of $\RS$. 
Then for any basis $\E_1, \E_2, \E_3$ of the 3-space $\I_3$,
the corresponding Witt-dual basis $\E_1', \E_2', \E_3'$ of the 3-space $\J_3$ is unique.
\el

\bl
\label{lem:ij}
Fix a Witt decomposition $\RS=\I_3\oplus \J_3$.
For any null 3-space $\S_3$ of $\RS$, if the dimension $n$ of $\S_3\cap \I_3$ is even (or odd), then
the dimension $m$ of $\S_3\cap \J_3$ is odd (or even). 
\el

{\it Proof.}
When $n=0$, then $\RS=\I_3\oplus \S_3$. Let $\E_1, \E_2, \E_3$ be a basis of $\I_3$, and let the 
corresponding Witt-pairing bases in $\J_3$ and $\S_3$ be respectively $\E_1', \E_2', \E_3'$ and
$\s_1, \s_2, \s_3$. Then $\s_i=\E_i'+\sum_j s_{ij}\E_j$ where $s_{ij}=-s_{ji}$. Let 
\[
\S=\left(\ba{ccc}
0 &\ s_{12} &\ s_{13} \\
-s_{12} &\ 0 &\ s_{23} \\
-s_{13} &\ -s_{23} &\ 0
\ea\right).
\]
If $\S=0$ then $\S_3=\J_3$ and $m=3$. If $\S\neq 0$, then its rank is 2, so its kernel has dimension 1.
Let $\X=x_1\s_1+x_2\s_2+x_3\s_3\in  \S_3\cap \J_3$, then 
$0=\sum_i x_i(\s_i-\E_i')=\sum_{i,j} x_i s_{ij}\E_j$, so $(x_1,x_2,x_3)^T$ is in the kernel of $\S^T$.
This proves $m=1$.

When $n=2$, let $\I_3$ be spanned by $\E_1, \E_2, \E_3$, and let $\S_3$ be spanned by $\E_1, \E_2, \s$.
Then $\s=s_3\E_3+\sum_i s_i'\E_i'$. Since $\S_3$ is a null 3-space, by $\s\cdot \E_1=\s\cdot \E_2=\s^2=0$,
we get $s_1'=s_2'=s_3s_3'=0$. If $s_3'=0$ then $\S_3=\I_3$, violating the assumption that $n=2$. So
$s_3'\neq 0$ and $s_3=0$. We get $\s=s_3'\E_3'$ and the 1-space spanned by $\E_3'$ is $\S_3\cap \J_3$.
This proves $m=1$.

When $n=1$, let $\S_3$ be spanned by $\E_1, \s, \t$, such that $\s=s_2\E_2+s_3\E_3+s_2'\E_2'+s_3'\E_3'$,
and $\t=t_2\E_2+t_3\E_3+t_2'\E_2'+t_3'\E_3'$. From $\s^2=\t^2=\s\cdot \t=0$, we get 
\be
s_2s_2'+s_3s_3'=t_2t_2'+t_3t_3'=s_2t_2'+s_3t_3'+t_2s_2'+t_3s_3'=0.
\label{inner:st}
\ee
Since the 1-space spanned by $\E_1$ is $\I_3\cap \S_3$, $\left|\ba{cc}
s_2' &\ s_3' \\
t_2' &\ t_3'
\ea\right|\neq 0$. 
In the 2-space spanned by vectors $(s_2', s_3')^T$ and $(s_2', s_3')^T$, by (\ref{inner:st}),
we have $(s_2, s_3)^T=\lambda (s_3', -s_2')^T$, and $(t_2, t_3)^T=\mu (t_3', -t_2')^T$, and
$\lambda=\mu$. If $\lambda=0$ then $\S_3$ is spanned by $\e_1, \e_2', \e_3'$, so $m=2$. 
If $\lambda\neq 0$, then 
$\left|\ba{cc}
s_2 &\ s_3 \\
t_2 &\ t_3
\ea\right|\neq 0$, so $\S_3\cap \J_3=\{0\}$, and $m=0$.

When $n=3$, the conclusion $m=0$ is trivial.
\endproof

\bc
\label{cor:standard}
For a Witt decomposition $\RS=\I_3\oplus \J_3$ with Witt basis $\E_1, \E_2, \E_3$ and
$\E_1', \E_2', \E_3'$, if $\S_3\cap \I_3$ is the 1-space spanned by $\E_1$, then $\S_3$ is spanned by
$\E_1, \s, \t$, where  
\be
\s=\lambda \E_2+\mu \E_3',\ \ \
\t=\lambda \E_3-\mu \E_2',
\ee 
for some $\mu\neq 0$ and $\lambda$. If $\S_3\cap \I_3$ is the 2-space spanned by $\E_1, \E_2$, 
then $\S_3$ is spanned by $\E_1, \E_2, \E_3'$.
\ec

\bl
Let the intersection of null 3-spaces $\P_3, \Q_3$ be a 1-space, then for any null 3-space $\S_3$ of $\RS$,
the dimension $n$ of $\S_3\cap \P_3$ and the dimension $m$ of $\S_3\cap \Q_3$ have the same parity.
\el

{\it Proof.} 
Let $\P_3\cap \Q_3$ be the 1-space spanned by vector $\E_1$. 
Assume $\S_3\neq \P_3$ and $\S_3\neq \Q_3$. 

When $n=0$, for the Witt decomposition $\RS=\P_3\oplus \S_3$, by Lemma \ref{lem:ij}, 
$m=0$ or 2.

When $n=1$ and $\P_3\cap \Q_3=\P_3\cap \S_3$, let $\E_1, \E_2, \E_3$ be a basis of $\P_3$, and let
$\E_1', \E_2', \E_3'$ be the Witt-pairing basis of $\P_3'$ for a Witt decomposition $\RS=\P_3\oplus \P_3'$.
By Corollary \ref{cor:standard}, $\Q_3$ is spanned by 
$\E_1, \lambda_q \E_2+\mu_q \E_3', \lambda_q \E_3-\mu_q \E_2'$, where $\mu_q\neq 0$, while $\S_3$ is spanned
by $\E_1, \lambda_s \E_2+\mu_s \E_3', \lambda_s \E_3-\mu_s \E_2'$, where $\mu_s\neq 0$. Since $\S_3\neq \Q_3$,
$\lambda_q:\mu_q\neq \lambda_s:\mu_s$, so the 1-space $\E_1$ is the only intersection of $\S_3$ and $\Q_3$.
This proves $m=1$.

When $n=1$ but $\P_3\cap \Q_3\neq \P_3\cap \S_3$, let $\E_1, \E_2, \E_3$ be a basis of $\P_3$ such that
$\P_3\cap \S_3$ is the 1-space spanned by $\E_2$. Let 
$\E_1', \E_2', \E_3'$ be the Witt-pairing basis of $\P_3'$ for a Witt decomposition $\RS=\P_3\oplus \P_3'$.
By Corollary \ref{cor:standard}, $\Q_3$ is spanned by 
$\E_1, \lambda_q \E_2+\mu_q \E_3', \lambda_q \E_3-\mu_q \E_2'$, where $\mu_q\neq 0$, 
while $\S_3$ is spanned
by $\E_2, \lambda_s \E_1+\mu_s \E_3', \lambda_s \E_3-\mu_s \E_1'$, where $\mu_s\neq 0$. 
Then $\S_3\cap \Q_3$ is the 1-space spanned by $\lambda_s\mu_q\E_1+\lambda_q\mu_s\E_2+\mu_s\mu_q\E_3'$.
Again $m=1$.

When $n=2$ and $\P_3\cap \Q_3\subset \P_3\cap \S_3$, let $\E_1, \E_2, \E_3$ be a basis of $\P_3$ such that
$\P_3\cap \S_3$ is spanned by $\E_1, \E_2$. Let 
$\E_1', \E_2', \E_3'$ be the Witt-pairing basis of $\P_3'$ for a Witt decomposition $\RS=\P_3\oplus \P_3'$,
then $\S_3$ is spanned by $\E_1, \E_2, \E_3'$, and $\Q_3$ is spanned by $\E_1, 
\lambda \E_2+\mu \E_3', \lambda \E_3-\mu \E_2'$ where $\mu\neq 0$. Then $\Q_3\cap \S_3$ is spanned by
$\E_1, \lambda \E_2+\mu \E_3'$, and $m=2$.

When $n=2$ but $\P_3\cap \Q_3$ is not in $\P_3\cap \S_3$, let $\E_1, \E_2, \E_3$ be a basis of $\P_3$ such that
$\P_3\cap \S_3$ is spanned by $\E_1, \E_2$, and $\P_3\cap \Q_3$ is spanned by $\E_3$. 
Let 
$\E_1', \E_2', \E_3'$ be the Witt-pairing basis of $\P_3'$ for a Witt decomposition $\RS=\P_3\oplus \P_3'$.
Then $\S_3$ is spanned by $\E_1, \E_2, \E_3'$, and $\Q_3$ is spanned by $\E_3, 
\lambda \E_1+\mu \E_2', \lambda \E_2-\mu \E_1'$ where $\mu\neq 0$. Obviously $\Q_3\cap \S_3=\{0\}$,
and $m=0$.
\endproof

\bp \label{prop:def:linegeo}
The set of null 3-spaces can be decomposed into two subsets (connected components):
in each subset the dimension of the intersection subspace of any two different elements 
is 1, and between the two subsets, the dimension of the intersection subspace of any two elements,
one from each subset, is 0 or 2.
\ep

As a corollary, the following concepts of points and planes in the set of
null 3-spaces of $\RS$ are well defined, 
and any pair of non-incident point and plane form a Witt decomposition of $\RS$.

\bdf 
\label{def:pts}
For fixed Witt decomposition $\RS=\I_3\oplus \J_3$, if we call 
$\I_3$ a point (the origin), and call $\J_3$ a plane (the plane at infinity), then
for any null 3-space $\S_3$ of $\RS$, let $n$ be the dimension of the vector space $\S_3\cap \I_3$,
if $n$ is even, then $\S_3$ is called a plane, and if $n$ is odd, then $\S_3$ is called a point.
\edf

The above concepts of points and planes have the following background. 
The classical
{\it Pl\"ucker map} changes a pair of points of Euclidean affine space ${\cal E}^3$ in their homogeneous coordinates form
$\X=(x_0, x_1, x_2, x_3)^T$ and 
$\Y=(y_0, y_1, y_2, y_3)^T$
with respect to the basis $\e_0, \e_1, \e_2, \e_3$ of $\RF$, to a vector $\X\Y:=\X\wedge \Y$ of $\RS$ where the induced basis is
$\e_{ij}=\e_i\e_j$ for $0\leq i<j\leq 3$,
and the image of the two points is represented by its {\it Pl\"ucker coordinates} with respect to the induced basis.

Henceforth we always denote the outer product of the Grassmann algebra 
$\Lambda(\RF)$ generated by $\RF$ by the juxtaposition of participating elements, while denoting the
the outer product of the Grassmann algebra $\Lambda(\RS)$ by the wedge symbol.

Let $\e_0, \e_1, \e_2, \e_3$ be an orthonormal basis of $\RF$. Let
\be\ba{lll}
\E_1=\e_{01},&\ \ \E_2=\e_{02},&\ \ \E_3=\e_{03}; \\

\E_1'=\e_{23},&\ \ \E_2'=\e_{31},&\ \ \E_3'=\e_{12}.
\ea
\label{induced:basis}
\ee
Then
\be
\I_3:=\E_1\wedge \E_2\wedge \E_3, \ \ \
\J_3:=\E_1'\wedge \E_2'\wedge \E_3' 
\ee
are two null 3-spaces of $\RS$ forming a Witt decomposition. 

Let $\e_0$ represent the origin of the Euclidean affine 3-space
${\cal E}^3$, and let $\RT=\langle \e_1, \e_2, \e_3\rangle$ represent the plane at infinity. 
The basis $\e_0, \e_1, \e_2, \e_3$ induces a Witt basis (\ref{induced:basis}). 
The two 3-vectors $\I_3, \J_3$ are invariant under any special linear transformation of $\RT$.

For $\X, \Y\in \RF$,
the vector $\X\Y\in \RS$, if not zero, is a null vector. Conversely, any null vector of $\RS$ is the image of either an affine line
or a line at infinity of ${\cal E}^3$ under the Pl\"ucker map. A 2-space $\S_2$ of $\RS$ spanned by null vectors, 
when interpreted geometrically so that its null
1-spaces are lines in space, has two kinds: (1) a pair of non-intersection lines in space, when the signature of $\S_2$
is ${\mathbb R}^{1,1}$; (2) a pair of incident point and plane, {\it i.e.}, a pencil of
lines incident at a fixed point and at the same time lying
on a fixed plane, when $\S_2$ is null.
A null 3-space $\S_3$ of $\RS$ when interpreted in line geometry, represents
either a point or a plane, {\it i.e.}, either all lines incident at a fixed point, or all lines lying
on a fixed plane.

Definition \ref{def:pts} is based on a fixed basis $\e_0, \e_1, \e_2, \e_3$ of $\RF$ and the induced Witt decomposition. 
Now Proposition \ref{prop:def:linegeo} tells us that no matter what the underlying 4-space $\RF$ is and what the induced Witt decomposition
of $\RS$ could be, as long as the fixed null 3-space $\E_1\wedge \E_2\wedge \E_3$ of $\RS$
is classified as a ``point", then in a
new line geometry whose abstract ``lines" are the null 1-spaces of $\RT$,
any ``point" defined with respect to the basis $\e_0, \e_1, \e_2, \e_3$ of the original $\RF$ is 
always classified as a point in the new line geometry.

Fix the underlying space ${\mathbb R}^4$ of the homogeneous coordinates of Euclidean affine space ${\cal E}^3$,
and fix a basis $\e_0, \e_1, \e_2, \e_3$ of it. The {\it projective transformation group} of ${\cal E}^3$ is the union
$SL(4)\cup SL^-(4)$, where $SL(4)$ is the linear transformations of determinant 1, while $SL^-(4)$ is the
linear transformations of determinant $-1$. In the dual space
$(\RF)^*$ of ${\mathbb R}^4$ equipped with the corresponding dual basis $\e_0^*, \e_1^*, \e_2^*, \e_3^*$ such that
the pairing between $\e_i$ and $\e_j^*$ is $\delta_{ij}$, 
the corresponding linear transformation of $\A\in GL(4)$ is $\A^{-T}$. The pair 
$(\A, \A^{-T})$ acts upon $\RF\times (\RF)^*$, and is still called a general linear transformation.

Any non-singular linear mapping from ${\mathbb R}^4$ to $(\RF)^*$ is called a 
{\it projective polarity}. The set of all projective polarities is denoted by $GP(4)$.
Such a mapping is called a {\it special polarity} if its matrix form $\A$ has determinant 1.
The matrix $\A^{-T}$ represents the corresponding linear mapping from $(\RF)^*$
to ${\mathbb R}^4$. The pair $(\A, \A^{-T})$ acts upon $\RF\times ({\mathbb R}^4)^*$, and is still called a 
projective polarity. 
The set of special polarities is denoted by $SP(4)$, and the set of
projective polarities with determinant $-1$ is denoted by $SP^-(4)$. 

Obviously, the set
\be
UR(4):=SL(4)\cup SL^-(4)\cup SP(4)\cup SP^-(4)
\ee
is a group, and $SL(4), SL(4)\cup SL^-(4), SL(4)\cup SP(4)$ are three subgroups.
$UR(4)$ is called the group of {\it unitary regularities}.

On the other hand, consider some general linear transformations of $\RS$.
Fix a Witt decomposition $\I_3\oplus \J_3$ and the corresponding Witt
basis $\E_1, \E_2, \E_3$, $\E_1', \E_2', \E_3'$ of $\RS$. 
Let $\cal J$ be the linear transformation in ${\mathbb R}^{3,3}$ interchanging $\E_i$ and $\E_i'$ for $i=1,2,3$,
{\it i.e.}, its matrix form is  
\be
{\cal J}:=\left(\ba{cc}
0 &\ \I_{3\times 3} \\

\I_{3\times 3} &\ 0
\ea\right).
\ee
Let $\cal T$ be the linear transformation in ${\mathbb R}^{3,3}$ changing $\E_i$ to $-\E_i$ while preserving
$\E_i'$ for $i=1,2,3$, {\it i.e.}, its matrix form is 
\be
\ba{lll}
{\cal T}:=\left(\ba{cc}
-\I_3 & 0\\

0 & \I_3
\ea\right).
\ea
\ee
Notice that ${\cal T}$ interchanges positive vectors and negative vectors of $\RS$.
The following is obvious:
\be
{\cal J}{\cal T}=-{\cal T}{\cal J}=\left(\ba{cc}
0 & \I_3\\

-\I_3 & 0
\ea\right).
\ee

The group of non-singular linear transformations $\B$ in ${\mathbb R}^{3,3}$ satisfying
\be
\B^T{\cal J}\B=\pm{\cal J}
\ee
is called the {\it group of linear regularities} in ${\mathbb R}^{3,3}$, denoted by
$RL(3,3)$. When the sign is positive, the corresponding subset forms the group of orthogonal transformations
$O(3,3)$; when the sign is negative, the corresponding subset is called the {\it anti-orthogonal transformations},
denoted by $AO(3,3)$. 

The subgroup of special orthogonal transformations $SO(3,3)$ has two connected components:
the component containing the identity transformation $\I_{6\times 6}$ is denoted
by $SO_0(3,3)$, while the component containing $-\I_{6\times 6}$ is denoted by $SO_1(3,3)$.

${\cal J}$ is an orthogonal transformations of determinant $-1$.
The subset of orthogonal transformations of determinant $-1$ is denoted by $SO^-(3,3)$. It also
has two connected components: the component containing matrix $\cal J$ is denoted
by $SO_0^-(3,3)$, while the component containing $-\cal J$ is denoted by $SO_1^-(3,3)$.

${\cal T}$ is an orthogonal transformations of determinant $-1$.
The set of anti-orthogonal transformations of determinant $-1$ is denoted by $SAO^-(3,3)$. It 
has two connected components: the component containing matrix $\cal T$ is denoted
by $SAO_0^-(3,3)$, while the component containing $-\cal T$ is denoted by $SAO_1^-(3,3)$.

${\cal J}{\cal T}$ is an anti-orthogonal transformations of determinant $1$.
The set of anti-orthogonal transformations of determinant 1 is denoted by $SAO(3,3)$. It 
has two connected components: the component containing matrix ${\cal J}{\cal T}$ is denoted
by $SAO_0(3,3)$, while the component containing $-{\cal J}{\cal T}$ is denoted by $SAO_1(3,3)$.

We have
\be\ba{lll}
RL(3,3) &=& SO_0(3,3)\cup SO_1(3,3)\cup SO_0^-(3,3)\cup SO_1^-(3,3)\\
&& \cup
SAO_0(3,3)\cup SAO_1(3,3)\cup SAO_0^-(3,3)\cup SAO_1^-(3,3).
\ea
\ee

The two groups $UR(4)$ and $RL(3,3)$ are related by the {\it Pl\"ucker transform} and
{\it dual Pl\"ucker transform} defined as follows. Let $\e_0, \e_1, \e_2, \e_3$ be a fixed basis of $\RF$.
Let $\Lambda^3(\RF)$ be the realization space of $(\RF)^*$, whose basis
\be
\ce_0=\e_1\wedge \e_2\wedge \e_3,\ \ \
\ce_1=-\e_0\wedge \e_2\wedge \e_3,\ \ \
\ce_2=-\e_0\wedge \e_3\wedge \e_1,\ \ \
\ce_3=-\e_0\wedge \e_1\wedge \e_2
\ee
satisfy for all positive permutations $ijk$ of 123, the following:
\be\ba{llclll}
\e_0\vee \ce_0 &=& \e_i\vee \ce_i &=& 1,\\
\ce_0\vee \ce_i &=& \e_{jk} &=& \E_i',\\
\ce_i\vee \ce_j &=& \e_{0k} &=& \E_k.
\ea
\label{vee:rules}
\ee
Here ``$\vee$" is the {\it meet product} in the Grassmann-Cayley algebra generated by $\RF$.
The pairing between $\e_i$ and $\ce_j$ is defined by
$\e_p\vee \e_q^*=\delta_{pq}$. Let
(\ref{induced:basis}) be the induced Witt basis of $\RS$. 

\bdf
The {\it Pl\"ucker transform} from $GL(4)$ to $GL(3,3)$ is defined by $\A\in GL(4)\mapsto \wedge^2 \A\in GL(3,3)$, where
\be
(\wedge^2 \A)\e_{ij}=(\A\e_i)\wedge (\A\e_j)\in \RS.
\ee
The {\it dual Pl\"ucker transform} from the set of projective polarities $GD(4)$ to $GL(3,3)$ is defined for any projective polarity
$\D$ as 
\be
(\vee^2 \D)\e_{ij}=(\D\e_i)\vee (\D\e_j)\in \RS.
\ee
\edf

For example, for the affine transformation
\[
\A: \left(\ba{c}
x_0 \\
\x
\ea\right)\in \RF \mapsto
\left(\ba{cc}
1 &\ \, 0\\
\t &\ \, \L
\ea\right)
\left(\ba{c}
x_0 \\
\x
\ea\right)\in \RF,
\]
where $\t\in \RT$ and $\L\in GL(3)$, 
the matrix form of $\wedge^2\A$ with respect to the basis
$\E_1, \E_2, \E_3,
\E_1', \E_2', \E_3'$ 
is
\be
\left(\ba{cc}
\L & \ \ \ 0 \\

\t\times \L & \ \ \ \L^{-T}
\ea\right).
\ee

The following result is direct.

\bl 
For any $\A\in GL(4)$,
\be
\det(\wedge^2 \A)=(\det(\A))^3.
\ee
For any $\D\in GD(4)$,
\be
\det(\vee^2\D)=-(\det(\D))^3.
\ee
Furthermore, let $\A_1, \A_2\in GL(4)$ and $\D_1, \D_2\in GD(4)$, and let ``$\circ$" denote the composition of mappings, 
then
\be\ba{ll}
\wedge^2(\A_1\circ \A_2)=(\wedge^2 \A_1)\circ (\wedge^2 \A_2), &\ 
\vee^2(\A_1\circ \D_2)=(\wedge^2 \A_1)\circ (\vee^2 \D_2), \\

\vee^2(\D_1\circ \A_2)=(\vee^2 \D_1)\circ (\wedge^2 \A_2), &\ 
\wedge^2(\D_1\circ \D_2)=(\vee^2 \D_1)\circ (\vee^2 \D_2).
\ea
\label{general:homo}
\ee
\el

The Pl\"ucker transform maps $SL(4)$ onto $SO_0(3,3)$, and the kernel is
$\pm \I_{4\times 4}$. It also maps $SL^-(4)$ onto $SAO^-_0(3,3)$, such that $\pm \A\in SL^-(4)$ are mapped to the same
image of $SAO^-_0(3,3)$. 
For example, the matrix ${\rm diag}(-1,1,1,1)$ is mapped to ${\cal T}\in SAO^-_0(3,3)$, and the pre-images of $\cal T$ are
$\pm {\rm diag}(-1,1,1,1)$. The branches $SO_1(3,3)$ and $SAO^-_1(3,3)$ have no pre-image in
$SL(4)\cup SL^-(4)$.

Similarly, 
the dual Pl\"ucker transform maps $SP^-(4)$ onto $SO^-_0(3,3)$, 
such that any $\B\in SO^-_0(3,3)$ 
has two pre-images $\pm \A\in SP^-(4)$. For example, 
the mapping $\D: \e_i\mapsto \e_i^*$ for $i=0,1,2,3$ is mapped to $\cal J$,
and the pre-images of $\cal J$ are $\pm \D$. The dual Pl\"ucker transform also maps $SP(4)$ onto
$SAO_0(3,3)$, and the pre-images of any $\B\in SAO_0(3,3)$ are of the form $\pm \A\in SP(4)$. 
The branches $SAO_1(3,3)$ and $SO^-_1(3,3)$ have no pre-image in
$SP(4)\cup SP^-(4)$.

\bp
The Pl\"ucker transform is a double-covering homomorphism from $SL(4)$ to $SO_0(3,3)$.
The Pl\"ucker transform and the dual Pl\"ucker transform provide a double-covering homomorphism 
in the sense of (\ref{general:homo}),
from 
$UR(4)$ to the following subgroup of $RL(3,3)$:
\be
RL_0(3,3):=SO_0(3,3)\cup SO^-_0(3,3)\cup SAO_0(3,3)\cup SAO^-_0(3,3).
\ee
\ep

In the setting of Clifford algebra $\CL(3,3)$, the Clifford product is always denoted by juxtaposition of 
participating elements. Any element of the Pin group $Pin(3,3)$ is generated by invertible vectors of unit magnitude
of ${\mathbb R}^{3,3}$ under the Clifford product. 
Any element of the subgroup $Spin(3,3)$ is the Clifford product
of even number of unit vectors. $Spin^-(3,3)$ is the subset of elements that are the
Clifford product of odd number of unit vectors. $Spin(3,3)$ has two connected components, and
the component contain the identity element is denoted by $Spin_0(3,3)$, the other component is denoted by
$Spin_1(3,3)$.

\bp
Any element of $Spin_0(3,3)$ is the Clifford product of even number of negative vectors and even number of
positive vectors; any element of $Spin_1(3,3)$ is the Clifford product of odd number of negative vectors 
and odd number of positive vectors. In particular, $\pm 1$ are in $Spin_0(3,3)$, while $\pm \I_{3,3}$ are in $Spin_1(3,3)$,
where 
\be
\I_{3,3}:=\E_1\wedge \E_2\wedge \E_3\wedge \E_1'\wedge \E_2'\wedge \E_3'=\E_{1231'2'3'},
\ee
is a pseduscalar of $\Lambda(\RS)$ satisfying $\I_{3,3}^2=1$.
\ep

{\it Proof.}\
Denote the set of positive vectors by $P(3,3)$, and denote the set of negative vectors by
$N(3,3)$.
First we prove that each of $P(3,3), N(3,3)$ is connected. By symmetry, we only consider $P(3,3)$.
Let $\v_1, \v_2$ be two positive vectors, then $\v_1\wedge \v_2$ is one of ${\mathbb R}^{2,0,0},
{\mathbb R}^{1,0,1}, {\mathbb R}^{1,1,0}$. Hence there exists a positive vector $\v_3$ that is orthogonal to both 
$\v_1, \v_2$. On Euclidean plane $\v_1\wedge \v_3$, $\v_1$ and $\v_3$ are connected by positive vectors;
on Euclidean plane $\v_2\wedge \v_3$, $\v_2$ and $\v_3$ are connected by positive vectors. So $\v_1, \v_2$ are connected.

Denote the set of elements that are the Clifford product of either two positive vectors or two negative vectors by 
$P_2(3,3)$, and denote the set of elements that are the Clifford product of either one negative vector and one positive vector,
or one positive vector and one negative vector, by $N_2(3,3)$. 
Next we prove that each of $P_2(3,3), N_2(3,3)$ is connected.
Let $\E_+, \E_-$ be a pair of unit positive vector and unit negative vector in $\RS$. Let 
\be
\v_+=\E_+-\lambda \E_-,\ \ \ 
\v_-=\E_-+\lambda \E_+,
\ee
where $0<\lambda<1$. They are respectively a positive vector and a negative vector.
By the continuity of the Clifford multiplication, we only need prove that 
$\E_-\v_-$ equals the Clifford product of two positive vectors, and $\E_-\E_+$ equals the 
Clifford product of a positive vector and a negative vector. Both are true because
\be\ba{lll}
\E_-\v_- &=& -\{\E_-(\E_-\E_+)\}\{(\E_+\E_-)(\E_-+\lambda \E_+)\}
=(-\E_+)(\E_+-\lambda \E_-),\\

\E_-\E_+ &=& (-\E_+)\E_-.
\ea
\ee

Obviously $\pm 1\in P_2(3,3)$. By $(\E_+\E_-)(\E_+\E_-)=1$ and the continuity of the Clifford multiplication, we 
get that the Clifford product of four invertible vectors, where the number of negative vectors is even, must be 
in the same connected component with $P_2(3,3)$. By $(\E_+\E_-)(\E_+\E_+)=\E_+\E_-$ and the continuity of the Clifford multiplication, 
we get that the Clifford product of four invertible vectors, where the number of negative vectors is odd, must be 
in the same connected component with $N_2(3,3)$. 

Denote by $Spin_0(3,3)$ the connected component of $Spin(3,3)$ containing $P_2(3,3)$, and denote by $Spin_1(3,3)$ the 
connected component containing $N_2(3,3)$. By induction on the number of invertible vector factors 
in the factorization of an element of $Spin(3,3)$, we get that $Spin_0(3,3)$ contains all elements that are the Clifford product
of even number of negative vectors and even number of positive vectors, while $Spin_1(3,3)$ contains all elements that are the Clifford product
of odd number of negative vectors and odd number of positive vectors.
Since $Spin(3,3)=Spin_0(3,3)\cup Spin_1(3,3)$ and $Spin(3,3)$ is known to have two connected components, 
$Spin_0(3,3)$ and $Spin_1(3,3)$ are not the same connected component.
\endproof

\bl
${\cal J}\in SO^-(3,3)$ is double covered by $\pm (\E_1-\E_1')(\E_1-\E_1')(\E_1-\E_1')\in Pin^-(3,3)$, and
$-{\cal J}$ is double covered by $\pm (\E_1+\E_1')(\E_1+\E_1')(\E_1+\E_1')\in Pin^-(3,3)$.
\el

\bp
The set $Spin^-(3,3)$ has two connected components. The component double-covering $\cal J$ is denoted by 
$Spin^-_0(3,3)$, whose elements are each the Clifford product of odd number of negative vectors and even number of
positive vectors. The component double-covering $-{\cal J}$ is denoted by 
$Spin^-_1(3,3)$, whose elements are each the Clifford product of even number of negative vectors and odd number of
positive vectors. In particular, all positive vectors are in $Spin^-_1(3,3)$, while all
negative vectors are in $Spin^-_0(3,3)$.
\ep

The group $O(3,3)$ is double-covered by $Pin(3,3)$, and the covering homomorphism is given as following:
for any $\pm \U\in Pin(3,3)$,
\be
Ad^*_\U \X:=\epsilon\U\X\U^{-1},\ \hbox{ for } \X\in {\mathbb R}^{3,3}
\ee
is a transformation belonging to $O(3,3)$, 
where $\epsilon=1$ if $\U\in Spin(3,3)$, and $\epsilon=-1$ if $\U\in Spin^-(3,3)$. 

Denote 
\be
Pin_0(3,3)=Spin_0(3,3)\cup Spin^-_0(3,3),\ \ \,
Pin_1(3,3)=Spin_1(3,3)\cup Spin^-_1(3,3).
\ee 
$Pin_0(3,3)$ is a subgroup of $Pin(3,3)$, and double covers 
$SO_0(3,3)\cup SO^-_0(3,3)$. In particular,
$SO_0(3,3)$ is double covered by $Spin_0(3,3)$.
Now that the Pl\"ucker transform provides another 
double-covering homomorphism of $SO_0(3,3)$ by $SL(4)$, we get the classical result that the two groups
$SL(4)$ and $Spin_0(3,3)$ are isomorphic.

Below we extend the above double-covering map to $SAO(3,3)\cup SAO^-(3,3)$. 
We have seen that ${\cal T}\in SAO^-(3,3)$. 
For any $\B\in O(3,3)$, obviously
${\cal T}\B{\cal T}\in O(3,3)$, so any element of $AO(3,3)$ must have a unique
matrix form ${\cal T}\B$ for some $\B\in O(3,3)$. In fact, we have the following:

\bdf
Define the following isomorphism in group $Pin(3,3)$:
For any 
$\U=\Y_1\Y_2\cdots \Y_r\in Pin(3,3)$ where $\Y_i\in \RS$, 
\be
\U^{\cal T}:=({\cal T}\Y_1)({\cal T}\Y_2)\cdots ({\cal T}\Y_r).
\ee
\edf

\bl
For any $\U\in Pin(3,3)$ and any $\X\in \RS$,
\be
{\cal T}(Ad^*_\U \X)=Ad^*_{\U^{\cal T}} ({\cal T}\X).
\label{change:T}
\ee
\el

{\it Proof.}
If (\ref{change:T}) is true for $\U=\Y_1$, by induction it is true for any other element of $Pin(3,3)$.
It holds for $\U=\Y_1$ by direct verification. 
\endproof

\bdf
By defining a formal associative product between ${\cal T}$ and $Pin(3,3)$ satisfying the following
commutativity, a new group is generated, denoted by
$Pin^{sp}(3,3)$: for any $\U\in Pin(3,3)$,
\be
{\cal T}\circ \U = \U^{\cal T}\circ {\cal T}.
\label{assoc:gp}
\ee
Let ${\cal T}Pin(3,3)$ be the coset of $Pin(3,3)$ with respect to ${\cal T}$, then
\be
Pin^{sp}(3,3)=Pin(3,3)\cup {\cal T}Pin(3,3).
\ee
The {\it adjoint action} of ${\cal T}Pin(3,3)$ upon $\RS$ is defined as follows: for any $\U\circ {\cal T}\circ \V\in
{\cal T}Pin(3,3)$, where $\U, \V\in Pin(3,3)$, for any $\X\in \RS$,
\be
Ad_{\U\circ {\cal T}\circ \V}^* \X :=Ad_{\U}^*({\cal T}(Ad_{\V}^*\X)).
\ee
\edf

$Pin^{sp}(3,3)$ double covers $O(3,3)\cup AO(3,3)$ by the adjoint action. For $i=0,1$, let 
${\cal T}Pin_i(3,3)$ be the coset of $Pin_i(3,3)$ with respect to ${\cal T}$ for $i=0,1$. 
Then 
\be
Pin^{sp}_0(3,3)=
Pin_0(3,3)\cup {\cal T}Pin_0(3,3)
\ee
double covers $RL_0(3,3)$. Since $UR(4)$ also double covers $RL_0(3,3)$ by the Pl\"ucker transform and dual
Pl\"ucker transform, $UR(4)$ and $Pin^{sp}_0(3,3)$ are isomorphic.

Below we realize $\cal T$ in $\CL(3,3)$.
Let
\be
\K_2 = \E_{11'}+\E_{22'}+\E_{33'}:=\e_{01}\wedge \e_{23}+\e_{02}\wedge \e_{31}+\e_{03}\wedge \e_{12}.
\ee
It is called the {\it symplectic form} of $\RS$ with respect to the Witt decomposition $\RS=\I_3\oplus \J_3$.

\bl
$\K_2$ is invariant under any general linear transformation $\C$ in the 3-space $\I_3$ and the
associated linear transformation $\C^{-T}$ in the 3-space $\J_3$. In other words, it is independent of the
the choice of Witt basis of the fixed Witt decomposition.
\el

{\it Proof.} Let $\e_1, \e_2, \e_3, \e_1', \e_2', \e_3'$ be a Witt basis of $\I_3\oplus \J_3$, and let
$\A\in GL(\I_3)$ such that $\A\e_i=\a_i$. Let $\a_i=(a_{1i},a_{2i},a_{3i})^T$, and let
$\a_i'=(a_{1i}',a_{2i}',a_{3i}')^T$, where $a_{ij}$ is the minor of $\A$ by 
removing the $i$-th row and $j$-th column. Then $\A^{-T}\e_i'=\a_i'/\det(\A)$, 
and $\a_i^T\a_j'=\delta_{ij}\det(\A)$. We have
\[\ba{ll}
& \a_1\wedge \a_1'+\a_1\wedge \a_1'+\a_1\wedge \a_1' \\

=& \phantom{-}\bigstrut
 (a_{11}a_{11'}+a_{21}a_{21'}+a_{31}a_{31'})\e_1\wedge \e_1' 
+(a_{11}a_{12'}+a_{21}a_{22'}+a_{31}a_{32'})\e_1\wedge \e_2' \\

& 
+(a_{11}a_{13'}+a_{21}a_{23'}+a_{31}a_{33'})\e_1\wedge \e_3'
+(a_{12}a_{11'}+a_{22}a_{21'}+a_{32}a_{31'})\e_2\wedge \e_1'\\

&
+(a_{12}a_{12'}+a_{22}a_{22'}+a_{32}a_{32'})\e_2\wedge \e_2'
+(a_{12}a_{13'}+a_{22}a_{23'}+a_{32}a_{33'})\e_2\wedge \e_3'
\\

&
+(a_{13}a_{11'}+a_{23}a_{21'}+a_{33}a_{31'})\e_3\wedge \e_1'
+(a_{13}a_{12'}+a_{23}a_{22'}+a_{33}a_{32'})\e_3\wedge \e_2'\\

&
+(a_{13}a_{13'}+a_{23}a_{23'}+a_{33}a_{33'})\e_3\wedge \e_3'
\\

=& \det(\A) (\e_1\wedge \e_1'+\e_2\wedge \e_2'+\e_3\wedge \e_3').\bigstrut
\ea
\]
\vskip -.8cm
\endproof

Some simple facts about $\K_2$:
\bi
\item
$
\K_2^2=3-2\K_2\I_{3,3}.
$

\item $\K_2$ is invertible:
$
\K_2^{-1}=(\K_2+2\I_{3,3})/3.
$

\item
For any $\X\in {\mathbb R}^{3,3}$,
$
\K_2\X\K_2=(1-2\I_{3,3})\X.
$

\item For any $\X, \Y\in {\mathbb R}^{3,3}$,
if $\X\wedge \K_2=\Y\wedge \K_2$, then $\X=\Y$.
\ei

\bp
For any $\X\in \RS$,
\be
{\cal T}\X=\X\cdot \K_2.
\ee
\ep

We have seen that $UR(4)$ is isomorphic to only half of the group $Pin^{sp}(3,3)$. To study 3-D projective transformations and
polarities, the homogeneous model $\RF$ of 3-D projective geometry does not allows for the whole group $Pin^{sp}(3,3)$ to be used.
To overcome this drawback, we need to proceed to use the null 3-space representation of projective points and planes, instead of
returning to $\RF$. 
This ideal leads to the following new model of 3-D projective geometry.

\bdf
The {\it Pl\"ucker model} of 3-D projective geometry refers to the study of 3-D projective transformations and
polarities by using the adjoint action of the whole group $Pin^{sp}(3,3)$ upon the  
null 3-spaces of $\RS$ as the representation space of 3-D points and planes. 
For any element $\C_3=\X_1\wedge \X_2\wedge\X_3\in \Lambda^3(\RS)$ where $\X_i\in \RS$,
for any $\U\in Pin^{sp}(3,3)$, the adjoint action is the following:
\be
(\wedge^3 Ad^*_\U)\C_3:=(Ad^*_\U \X_1)\wedge (Ad^*_\U \X_2)\wedge (Ad^*_\U \X_3).
\ee
\edf

\bdf
The 3-D {\it projective regularity group} $PR(3)$ is defined as the quotient of $UR(4)$ modulo the equivalence relation
``$*$" in which two elements are equivalent if and only if they differ by scale. Usually we write
\be
PR(3):=UR(4)/*.
\ee
\edf

In the next section, it will be shown that via the Pl\"ucker model, group $Pin^{sp}(3,3)$ quadruple-covers group $PR(3)$,
and the kernel of the covering homomorphism is $\{\pm 1, \pm \I_{3,3}\}$.

\section{Some properties of the Pl\"ucker transform}
\setcounter{equation}{0}

In this section we investigate two topics: null 3-vector representations of 3-D points and planes and their 
covariance under the group action of $GL(4)$ upon $\RF$, and the adjoint
action of $Pin^{sp}(3,3)$ upon the null 3-vector representations.

{\bf Notations of algebra.}
An orthonormal basis $\e_0, \e_1, \e_2, \e_3$ is fixed in $\RF$, which induces a fixed 
Witt basis $\E_i=\e_{0i}$ and $\E_i'=\e_{jk}$ for positive
permutation $ijk$ of 123, where $i=1,2,3$. Denote
\be
\E_{i_1\ldots i_r}=\E_{i_1}\wedge \cdots \wedge \E_{i_r},
\ee
for $i_j\in \{1,2,3,1',2',3'\}$. Denote $\e_{123}=\e_1\e_2\e_3$, and for any $\x\in \RT$, denote
\be
\x^\perp:=\x\cdot \e_{123}.
\ee
For example, $\e_i^\perp=\e_{jk}$ for positive
permutation $ijk$ of 123.

For two 3-vectors $\P_3, \Q_3\in \Lambda(\RF)$, denote by $\P_3\Q_3$ the meet product
$\P_3\vee \Q_3$. This notation does not conflict with the same usage of the juxtaposition 
for representing the outer product of two vectors of $\RF$, as the outer product of two 3-vectors 
in $\Lambda(\RF)$ is always zero, so is the meet product of two vectors.

{\bf Notations of geometry.}
Let $\x\in \RT$, then vector $(x_0, \x)=x_0\e_0+\x\in \RF$ represents an affine point if and only if $x_0\neq 0$.
A plane is determined by the equation $x_0(-d)+\n\cdot \x=0$ for point $(x_0, \x)$ on it, where $\n\in \RT$ and
$d\in \mathbb R$, and at least one of them is nonzero. when $d\n\neq 0$, the equation represents the affine plane 
normal to vector $\n$ and with signed distance $d/|\n|$ from the origin along direction $\n$. 
When $d=0$, the plane passes through the origin; when $\n=0$, the plane is the plane at infinity.
So the pair $(\n, -d)\in \RT\times $ represents the plane, and the representation is unique up to scale.

\bp
Point $(x_0, \x)\in \RT$ has the following null 3-vector form in $\Lambda^3(\RS)$:
\be
(\e_0\x)\wedge ((\e_0\x)\cdot \J_3)-x_0 (\e_0\x)\wedge \K_2+x_0^2 \I_3,
\label{expr:pt}
\ee
and plane $(\n, -d)$ has 3-vector form
\be
\n^\perp\wedge (\n^\perp\cdot \I_3)-d \n^\perp\wedge \K_2 +d^2 \J_3,
\label{expr:plane}
\ee
\ep

{\it Proof.}
When $x_0\e_0+\x$ is an affine point, where $\x=x_1\e_1+x_2\e_2+x_3\e_3$, 
then $x_0\neq 0$, and the following three lines
pass through the point, and form a basis of the 3-space of lines through the point:
$(x_0\e_0+\x)\e_1, (x_0\e_0+\x)\e_2, (x_0\e_0+\x)\e_3$.
We have
\[\ba{ll}
& ((x_0\e_0+\x)\e_1)\wedge ((x_0\e_0+\x)\e_2)\wedge ((x_0\e_0+\x)\e_3) \\

=& (x_0 \e_{01}-x_2\e_{12}+x_3\e_{31})\wedge 
(x_0 \e_{02}+x_1\e_{12}-x_3\e_{23})\wedge
(x_0 \e_{03}-x_1\e_{31}+x_2\e_{23}) \bigstrut\\

=& 
x_0^3 \I_3 -x_0^2\{x_1\e_{01}\wedge (\e_{02}\wedge \e_{31}+\e_{03}\wedge \e_{12})
+x_2\e_{02}\wedge (\e_{01}\wedge \e_{23}+\e_{03}\wedge \e_{12})
 \bigstrut\\
 
& \hfill
+x_3\e_{03}\wedge (\e_{01}\wedge \e_{23}+\e_{02}\wedge \e_{31})
\} \\

& 
+x_0(x_1\e_{01}+x_2\e_{02}+x_3\e_{03})
\wedge \{
(-x_1\e_{12}\wedge \e_{31}+x_2\e_{12}\wedge \e_{23}-x_3\e_{31}\wedge \e_{23}
\}\\

=& x_0^3 \I_3 -x_0^2 (\e_0\x)\wedge \K_2 +x_0 (\e_0\x)\wedge 
\{(\e_0\x)\cdot (\e_{23}\wedge \e_{31}\wedge \e_{12})\}. \bigstrut
\ea\]

Removing the common factor $x_0$, we get (\ref{expr:pt}). When $x_0=0$, then $\e_0\x$ is
the line through $\x$ and the origin, and $(\e_0\x)\cdot \J_3$ are spanned by two lines at infinity
that meet line $\e_0\x$, {\it i.e.}, through the point at infinity $\x$. So 
$(\e_0\x)\wedge ((\e_0\x)\cdot \J_3)$ is the 3-vector representing the point at infinity $\x$.

For an affine plane $(\n, -d)$ where $\n^2=1$,
let $\a, \b, \n$ be an orthonormal frame of $\RT$. 
Then $\n^\perp=\a\b$, and 
\[\ba{lll}
\I_3 &=& (\e_0\n)\wedge (\e_0\a)\wedge (\e_0\b), \\
\J_3 &=& (\n\a)\wedge (\n\b)\wedge (\a\b), \\
\K_2 &=& (\e_0\n)\wedge (\a\b)-(\e_0\a)\wedge (\n\b)+(\e_0\b)\wedge (\n\a).
\ea\]
The following three
lines span the 3-space of lines on the plane: $(\e_0+d\n)\a, (\e_0+d\n)\b, \a\b$. We have
\[\ba{ll}
& ((\e_0+d\n)\a)\wedge ((\e_0+d\n)\b)\wedge (\a\b)\\

=& (\e_0\a)\wedge (\e_0\b)\wedge (\a\b) +d^2 (\n\a)\wedge (\n\b)\wedge (\a\b)\bigstrut  \\

&\hfill
+d (\a\b)\wedge \{(\e_0\a)\wedge (\n\b)-(\e_0\b)\wedge (\n\a)\}
\\

=& (\a\b)\wedge ((\a\b)\cdot \I_3)-d (\a\b)\wedge \K_2+d^2 \J_3.\bigstrut  
\ea\]
Thus we get (\ref{expr:plane}) for affine plane. When $\n=0$, then $\J_3$ is the plane at infinity.
\endproof

\bdf
\label{def:represent:p}
The map
\be
{\cal F}:
x_0\e_0+\x \ \ \mapsto\ \ (\e_0\x)\wedge ((\e_0\x)\cdot \J_3)-x_0 (\e_0\x)\wedge \K_2+x_0^2 \I_3
\label{def:F}
\ee
is called the {\it Pl\"ucker representation} of the point. The plane $(\n, d)$ can be represented by vector
$d\e_0+\n$; for the plane, the map
\be
{\cal F}':
d\e_0+\n \ \ \mapsto\ \ \n^\perp\wedge (\n^\perp\cdot \I_3)+d \n^\perp\wedge \K_2 + d^2\J_3
\label{def:Fp}
\ee
is called the {\it Pl\"ucker representation} of the plane. 
\edf

{\it Remark.}\ The Pl\"ucker representations of points and points at infinity 
$\e_0, \e_1$, $\e_2, \e_3$ are respectively $\E_{123}, \E_{12'3'},
\E_{23'1'}, \E_{31'2'}$. Furthermore,
\be
\ba{lll}
{\cal F}(x_0\e_0+\x) &=& {\cal F}(x_0\e_0)+{\cal F}(\x)-((x_0\e_0)\x)\wedge \K_2, \\

{\cal F}'(d\e_0+\n) &=& {\cal F}'(d\e_0)+{\cal F}'(\n)+(d \n^\perp)\wedge \K_2.\bigstrut
\ea\ee

Definition \ref{def:represent:p} relies upon the specific Witt decomposition $\I_3\oplus \J_3$. The following property gives
an alternative definition of $\cal F$ that is independent of the Witt decomposition.

\bp
\label{prop:fund:ff}
Let $\X, \a_1, \a_2, \a_3$ be a basis of $\RT$. Then
\be
{\cal F}(\X)=\frac{(\X\a_1)\wedge (\X\a_2)\wedge (\X\a_3)}{[\X\a_1\a_2\a_3]}.
\label{good:F}
\ee
In particular, ${\cal F}(\X)$ is independent of the choice of $\a_1, \a_2, \a_3$ as long as the latter
forms a basis of the 4-space together with $\X$.
Furthermore, 
\be
{\cal F}'=(\wedge^3 {\cal J})\circ {\cal F}. \label{fpf}
\ee
\ep

{\it Proof.} (\ref{fpf}) is obvious by (\ref{def:F}) and (\ref{def:Fp}). We only need to prove (\ref{good:F}).
When $\X$ is an affine point, by definition,
\[
{\cal F}(\X)=\frac{(\X\e_1)\wedge (\X\e_2)\wedge (\X\e_3)}{[\X\e_1\e_2\e_3]}.
\]
For $i=1,2,3$, let
$
\a_i=d_i\X+c_{1i}\e_1+c_{2i}\e_2+c_{3i}\e_3.
$
Then $\det(c_{ji})_{i,j=1..3}\neq 0$, as
\[
(\X\ \, \a_1\ \a_2\ \a_3)=(\X\ \, \e_1\ \e_2\ \e_3)
\left(\ba{ll}
1 &\ \ \d^T \\
0 &\ \ \C
\ea\right),
\]
where $\d^T=(d_1,d_2,d_3)$, and $\C=(c_{ji})_{i,j=1..3}$.
By direct verification,
\[
\frac{(\X\a_1)\wedge (\X\a_2)\wedge (\X\a_3)}{[\X\a_1\a_2\a_3]}
=\frac{\det(\C)(\X\e_1)\wedge (\X\e_2)\wedge (\X\e_3)}{\det(\C)[\X\e_1\e_2\e_3]}.
\]

When $\X$ is a point at infinity, let $\X=\lambda \x$, where $\x$ is a unit vector, 
and let $\x, \a, \b$ be an orthonormal basis of 
the 3-space spanned by $\e_1, \e_2, \e_3$,
then $[\X\e_0\a\b]=-\lambda$.
By definition,
\[
{\cal F}(\X) = \lambda^2 (\e_0\x)\wedge ((\e_0\x)\cdot (\a\b\wedge \b\x\wedge \x\a))
=
\frac{(\X\e_0)\wedge (\X\b)\wedge (\X\a)}{[\X\e_0\b\a]}.
\]
\vskip -.8cm
\endproof

The following proposition answers the following question: Given ${\cal F}(\X)$ for $\X=x_0\e_0+\x$, how to recover $\X$?

\bp
For any $\X=x_0\e_0+\x\in \RF$, and any $\Pi=d\ce_0+\cn\in (\RF)^*$, 
\be\ba{lll}
{\cal F}(\X)\cdot \K_2 &=& -2x_0(\e_0\x),\\
{\cal F}''(\Pi)\cdot \K_2 &=& -2d\,(\ce_0\cn).
\ea
\ee
\ep

{\it Proof.} Let $\x=x_1\e_1+x_2\e_2+x_3\e_3$, then
\[
{\cal F}(\X)\cdot \K_2 =
\sum_{i=1}^3 x_i (\E_i\wedge (\e_0\x))\cdot \J_3 -3x_0(\e_0\x)+x_0\sum_{i=1}^3 x_i \E_i 
= -2x_0(\e_0\x).
\]
The second result can be proved similarly. 
\endproof

The following is on the covariance of ${\cal F}(\X)$: if $\X\in \RF$ undergoes a general linear transformation, 
how does ${\cal F}(\X)$ change in $\Lambda^3(\RS)$ accordingly?

\bp
\label{prop:inv}
Let $\A$ be a non-singular linear transformation in the 4-space spanned by $\e_0, \e_1, \e_2, \e_3$.
Then for any vector $\X=x_0\e_0+\x$ of the 4-space, 
\be
\det(\A) {\cal F}(\A\X)=\wedge^3(\wedge^2\A) {\cal F}(\X),
\label{relation:p1}
\ee
and
\be
\det(\A) {\cal F}'(\A^{-T}\X)=\wedge^3({\cal J}\circ (\wedge^2\A^{-T})\circ {\cal J}) {\cal F}'(\X).
\label{relation:p2}
\ee
\ep

{\it Proof.}
Without loss of generality, let $\X$ represent an affine point in space, then 
\[\ba{lll}
\wedge^3(\wedge^2\A) {\cal F}(\X) 
&=& \ds\frac{(\wedge^2\A)(\X\e_1)\wedge (\wedge^2\A)(\X\e_2)\wedge (\wedge^2\A)(\X\e_3)}{[\X\e_1\e_2\e_3]}\\

&=& \ds \det(\A)\frac{(\A\X)(\A\e_1)\wedge (\A\X)(\A\e_2)\wedge (\A\X)(\A\e_3)}
{[(\A\X)(\A\e_1)(\A\e_2)(\A\e_3)]} \Bigstrut\\

&=& \det(\A){\cal F}(\A\X). \bigstrut
\ea\]
(\ref{relation:p2}) is obvious from (\ref{relation:p1}), where $\A^{-T}$ is the induced
linear transformation in the vector space of planes.
\endproof

Alternatively, a plane $(\n, d)$ can be represented by a 3-vector in $\Lambda({\mathbb R}^4)$.
Let $\n=a\e_1+b\e_2+c\e_3$, then $\Pi=d\ce_0+a\ce_1+b\ce_2+c\ce_3$ represents the plane.
A point $\X=x_0\e_0+x_1\e_1+x_2\e_2+x_3\e_3$ is on the plane if and only if
$\X\vee \Pi=0$, {\it i.e.}, $x_0d+x_1a+x_2b+x_3c=0$.

\bp
Let $\A\in GL(4)$ act upon $\RF$. 
In $\Lambda^3(\RF)$ and $\Lambda^2(\RS)$, 
the two matrices $\wedge^3\A=\det(\A)\A^{-T}$ and $\B=\vee^2 (\wedge^3\A)$ satisfy
\be
\det(\B)=(\det(\det(\A)\A^{-T}))^3=\det(\A)^{9}.
\ee
\ep

Let $\cn=a\ce_1+b\ce_2+c\ce_3$.
When $\Pi=d\ce_0+\cn$ is an affine plane not through the origin,  
then $d\cn\neq 0$, and the plane meets the following three planes
passing through the origin: $\ce_1, \ce_2, \ce_3$, so lines $\Pi\ce_i$ for $i=1,2,3$ 
form a basis of the 3-space of lines on plane $\Pi$. By (\ref{vee:rules}),
\[\ba{ll}
& (\Pi\ce_1)\wedge (\Pi\ce_2)\wedge (\Pi\ce_3) \\

=& (d \e_{23}-b\e_{03}+c\e_{02})\wedge 
(d \e_{31}+a\e_{03}-c\e_{01})\wedge
(d \e_{12}-a\e_{02}+b\e_{01}) \bigstrut\\

=& 
d^3 \J_3 +d^2\{a\e_{23}\wedge (\e_{02}\wedge \e_{31}+\e_{03}\wedge \e_{12})
+b\e_{31}\wedge (\e_{01}\wedge \e_{23}+\e_{03}\wedge \e_{12})
 \bigstrut\\
 
& \hfill
+c\e_{12}\wedge (\e_{01}\wedge \e_{23}+\e_{02}\wedge \e_{31})
\} \\

& 
+d(a\e_{23}+b\e_{31}+c\e_{12})
\wedge \{
(a\e_{02}\wedge \e_{03}+b\e_{03}\wedge \e_{01}+c\e_{01}\wedge \e_{02}
\}\\

=& d^3 \J_3 +d^2 (\ce_0\Pi)\wedge \K_2 +d (\ce_0\Pi)\wedge 
\{(\ce_0\Pi)\cdot \I_3\}. \bigstrut
\ea\]

Removing the common factor $d$, we get the following alternative form of (\ref{expr:plane}):
\be
(\ce_0\Pi)\wedge ((\ce_0\Pi)\cdot \I_3)+d (\ce_0\Pi)\wedge \K_2 +d^2 \J_3.
\label{expr:plane:2}
\ee
When $\Pi$ passes through the origin, then $d=0$, and (\ref{expr:plane:2}) still represents the
plane. When $\Pi$ is the plane at infinity, then $\Pi=d\ce_0$, so $\ce_0\Pi=0$,
and (\ref{expr:plane:2}) once again represents the plane.

\bdf
When a plane $(\n, d)$ is represented by a 3-vector
$\Pi=d\ce_0+\cn$, the map
\be
{\cal F}'':
d\ce_0+\cn \ \ \mapsto\ \ (\ce_0\Pi)\wedge ((\ce_0\Pi)\cdot \I_3)+d (\ce_0\Pi)\wedge \K_2 + d^2\J_3
\ee
is called the {\it dual Pl\"ucker representation} of the plane. 
\edf

{\it Remark.}\
\be
{\cal F}''(d\ce_0+\cn) = {\cal F}''(d\ce_0)+{\cal F}''(\cn)+((d\ce_0)\cn)\wedge \K_2.
\ee

\bl
For the plane $(\n, d)$ represented by 3-vector $\Pi=d\ce_0+\cn$,
\be
d=\Pi\ce_1\ce_2\ce_3=\Pi\vee \ce_1\vee \ce_2\vee \ce_3.
\ee
\el

\bl Let $\A\in GL(4)$ act upon $\RF$. Then
in ${\mathbb R}^{3,3}$,
\be
\wedge^2\A= \det(\A)\, (\vee^2\hskip -.05cm \A^{-T}).
\ee
\el

{\it Proof.} Let $\A=(\a_0\ \a_1 \ \a_2\ \a_3)$.
On one hand,
$\wedge^2\A (\e_{0i})=\A(\e_0)\wedge \A(\e_i)=\a_0\wedge \a_i$, and
$\wedge^2\A (\e_{ij})=\a_i\wedge \a_j$. On the other hand, for positive permutation $ijk$ of 123,
\[\ba{lll}
\vee^2 \A^{-T}(\e_{0i}) &=&  \vee^2 (\det(\A)^{-1} \wedge^3\hskip -.1cm \A) (\ce_j\vee \ce_k)\\

&=& \det(\A)^{-2} (\a_0\wedge \a_k\wedge \a_i)\vee (\a_0\wedge \a_i\wedge \a_j)\\

&=& \det(\A)^{-1} \a_0\wedge \a_i; \\
\\

\vee^2 \A^{-T}(\e_{ij}) &=&  \vee^2 (\det(\A)^{-1} \wedge^3\hskip -.1cm \A) (\ce_0\vee \ce_k)\\

&=& -\det(\A)^{-2} (\a_i\wedge \a_j\wedge \a_k)\vee (\a_0\wedge \a_i\wedge \a_j)\\

&=& \det(\A)^{-1} \a_i\wedge \a_j.
\ea
\]
\endproofs

\bp
For $\X\in \RF$ and $\Pi\in \Lambda^3(\RF)$,
\be
{\cal F}(\X)\cdot {\cal F}''(\Pi)=-[\X\Pi]^2.
\ee
\ep

{\it Proof.} Let $\X=x_0\e_0+\x$ and $\Pi=d\ce_0+\cn$, then
\[\ba{lll}
{\cal F}(\X)\cdot {\cal F}''(\Pi) &=& ({\cal F}(\X)+x_0^2 {\cal F}(\e_0)-x_0(\e_0\x)\wedge \K_2)\\
&&\hskip .5cm
\cdot ({\cal F}''(\cn)+d^2 {\cal F}''(\ce_0)+d(\ce_0\cn)\wedge \K_2)\bigstrut\\

&=& \bigstrut
-(\x\cdot \n)^2-x_0^2d^2-2x_0d\,(\x\cdot \n) \\

&=& \bigstrut
-[\X\Pi]^2.
\ea
\]
\endproofs

The proof of the following proposition is similar to those of Proposition \ref{prop:fund:ff}
and Proposition \ref{prop:inv}, and is omitted.

\bp
Let $\Pi, \ca_1, \ca_2, \ca_3$ be a basis of $(\RF)^*$. Then
\be
{\cal F}''(\Pi)=\frac{(\Pi\ca_1)\wedge (\Pi\ca_2)\wedge (\Pi\ca_3)}{\Pi\ca_1\ca_2\ca_3}.
\label{good:F2}
\ee
In particular, ${\cal F}''(\Pi)$ is independent of the choice of $\ca_1, \ca_2, \ca_3$ as long as the latter
forms a basis of the 4-space together with $\Pi$.

Let $\A$ be a non-singular linear transformation in $\RF$. It
induces a linear transformation $\A^{-T}=det(\A)^{-1} \wedge^3\hskip -.1cm\A$
in $(\RF)^*$.
For any $\Pi=d\ce_0+\cn\in (\RF)^*$, 
\be
\det(\A^{-T}) {\cal F}''(\A^{-T}\Pi)=\wedge^3(\vee^2\A^{-T}) {\cal F}''(\Pi).
\label{relation:ff2}
\ee
\ep

Below we extend the previous results from $GL(4)$ to $GD(4)$.
The following is a simple fact from linear algebra on 
the influence of coordinate transformations in $\RF$ upon the matrix form of
$\A: \RF\longrightarrow (\RF)^*$. 

\bp \label{lem:ctrans}
Let $\C=(\c_0\ \c_1\ \c_2\ \c_3)\in GL(4)$ be a coordinate transformation in $\RF$,
and let $\A$ be an invertible linear mapping from $\RF$ to $(\RF)^*$,
whose matrix form with respect to the basis $\e_i$'s of $\RF$
and $\e_j^*$'s of $(\RF)^*$ is $\A_{4\times 4}$. Then 
with respect to the basis $\c_i$'s of $\RF$ and its dual basis $\c_j^*$'s of $(\RF)^*$, where
$\c_j^*=(\wedge^3\C)\e_j^*/\det(\C)$, 
the matrix form of $\A$ is $\C^T\A_{4\times 4}\C$.
\ep

{\it Proof.}\
The coordinate transformation $\C$ in $\RF$ induces the following coordinate transformation in
$\Lambda^3(\RF)$:
\be
(\wedge^3\C)\e_i^*=\sum_{j=0}^3 \det(\C)(\C^{-T})_{ji} \e_j^*.
\ee
Let $\A\e_i=\sum_{j=0}^3 a_{ij}\e_j^*$, and let $\c_i=(c_{0i},c_{1i},c_{2i},c_{3i})^T$. Then
\[\ba{lll}
\ds\A\c_i = \sum_{k=0}^3 c_{ki}\A\e_k 
=\sum_{j,k=0}^3 c_{ki}a_{jk}\e_j^* 
&=& \ds \det(\C)^{-1} \sum_{j,k,l=0}^3 c_{ki}a_{jk} (\C^T)_{lj} (\wedge^3\C)\e_l^*\\

&=& \ds \sum_{l=0}^3 (\C^T\A\C)_{li}\c_l^*.
\ea
\]
\vskip -.8cm
\endproof

{\bf Notation.}\
For any $\A_{4\times 4}=(a_{ij})_{ij,=0..3}$, denote its $2\times 2$ minor
composed of the $u,v$ rows and $p,q$ columns by $\A^{uv}_{pq}$, {\it i.e.},
\be
\A^{uv}_{pq}:=\left|\ba{ll}
a_{up} & a_{uq} \\
a_{vp} & a_{vq}
\ea
\right|.
\ee

Let $\A\in GD(4)$, {\it i.e.}, a non-singular linear map from $\RF$ to $(\RF)^*$. 
Let the matrix form of $\A$ 
with respect to the basis $\e_0, \e_1, \e_2, \e_3$ of $\RF$ and the dual basis
$\ce_0, \ce_1, \ce_2, \ce_3$ of $(\RF)^*$ be
$\A_{4\times 4}=(\a_0\  \a_1\  \a_2\  \a_3)$, where $\a_i=(a_{0i}, a_{1i}, a_{2i}, a_{3i})^T$. 
Then $\A\e_i=\a_i^*:=\sum_{j=0}^3 a_{ji}\e_j^*$, and
\be
(\vee^2\A)(\e_{ij})=\a_i^*\a_j^*=\sum\limits_{0\leq p<q\leq 3} \A^{pq}_{ij} \e_{rs},
\label{vee:matrix1}
\ee
where $pqrs$ is a positive permutation of 0123.

\bl
For the $4\times 4$ identity matrix $\I_{4\times 4}$ representing an element of $GD(4)$, 
$\vee^2\I_{4\times 4}=\cal J$. 
\el

\bp
Let $\A\in GD(4)$ have matrix form $\A_{4\times 4}$, and Let $\B\in GL(3,3)$ be the matrix form of
$\vee^2\A$. When matrix $\A_{4\times 4}$ is taken as a linear transformation of $\RF$, then
$\wedge^2\A_{4\times 4}$ has the matrix form ${\cal J}\B\in GL(3,3)$. 
\ep

{\it Proof.}
Let $\tilde{\A}$ be the non-singular linear transformation in $\RF$ having the same matrix form 
$\A_{4\times 4}=(\a_0\  \a_1\  \a_2\  \a_3)$ with $\A$. Then
the matrix form $\tilde{\B}$
of $\wedge^2 \tilde{\A}$ is given by
\be
(\wedge^2\tilde{\A})(\e_{ij})=\sum\limits_{0\leq p<q\leq 3} \A^{pq}_{ij} \e_{pq}.
\label{matrix:tilde}
\ee
Obviously $\tilde{\B}={\cal J}\B$. 
\endproof

\bp
\label{prop:inv:dual}
Let $\A\in GD(4)$.
Then for any $\X\in \RF$, 
\be
\det(\A) {\cal F}''(\A\X)=\wedge^3(\vee^2\A) {\cal F}(\X).
\label{relation:p1:2}
\ee
\ep

{\it Proof}. We have
\[
\wedge^3(\vee^2\A) {\cal F}(\X)=\frac{(\A\X)(\A\e_1)\wedge (\A\X)(\A\e_2)\wedge (\A\X)(\A\e_3)}{[\X\e_1\e_2\e_3]}.
\]
Since 
$
(\A\e_0)(\A\e_1)(\A\e_2)(\A\e_3)=\a_0^*\a_1^*\a_2^*\a_3^*=\det(\A)\e_0^*\e_1^*\e_2^*\e_3^*=\det(\A),
$
we have $\det(\A)[\X\e_1\e_2\e_3]=(\A\X)(\A\e_1)(\A\e_2)(\A\e_3)$. On the other hand,
\[
{\cal F}''(\A\X)=
\frac{(\A\X)(\A\e_1)\wedge (\A\X)(\A\e_2)\wedge (\A\X)(\A\e_3)}{(\A\X)(\A\e_1)(\A\e_2)(\A\e_3)}.
\]
\endproofs

For any $\U\in Pin^{sp}(3,3)$, $\wedge^3 Ad^*_\U$ maps a null vector to a null vector. If $\U\in Spin(3,3)\cup 
{\cal T}Spin(3,3)$, then $\wedge^3 Ad^*_\U$ maps the set of null 3-spaces representing points to the same set; 
if $\U\in Spin^-(3,3)\cup {\cal T}Spin^-(3,3)$, then $\wedge^3 Ad^*_\U$ interchanges the set of null 3-spaces representing points and
the set of null 3-spaces representing planes. 

\bl
For any $\X\in \RF$ and $\Pi\in (\RF)^*$, 
\be
{\cal F}(\X)\cdot \I_{3,3}= {\cal F}(\X), 
\hskip 1cm
{\cal F}''(\Pi)\cdot \I_{3,3}= -{\cal F}''(\Pi).
\label{3:different}
\ee
Furthermore, for any $\P\in \RS$,\ 
$
Ad_{\I_{3,3}}^*\P=-\P.
$
\el

Given a null 3-vector, (\ref{3:different}) can be used to determine whether it represents a point or plane of ${\cal E}^3$.

{\bf Notation.} For two algebraic elements $\a, \b$, if they differ by scale, we write
\be
\a\seq \b.
\ee

In $\Lambda^3(\RS)$, the representation of a point or plane by a null 3-vector is unique up to scale.
So $Ad_{\I_{3,3}}^*$ leaves each point and plane of ${\cal E}^3$ invariant, and induces the projective transformation represented by the
the identity transformation $\I_{4\times 4}$ of $\RF$. As a consequence, the adjoint action of 
$Pin^{sp}(3,3)$ upon the null 3-vectors of $\Lambda^3(\RS)$ induces a quadruple-covering homomorphism upon the group
$PR(3)$, and the kernel is composed of 4 elements: $\pm 1, \, \pm \I_{3,3}$. 

\bdf
For $\U\in Pin^{sp}(3,3)$, 
The {\it $4\times 4$ matrix representation} of $Ad_{\U}^*$ refers to a matrix form of the linear map induced by $Ad_{\U}^*$ via the
null 3-vector representation of points and planes, from $\RF$ to $\RF$ or $(\RF)^*$. The matrix form is unique up to scale.
\edf

As $Pin^{sp}(3,3)$ double covers $RL(3,3)$, any element $\B\in RL(3,3)$ also has $4\times 4$ matrix representation.
By the covariance of $\cal F$ and ${\cal F}''$, the $4\times 4$ matrix representation after rescaling
is mapped to one of $\pm \B$ by the Pl\"ucker transform or dual Pl\"ucker transform. 
We study how to compute the $4\times 4$ matrix form $\A$ of $\B$. 
Without loss of generality, let $\B\in SO_0(3,3)$. Then $\B$
is the image of $\A\in SL(4)$ under the Pl\"ucker transform.

{\bf Notations used only in this section.} 
Denote
\be
P(0):=123,\ \ \ P(1):=12'3',\ \ \ P(2):=23'1',\ \ \ P(3):=31'2'. 
\ee
Let $ijk$ be a positive permutation of 123, denote
\be
Q(i):=j'k', \ \ \
Q(ij):=k'.
\ee

\bp
Let $\A=(\a_1\ \a_2\ \a_3\ \a_4)=
(a_{pq})_{p,q=0..3}\in SL(4)$, and 
let $\B=(b_{\alpha\beta})_{\alpha,\beta=1,2,3,1',2',3'}$ be the matrix of $\wedge^2 \A$ 
with respect to the basis $\E_1, \E_2, \E_3$, $\E_{1'}, \E_{2'}, \E_{3'}$, {\it i.e.},
for $1\leq i,j\leq 3$, let $ipq$ and $juv$ be both positive permutations of 123, then
\be
b_{ij}=\A^{0i}_{0j}, \ \ \
b_{i'j}=\A^{pq}_{0j},\ \ \
b_{ij'}=\A^{0i}_{uv},\ \ \
b_{i'j'}=\A^{pq}_{uv}.
\ee
Denote by $\B^{uvw}_{pqr}$
the $3\times 3$ minor formed by the $u,v,w$ rows and $p,q,r$ columns of $\B$. Then for any $0\leq i,j\leq 3$,
\be
\det(\A) a_{ij}^2 = \B^{P(i)}_{P(j)}.\label{expr:square}
\ee
\ep

{\it Proof.} Without loss of generality, consider $j=1$.
The following diagram commutes:
\be
\ba{rcl}
\e_1 & \stackrel{\cal F}{\longrightarrow} & \E_{12'3'}\\
|\, && \,| \\
\A|\,& & \,|\wedge^3\hskip -.08cm \B \\
 \downarrow && \downarrow \\
\a_1 & \stackrel{\det(\A) \cal F}{\longrightarrow} & \C_3
\ea
\label{diagram:p1}
\ee
where
\be\ba{lll}
\C_3 &=& \det(\A)\{a_{01}^2 \I_3-a_{01} (\e_0\a_1)\wedge \K_2 + (\e_0\a_1)\wedge ((\e_0\a_1)\cdot \J_3)\}\\

&=& \det(\A)\{
a_{01}^2\E_{123}
+a_{01}a_{11}(\E_{313'}-\E_{122'})
+a_{01}a_{21}(\E_{121'}-\E_{233'}) \Bigstrut\\

&& 
+a_{01}a_{31}(\E_{232'}-\E_{311'})
+a_{11}^2\E_{12'3'}+a_{21}^2\E_{23'1'}+a_{31}^2\E_{31'2'}\\

&& \hfill
+a_{11}a_{21}(\E_{13'1'}+\E_{22'3'})
+a_{11}a_{31}(\E_{11'2'}+\E_{32'3'})
+a_{21}a_{31}(\E_{21'2'}+\E_{33'1'})
\}.
\ea
\label{expr:C3}
\ee

Since $\wedge^3 \B (\E_{12'3'})=\sum_{{ijk}} \B_{12'3'}^{ijk}\E_{ijk}$,
by comparing the coefficients of $\E_{123}$, $\E_{12'3'}, \E_{23'1'}, \E_{31'2'}$ 
in this expression and in $\C_3$, we get 
\[\ba{ll}
\B_{12'3'}^{123}=\det(\A) a_{01}^2, &\hskip .4cm
\B_{12'3'}^{12'3'}=\det(\A) a_{11}^2,\\

\B_{12'3'}^{23'1'}=\det(\A) a_{21}^2, &\hskip .4cm
\B_{12'3'}^{31'2'}=\det(\A) a_{31}^2.\bigstrut
\ea\]
\endproofs

{\it Remark.}
By (\ref{expr:square}), a matrix $\B\in SO(3,3)$ is in $SO_0(3,3)$ if and only if
$\B^{P(i)}_{P(j)}\geq 0$ for all $0\leq i,j\leq 3$.

\bp
Let $ijk$ be a positive permutation of 123, and let $l$ be one of 0,1,2,3, then
\be\ba{lllll}
\det(\A) a_{0l}a_{il} &=& -\B_{P(l)}^{ijj'} &=& -\B_{P(l)}^{ikk'}, \\
\det(\A) a_{il}a_{jl} &=& \phantom{-}\B_{P(l)}^{iQ(j)} &=& \phantom{-}\B_{P(l)}^{Q(i)j}; \bigstrut\\

\det(\A) a_{l0}a_{li} &=& -\B^{P(l)}_{ijj'} &=& -\B^{P(l)}_{ikk'}, \bigstrut\\
\det(\A) a_{li}a_{lj} &=& \phantom{-}\B^{P(l)}_{iQ(j)} &=& \phantom{-}\B^{P(l)}_{Q(i)j}.\bigstrut
\ea
\label{sym:10:2}
\ee
\ep

{\it Proof.} Without loss of generality, consider $l=1$.
In the diagram (\ref{diagram:p1}), by comparing the coefficients of
$\E_{122'}, \E_{313'}, \E_{121'}, \E_{233'}, \E_{311'}, \E_{232'}$ 
in $\sum \B_{12'3'}^{ijk}\E_{ijk}$ and $\C_3$ of (\ref{expr:C3}), we get
the first line of (\ref{sym:10:2}) for $l=1$. By comparing the coefficients of
$\E_{13'1'}, \E_{22'3'}, \E_{11'2'}, \E_{32'3'}, \E_{21'2'}, \E_{33'1'}$, 
we get the second line of (\ref{sym:10:2}) for $l=1$.
When replacing $\A$ with $\A^T$, by the trivial fact
\be
\wedge^2 (\A^T)= (\wedge^2 \A)^T,
\ee
we get the last two lines of (\ref{sym:10:2}). \endproof

\bp
Let $ijk$ and $pqr$ be both positive permutations of 123, then
\be\ba{lllll}
\det(\A) (a_{0p}a_{iq}+a_{0q}a_{ip}) &=& -\B^{ijj'}_{pQ(q)}-\B^{ijj'}_{Q(p)q}
&=& -\B^{ikk'}_{pQ(q)}-\B^{ikk'}_{Q(p)q}, \\

\det(\A) (a_{ip}a_{jq}+a_{iq}a_{jp}) &=& \phantom{-}\B^{iQ(j)}_{pQ(q)}+\B^{iQ(j)}_{Q(p)q}
&=& \phantom{-}\B^{Q(i)j}_{pQ(q)}+\B^{Q(i)j}_{Q(p)q}.
\Bigstrut
\ea
\label{expr:group3}
\ee
\ep

{\it Proof.}
Without loss of generality, consider $i=1$ and $j=2$.
The following diagram commutes:
\be
\ba{rcl}
\e_1+\e_2 & \stackrel{\cal F}{\longrightarrow} & \E_{12'3'}+\E_{23'1'}+\E_{13'1'}+\E_{22'3'}\\
|\, && \,| \\
\A|\,& & \,|\wedge^3\hskip -.08cm \B \\
 \downarrow && \downarrow \\
\a_1+\a_2 & \stackrel{\det(\A) \cal F}{\longrightarrow} & \D_3
\ea
\label{diagram:p2}
\ee
where for $\C_3$ of (\ref{expr:C3}),
\be\ba{lll}
\D_3 
&=& \C_3+\det(\A)\{
2a_{01}a_{02}\E_{123}
+(a_{01}a_{12}+a_{02}a_{11})(\E_{313'}-\E_{122'})\Bigstrut\\

&& 
+(a_{01}a_{22}+a_{02}a_{21})(\E_{121'}-\E_{233'}) 
+(a_{01}a_{32}+a_{02}a_{31})(\E_{232'}-\E_{311'})\\

&&
+2a_{11}a_{12}\E_{12'3'}+2a_{21}a_{22}\E_{23'1'}+2a_{31}a_{32}\E_{31'2'}\\

&& \hfill
+(a_{11}a_{22}+a_{12}a_{21})(\E_{13'1'}+\E_{22'3'})
+(a_{11}a_{32}+a_{12}a_{31})(\E_{11'2'}+\E_{32'3'})\\

&& \hfill
+(a_{21}a_{32}+a_{22}a_{31})(\E_{21'2'}+\E_{33'1'})
\}.
\ea
\label{expr:D3}
\ee

Since $\wedge^3\B (\E_{13'1'}+\E_{22'3'})=\sum_{{ijk}} (\B_{13'1'}^{ijk}+\B_{22'3'}^{ijk})\E_{ijk}$,
by comparing the coefficients of 
$\E_{313'}, \E_{122'}, \E_{121'}, \E_{233'}, \E_{232'}, \E_{311'},
\E_{13'1'}, \E_{22'3'}, \E_{11'2'}, \E_{32'3'}$, \\
$\E_{21'2'}, \E_{33'1'}$ 
in this expression and in $\D_3-\C_3$, we get for positive permutation $pqr$ of 123 the following:
\[\ba{rrrrr}
\det(\A) (a_{01}a_{p2}+a_{02}a_{p1}) &=& -\B^{pqq'}_{13'1'}-\B^{pqq'}_{22'3'}
&=& -\B^{prr'}_{13'1'}-\B^{prr'}_{22'3'}, \\

\det(\A) (a_{p1}a_{q2}+a_{q2}a_{p1}) &=& \B^{pr'p'}_{13'1'}+\B^{pr'p'}_{22'3'}
&=& \B^{q'r'q}_{13'1'}+\B^{q'r'q}_{22'3'}.
\Bigstrut
\ea
\]
\endproofs

The $4\times 4$ matrix representation $\A_{4\times 4}=(a_{pq})_{p,q=0..3}$ of $\B\in SO_0(3,3)$ is given explicitly by 
(\ref{expr:square}), (\ref{sym:10:2}), and (\ref{expr:group3}) as following: one entry of
the first row of $\A$ must be nonzero, say 
$a_{0q_0}\neq 0$ for some $0\leq q_0\leq 3$. From (\ref{expr:square}) we get $a_{0q_0}^2$; from
$a_{0q_0}$ and (\ref{sym:10:2}) we get $a_{iq_0}$ for all $1\leq i\leq 3$, and $a_{0q}$ for all $q\neq q_0$.
From $a_{0q_0}a_{iq}+a_{iq_0}a_{0q}$ of (\ref{expr:group3}) we get the rest components $a_{iq}$ of $\A$.
The last step can be made more explicit as follows.

By
\[\ba{lllll}
a_{0p}a_{iq}-a_{0q}a_{ip} &=& \A^{0i}_{pq}&=& b_{ir'}, \\

a_{ip}a_{jq}-a_{iq}a_{jp} &=& \A^{ij}_{pq}&=& b_{k'r'},\bigstrut
\ea
\]
we get

\bp
\be\ba{lll}
2a_{0p}a_{iq} &=& -\B^{ijj'}_{pQ(q)}-\B^{ijj'}_{Q(p)q}+b_{iQ(pq)}\det(\A)\\

&=& -\B^{ikk'}_{pQ(q)}-\B^{ikk'}_{Q(p)q}+b_{iQ(pq)}\det(\A), \bigstrut\\

2a_{0q}a_{ip} &=& -\B^{ijj'}_{pQ(q)}-\B^{ijj'}_{Q(p)q}-b_{iQ(pq)}\det(\A)\Bigstrut\\

&=& -\B^{ikk'}_{pQ(q)}-\B^{ikk'}_{Q(p)q}-b_{iQ(pq)}\det(\A), \bigstrut\\

2a_{ip}a_{jq} &=& \phantom{-} \B^{iQ(j)}_{pQ(q)}+\B^{iQ(j)}_{Q(p)q}
+b_{Q(ij)Q(pq)}\det(\A)\Bigstrut\\

&=& \phantom{-} \B^{Q(i)j}_{pQ(q)}+\B^{Q(i)j}_{Q(p)q}+b_{Q(ij)Q(pq)}\det(\A), \bigstrut\\

2a_{iq}a_{jp} &=& \phantom{-} \B^{iQ(j)}_{pQ(q)}+\B^{iQ(j)}_{Q(p)q}
-b_{Q(ij)Q(pq)}\det(\A)\Bigstrut\\

&=& \phantom{-} \B^{Q(i)j}_{pQ(q)}+\B^{Q(i)j}_{Q(p)q}-b_{Q(ij)Q(pq)}\det(\A). \bigstrut
\ea\ee
\ep

\section{From reflections in $\RS$ to rigid-body motions in ${\cal E}^3$}
\setcounter{equation}{0}

Any non-null vector $\X\in \RS$ generate a reflection of $\RS$ by its adjoint action. The induced action in the null 3-spaces
interchanges points and planes, so it is a 3-D projective polarity from $\RF$ to $(\RF)^*$. 

{\bf Notations.} For any $\x=(x,y,z)\in \RT$, denote
\be\ba{lll}
\E(\x) &:=& x\E_1+y\E_2+z\E_3,\\ 
\E'(\x) &:=& x\E_1'+y\E_2'+z\E_3'.
\ea
\ee

\bdf
Any vector $\X\in \RS$ must be of the form $\E(\x)+\E'(\y)$ for $\x, \y\in \RT$.
The vector 
\be
\left(\ba{c}\x \\ \y\ea\right)=(\x, \y)^T
\in {\mathbb R}^3\times {\mathbb R}^3
\ee
is called the {\it screw form} of vector $\X$.
\edf

One advantage of the screw form is its capability of representing two different inner products simultaneously by the same dot
symbol with no ambiguity. On one hand,
$\RT$ is naturally equipped with the standard Euclidean inner product. On the other hand, the inner product of two 
vectors $\X_1=\E(\x_1)+\E'(\y_1)$ and $\X_2=\E(\x_2)+\E'(\y_2)$ of $\RS$ is denoted by the same dot symbol. The two inner products
are distinguished by the participating vectors, with the lower-case letters denoting vectors of $\RT$, while the
capitals denoting vectors of $\RS$. 

We have
\be
\X_1\cdot \X_2=\left(\ba{c}\x_1 \\ \y_1\ea\right)\cdot \left(\ba{c}\x_2 \\ \y_2\ea\right)
=\x_1\cdot \y_2+\y_1\cdot \x_2.
\ee
Vector $\X_1$ is invertible if and only if $\x_1\cdot \y_1\neq 0$.
In $\Lambda(\RF)$, $\E'(\x)$ is the dual of $\E(\x)$. In $\RS$,
\be
\E(\x)\cdot \E'(\y)=\x\cdot \y.
\ee

{\bf Notations.} For two vectors $\x, \y\in {\mathbb R}^3$, $\x\cdot \y$ and $\x\times \y$ are the usual
Euclidean inner product and vector cross product in vector algebra.
For a $3\times 3$ matrix $\M=(\m_1\ \,\m_2\ \,\m_3)$ where each $\m_i\in \RT$, denote
\be
\x\times \M:=(\x\times \m_1\ \ \ \x\times \m_2\ \ \ \x\times \m_3).
\ee
Alternatively, matrix $\x\times \M$ is defined by
\be
(\x\times \M)\y=\x\times (\M\y),\ \hbox{ for } \y\in \RT.
\ee
For example, when $\x=(x,y,z)^T$,
\be
\x\times \I_{3\times 3}
=\left(\ba{rrr}
0 &\ \ -z &\ \ y \\
z &\ \ 0 &\ \ -x \\
-y &\ \ x &\ \ 0
\ea\right).
\ee

Let $\U=(\x, \y)^T\in \RS$ be invertible, then for any $\V=(\p, \q)^T\in \RS$,
\be
Ad^*_\U(\V)=\V-2\frac{\V\cdot \U}{\U\cdot \U}\U=
\left(\ba{c}
\ds \p-\frac{\x\cdot \q+\p\cdot \y}{\x\cdot \y}\x\\

\ds \q-\frac{\bigstrut \x\cdot \q+\p\cdot \y}{\x\cdot \y}\y
\ea\right).
\ee
By direct computation, we get

\bp \label{reflection:one}
For invertible $\U=(\x, \y)^T\in \RS$, 
\be
(\wedge^3 Ad^*_\U)\J_3=-\frac{{\cal F}(\x)}{\x\cdot \y},\hskip 1cm
(\wedge^3 Ad^*_\U){\cal F}(\x)=-(\x\cdot \y)\J_3.
\ee
In other words, the plane at infinity is mapped to the point at infinity $\x$ by the adjoint action of
$\U$. Furthermore,
with respect to the basis $\e_i$ of $\RF$ and $\ce_j$ of $(\RF)^*$, the $4\times 4$ matrix form
of $Ad^*_\U$ of unit determinant is
\be
\pm \frac{1}{\sqrt{|\x\cdot \y|}}
\left(\ba{cc}
0 &\ \ -{\y}^T \\
\y &\ \ \x\times \I_{3\times 3}
\ea\right).
\ee
The matrix is skew-symmetric. Its inverse-transpose is
\be
\pm \frac{1}{\sqrt{|\x\cdot \y|}}
\left(\ba{cc}
0 &\ \ -{\x}^T \\
\x &\ \ \y\times \I_{3\times 3}
\ea\right).
\ee
\ep

$Ad^*_\U$ induces an affine transformation of ${\cal E}^3$ if and only if it preserves the plane at infinity.
In $\Lambda(\RS)$, the plane at infinity is represented by $\J_3=\E_{1'2'3'}$. The following is a direct corollary
 of the above proposition.

\bc
Let $\U_i=(\x_i,\y_i)^T\in \RS$ for $i=1,2$ be invertible.
Then $Ad^*_{\U_1\U_2}$ induces an affine transformation of ${\cal E}^3$ if and only if 
$\x_1\times \x_2=0$. If $\x_1=\x_2=\x$, then
the affine transformation is volume-preserving if and only if $\x\cdot \y_1=\x\cdot \y_2$.
\ec

\bdf
Any element $\M\in \Lambda^2(\RF)=\RS$ induces a linear transformation of $\RF$ as following:
\be
{\cal L}(\M): \
\X\in \RF\mapsto \M\cdot \X\in \RF, 
\label{bilinear:2}
\ee
where the inner product is in $\CL(\RF)$.
It is called the {\it bilinear form} of vector $\M\in \RS$.
\edf

The matrix form of the transformation ${\cal L}(\M)$ is obviously skew-symmetric. 
The following lemma is by direct verification:

\bl
For any 
$\M=(\x, \y)^T\in \Lambda^2(\RF)$, 
\be
{\cal L}\left(\left(\ba{c}
\x \\
\y
\ea\right)\right)=-
\left(\ba{cc}
0 &\ \ -\x^T \\
\x &\ \ \y\times \I_{3\times 3}
\ea\right).
\ee
When $\M\in \RS$ is invertible, the above matrix
is exactly the $4\times 4$ matrix of $Ad^*_\M$ from the space of planes $(\RF)^*$ to the space of points $\RF$.
\el

The inner product in the space of $4\times 4$ matrices is that of the embedded space ${\mathbb R}^{16}$ of components.
For any $\M_1, \M_2\in \RS$,
\be
{\cal L}(\M_1)\cdot {\cal L}({\cal J}\M_2):=
\tr\left(({\cal L}(\M_1))^T({\cal L}({\cal J}\M_2))\right)
=2(\M_1\cdot \M_2).
\ee

\bc
Let $\M=(\x, \y)^T\in \RS$. When $\x\cdot \y\neq 0$, then ${\rm ker}({\cal L}(\M))=0$; when $\x\cdot \y=0$ but $\x, \y$ are not both zero, then
${\rm ker}({\cal L}(\M))$ is the 2-space of $\RF$ containing all the projective points
of the line with Pl\"ucker coordinates
$(\x, \y)^T$.
\ec

An element of $Spin(3,3)$ in factored form with respect to the Clifford product has the following nice property:

\bp
Let $\U=\V_1\V_2\cdots \V_{r}\in Pin(3,3)$, where each $\V_i\in \RS$ is invertible.
Let $\tilde{\U}$ be obtained from $\U$ by replacing the 
odd-positioned vector factors $\V_r, \V_{r-2}, \ldots$ counted from the right in the product
with ${\cal J}\V_r, {\cal J}\V_{r-2}, \ldots$. Let the result be $\tilde{\U}=\W_1\W_2\cdots \W_r$.
Then the $4\times 4$ matrix form of $Ad^*_{\U}$ is up to scale that of the following transformation:
\be
\X\in \RF \mapsto
\W_{1}\cdot(\W_{2}\cdot (\cdots\, \cdot(\W_r\cdot \X)))\in \RF,
\ee
where the inner product is in $\CL(\RF)$.
\ep

\bdf
Let $\U=\X_1\X_2\cdots \X_{r}\in Pin(3,3)$. The {\it inverse-transpose} of $\U$ is the following element of $Pin(3,3)$:
\be
\U^{\cal J}:=({\cal J}\X_1)({\cal J}\X_2)\cdots ({\cal J}\X_{r}).
\ee
\edf

\bl
Let $\U\in Pin(3,3)$. Then as matrices in $O(3,3)$, 
\be
(Ad^*_\U)^{-T} = Ad^*_{\U^{\cal J}}.
\ee
\el

{\it Proof}. For any $\M\in O(3,3)$, $\M^T{\cal J}\M=\cal J$. So $\M^{-T}={\cal J}\M{\cal J}$. For $\M=Ad^*_\U$,
for any $\X\in \RS$, 
\[
\M^{-T}\X={\cal J}(\M(\X^{\cal J}))=Ad^*_{\U^{\cal J}}((\X^{\cal J})^{\cal J})=Ad^*_{\U^{\cal J}}\X.
\]
\endproofs

Below we investigate 3-D projective transformations induced by two reflections in $R^{3,3}$.
Let $\U=(\x_1,\y_1)^T$ and $\V=(\x_2,\y_2)^T$, then $Ad^*_{\U^{\cal J}\V}$ has the following $4\times $ matrix form up to scale:
\be\ba{ll}
& \ds 
-\left(\ba{cc}
0 &\ \ -{\x_1}^T \\
\x_1 &\ \ \y_1\times \I_{3\times 3}
\ea\right)
\left(\ba{cc}
0 &\ \ -{\y_2}^T \\
\y_2 &\ \ \x_2\times \I_{3\times 3}
\ea\right)\\

=& \Bigstrut
\ds 
\left(\ba{cc}
\x_1\cdot \y_2 &\ \ (\x_1\times \x_2)^T \\
-\y_1\times \y_2 &\hskip .7cm (\x_2\cdot \y_1)\I_{3\times 3}+\x_1{\y_2}^T-\x_2{\y_1}^T
\ea\right).
\ea
\label{compose:UTV}
\ee
The resulting matrix has trace $2(\x_1\cdot \y_2+\x_2\cdot \y_1)$, so it
is in $sl(4)$ if and only if $\U\cdot \V=0$.

When $\x_1=\x_2=\x$, (\ref{compose:UTV}) becomes
\be
\left(\ba{cc}
\x\cdot \y_2 &\ \ 0 \\
-\y_1\times \y_2 &\ \ \ \ (\x\cdot \y_1)\I_{3\times 3}+\x(\y_2-\y_1)^T
\ea\right).
\label{affine:form}
\ee
For $\x=\e_1$, and $\y_1=z_1\e_1+z_2\e_2$, and $\y_2=y_1\e_1+y_2\e_2+y_3\e_3$, (\ref{affine:form})
becomes
\be
\left(\ba{cccc}
y_1 &\ \ 0 &\  \ 0 &\  \ 0 \\
-z_2y_3 &\  \ y_1 &\  \ y_2-z_2 &\  \ y_3\\
z_1y_3 &\  \ 0 &\  \ z_1 &\  \ 0 \\
z_2y_1-z_1y_2 &\ \  0 &\ \  0 &\  \ z_1
\ea\right),
\label{mat:factor}
\ee
where the submatrix $\A$ composed of the last three rows and the last three columns
is upper-triangular. $\A$ is orthogonal if and only if 
\[
y_1=\pm z_1,\ \ \ \, y_3=0,\ \ \ \,y_2=z_2, 
\]
{\it i.e.}, 
either $\y_1=\y_2$ or $\y_2=-z_1\e_1+z_2\e_2$. In the former case, (\ref{mat:factor}) equals
$\I_{4\times 4}$ up to scale; in the latter case,
(\ref{mat:factor}) becomes
\be
-z_1\left(\ba{cccc}
1 &\ 0 &\ 0 &\ 0 \\
0 &\ 1 &\ 0 &\ 0\\
0 &\ 0 &\ -1 &\ 0 \\
2z_2 &\ 0 &\ 0 &\ -1
\ea\right).
\ee
The $3\times 3$ submatrix composed of the last three rows and the last three columns represents a rotation of angle $\pi$ about the axis
in direction $\e_1$ and passing through point $z_2\e_3$. Thus we get the following:

\bp
Let $\x\in \RT$ be nonzero, and $\y\in \RT$ satisfies $\x\cdot \y=0$. The rotation of angle $\pi$ about the axis passing through
point $\y$ in direction $\x$ is generated by two reflections induced by the following vectors sequentially (from right to left):
\be
\left(\ba{c}
\x\\
\lambda \x+\x\times \y
\ea\right),\hskip .4cm
\left(\ba{c}
\x\\
-\lambda \x+\x\times \y
\ea\right),
\ee
where $\lambda\neq 0$ and is arbitrary.
Except for such rotations in ${\cal E}^3$, 
no other rigid-body motion is generated by two reflections in $\RS$.
\ep


There are five kinds of 2-planes in $\RS$ that are spanned by invertible vectors, 
classified by the signature of the 2-plane:
${\mathbb R}^{2,0,0}$, 
${\mathbb R}^{0,2,0}$, 
${\mathbb R}^{1,0,1}$, 
${\mathbb R}^{0,1,1}$, 
${\mathbb R}^{1,1,0}$. We use them to give the normal form of any element of $SL(4)$ that is the product of two invertible
skew-symmetric matrices, under the matrix similarity transformations of $SL(4)$.

\bdf
Two matrices $\M, \N\in SL(4)$ are said to be {\it conjugate in $SL(4)$}, if there exists a matrix $\C\in SL(4)$ such that
$\M=\C\N\C^{-1}$.
\edf

\bl
Let $\U, \U_0, \V, \V_0\in \RS$ be invertible, such that $\U^2=\U_0^2$, and $\U\wedge \V$ has
the same signature with $\U_0\wedge \V_0$. Then there exists $g\in SO_0(3,3)$ such that $g(\U)=\U_0$, and 
$g(\V)\in \U_0\wedge \V_0$. 
\label{g:lemma}
\el

{\it Proof.}
Obviously there are two 4-spaces ${\cal V}^4$ and
${\cal V}^4_0$ of $\RS$, each having
signature ${\mathbb R}^{2,2}$, such that $\U, \V\in {\cal V}^4$ and
$\U_0, \V_0\in {\cal V}^4_0$, and the orthogonal map $h$ defined on $\U\wedge \V$ by $h(\U)=\U_0$ and 
$h(\V)\in \U_0\wedge \V_0$ can be extended to an orthogonal map from ${\cal V}^4$ to
${\cal V}^4_0$, still denoted by $h$. 

Let 
${\cal V}^2$ be the orthogonal complement of ${\cal V}^4$ in $\RS$, and let
${\cal V}^2_0$ be the orthogonal complement of ${\cal V}^4_0$ in $\RS$. Let $\v_+, \v_-$ be
an orthogonal basis of ${\cal V}^2$ such that $\v_+^2=1$ and $\v_-^2=-1$, and let
$\v_{0+}, \v_{0-}$ be
an orthogonal basis of ${\cal V}^2_0$ such that $\v_{0+}^2=1$ and $\v_{0-}^2=-1$.
Extend $h$ to an orthogonal transformation of $\RS$ by setting
$h(\v_+)=\v_{0+}$ and $h(\v_-)=\v_{0-}$. Set $g=h$ on ${\cal V}^4$, and
\bi
\item if $h$ is the composition of the reflections with respect to even number of positive
vectors and even number of negative vectors, set $g(\v_+)=\v_{0+}$ and $g(\v_-)=\v_{0-}$;
 
\item if $h$ is the composition of the reflections with respect to odd number of positive
vectors and odd number of negative vectors, set $g(\v_+)=-\v_{0+}$ and $g(\v_-)=-\v_{0-}$;

\item if $h$ is the composition of the reflections with respect to even number of positive
vectors and odd number of negative vectors, set $g(\v_+)=\v_{0+}$ and $g(\v_-)=-\v_{0-}$;

\item if $h$ is the composition of the reflections with respect to odd number of positive
vectors and even number of negative vectors, set $g(\v_+)=-\v_{0+}$ and $g(\v_-)=\v_{0-}$.
\ei
Then $g$ is the composition of the reflections with respect to even number of positive
vectors and even number of negative vectors, so $g\in SO_0(3,3)$. \endproof

\bl
Let $\U=\X_1\X_2\dots \X_{2k}\in Spin(3,3)$, where each $\X_i\in \RS$ is invertible.
Let $g\in SO_0(3,3)$ be the Pl\"ucker transform of $\G\in SL(4)$. Denote
$g(\U):=(g\X_1)(g\X_2)\dots (g\X_{2k})$. Let the $4\times 4$ matrix forms of $Ad^*_\U$ and
$Ad^*_{g(\U)}$ be $\A, \C\in SL(4)$ respectively. Then
\be
\A=\G\C\G^{-1}.
\ee
\el

{\it Proof.} Denote $g\X_i=\Y_i$ for $1\leq i\leq 2k$.
Let $g=Ad^*_\W$ for some $\W\in Spin_0(3,3)$,
then $g^{-1}=Ad^*_{\W^{-1}}$. We have
\[\ba{lll}
Ad^*_{\U}(\X) &=& (g^{-1}(\Y_1\Y_2\cdots \Y_{2k}))\X(g^{-1}(\Y_1\Y_2\cdots \Y_{2k}))^{-1} \\

&=& \W(\Y_1\Y_2\cdots \Y_{2k})\W^{-1}\X\W(\Y_1\Y_2\cdots \Y_{2k})^{-1}\W^{-1}\bigstrut \\

&=& (Ad^*_{\W}\circ Ad^*_{g(\U)}\circ Ad^*_{\W^{-1}}) \X.\bigstrut\\
\ea
\]
\endproofs

\bp \label{normalform:2}
Let $\A$ be the product of two skew-symmetric matrices of $SL(4)$. Then $\A$ is conjugate in 
$SL(4)$ to one and only one of the following block diagonal matrices up to scale:
\bu
\item 
\be
\left(\ba{cccc}
c&\ \ -s&\ \ 0&\ \ 0\\
s&\ \ c&\ \ 0&\ \ 0\\
0&\ \ 0&\ \ c&\ \ -s\\
0&\ \ 0&\ \ s&\ \ c
\ea\right),\
\hbox{ where $c=\cos\theta$, and $s=\sin\theta\neq 0$.}
\label{matrix:1}
\ee

\item 
\be
\left(\ba{cccc}
c&\ \ s&\ \ 0&\ \ 0\\
-s&\ \ c&\ \ 0&\ \ 0\\
0&\ \ 0&\ \ c&\ \ -s\\
0&\ \ 0&\ \ s&\ \ c
\ea\right),\
\hbox{ where $c=\cos\theta$, and $s=\sin\theta\neq 0$.}
\label{matrix:2}
\ee

\item 
\be
\left(\ba{cccc}
1&\ \ -\lambda&\ \ 0&\ \ 0\\
0&\ \ 1&\ \ 0&\ \ 0\\
0&\ \ 0&\ \ 1&\ \ 0\\
0&\ \ 0&\ \ \lambda&\ \ 1
\ea\right),\
\hbox{ where $\lambda\neq 0$.}
\label{matrix:3}
\ee

\item 
\be
\left(\ba{cccc}
1&\ \ \lambda&\ \ 0&\ \ 0\\
0&\ \ 1&\ \ 0&\ \ 0\\
0&\ \ 0&\ \ 1&\ \ 0\\
0&\ \ 0&\ \ \lambda&\ \ 1
\ea\right),\
\hbox{ where $\lambda\neq 0$.}
\label{matrix:4}
\ee 

\item 
\be
{\rm diag}(\lambda^{-1}, \lambda^{-1}, \lambda, \lambda),\
\hbox{ where $\lambda> 0$.}
\label{matrix:5}
\ee

\item 
\be
{\rm diag}(\lambda^{-1}, \lambda^{-1}, -\lambda, -\lambda),\
\hbox{ where $\lambda> 0$.}
\label{matrix:6}
\ee 
\eu
\ep

{\it Proof.}
When $\A=\pm \I_{4\times 4}$, then it belongs to the class (\ref{matrix:5}). Below we assume $\A\neq \pm \I_{4\times 4}$.
Then $\A$ must be the $4\times 4$ matrix form of $Ad^*_{\U\V}$ for some unit vectors $\U, \V \in \RS$,
where $\U\wedge \V\neq 0$.

Case 1. If $\U\wedge \V$ has signature ${\mathbb R}^{2,0,0}$, choose
\be
\V_0=\frac{1}{\sqrt{2}}\left(\ba{c}
\e_2\\
\e_2
\ea\right),\ \ \
\U_0=\frac{c}{\sqrt{2}}\left(\ba{c}
\e_2\\
\e_2
\ea\right)+
\frac{s}{\sqrt{2}}\left(\ba{c}
\e_3\\
\e_3
\ea\right), 
\ee
where $c=\cos\theta$, and $s=\sin\theta\neq 0$. They are both unit vectors, and $\U_0\wedge \V_0$ has signature ${\mathbb R}^{2,0,0}$.
By Lemma \ref{g:lemma}, there exists an element $g\in SO_0(3,3)$ such that
$g(\V)=\V_0$, and $g(\U)\in \V_0\wedge \left(\ba{c}\e_3\\ \e_3\ea\right)$. 
Obviously there exists a parameter $\theta$ such that $g(\U)=\U_0$.
So the $4\times 4$ matrix form of $Ad^*_{\U\V}$ is conjugate in $SL(4)$ to that of $Ad^*_{\U_0\V_0}$.
By (\ref{compose:UTV}), the $4\times 4$ matrix form of $Ad^*_{\U_0\V_0}$ is
(\ref{matrix:1}).

Case 2. 
If $\U\wedge \V$ has signature ${\mathbb R}^{0,2,0}$, then choose unit vectors
\be
\V_0=\frac{1}{\sqrt{2}}\left(\ba{r}
\e_2\\
-\e_2
\ea\right),\ \ \ 
\U_0=
\frac{c}{\sqrt{2}}\left(\ba{r}
\e_2\\
-\e_2
\ea\right)+
\frac{s}{\sqrt{2}}\left(\ba{r}
\e_3\\
-\e_3
\ea\right), 
\ee
where $c=\cos\theta$, and $s=\sin\theta\neq 0$. Obviously $\U_0\wedge \V_0$ has signature ${\mathbb R}^{0,2,0}$, so
an element $g\in SO_0(3,3)$ maps $\V$ to $\V_0$, and maps $\U$ to $\U_0$ for some parameter $\theta$.
The $4\times 4$ matrix form of $Ad^*_{\U_0\V_0}$ is
(\ref{matrix:2}).

Case 3.  If $\U\wedge \V$ has signature ${\mathbb R}^{1,0,1}$, then choose unit vectors
\be
\V_0=\frac{1}{\sqrt{2}}\left(\ba{c}
\e_2\\
\e_2
\ea\right),\ \ \
\U_0=\frac{1}{\sqrt{2}}\left(\ba{c}
\e_2\\
\e_2
\ea\right)
+\frac{\lambda}{\sqrt{2}}
\left(\ba{c}
\e_3\\
0
\ea\right),
\ee
where $\lambda\neq 0$. Since $\U_0\wedge \V_0$ has signature ${\mathbb R}^{1,0,1}$, 
an element $g\in SO_0(3,3)$ maps $\V$ to $\V_0$, and maps $\U$ to $\U_0$ for some parameter $\lambda$.
The $4\times 4$ matrix form of $Ad^*_{\U_0\V_0}$ is
(\ref{matrix:3}).

Case 4.  If $\U\wedge \V$ has signature ${\mathbb R}^{0,1,1}$, then choose unit vectors
\be
\V_0=\frac{1}{\sqrt{2}}\left(\ba{r}
\e_2\\
-\e_2
\ea\right),\ \ \
\U_0=\frac{1}{\sqrt{2}}\left(\ba{r}
\e_2\\
-\e_2
\ea\right)
+\frac{\lambda}{\sqrt{2}}
\left(\ba{c}
\e_3\\
0
\ea\right),
\ee
where $\lambda\neq 0$. The $4\times 4$ matrix form of $Ad^*_{\U_0\V_0}$ is
(\ref{matrix:4}).

Case 5.  If $\U\wedge \V$ has signature ${\mathbb R}^{1,1,0}$, there are four subcases: (1) $\U, \V$ are both positive,
(2) $\U$ is negative while $\V$ is positive, (3) $\U$ is positive while $\V$ is negative, (4) $\U, \V$ are both negative.

Let $\V'\in \U\wedge \V$ be a unit vector orthogonal to
$\V$. Then $\V, \V'$ have opposite signatures, and
$\U\cdot (\V\wedge \V')$ is a unit vector of $\U\wedge \V$ having opposite signature to $\U$. Since
\[
\U\V=\U(\V\wedge \V')(\V\wedge \V')\V=-\V^2 \{\U\cdot (\V\wedge \V')\} \V',
\]
we have $Ad^*_{\U\V}=Ad^*_{\{\U\cdot (\V\wedge \V')\}\V'}$. So we only need to consider the first two subcases by 
assuming $\V^2=1$. Then $(\V')^2=-1$. Let $\U=a\V+b\V'$ where $b\neq 0$. Then $|a^2-b^2|=1$, so $a\neq \pm b$.

Choose unit vectors
\be
\V_0=\frac{1}{\sqrt{2}}\left(\ba{c}
\e_1\\
\e_1
\ea\right),\ \ \
\V_0'=\frac{1}{\sqrt{2}}\left(\ba{r}
\e_1\\
-\e_1
\ea\right).
\ee
An element $g\in SO_0(3,3)$ maps $\V$ to $\V_0$, and maps $\V'$ to $\V_0'$,
so it maps $\U\V$ up to scale to 
\be
\left(\ba{r}
(a+b)\e_1\\
(a-b)\e_1
\ea\right)
\left(\ba{c}
\e_1\\
\e_1
\ea\right).
\label{adjoint:case5}
\ee
The $4\times 4$ matrix form of the adjoint action of (\ref{adjoint:case5}) is
\be
{\rm diag}(a+b, a+b, a-b, a-b). 
\ee
Let
$
\lambda=\sqrt{\left|\frac{a-b}{a+b}\right|}.
$
Then $\lambda > 0$.

Subcase 5.1. When $\U^2=1$, then $a^2-b^2=1$. 
The $4\times 4$ matrix form of the adjoint action of (\ref{adjoint:case5}) is
${\rm diag}(\lambda^{-1}, \lambda^{-1}, \lambda, \lambda)$. 

Subcase 5.2. When $\U^2=-1$, then $a^2-b^2=-1$. 
The $4\times 4$ matrix form of the adjoint action of (\ref{adjoint:case5}) is
${\rm diag}(\lambda^{-1}, \lambda^{-1}, -\lambda, -\lambda)$. 
\endproof

{\it Remark}.
When $\lambda\neq 0$, then
\be
\left(\ba{c}
\x \\
\lambda\x
\ea\right)
\left(\ba{c}
\x \\
\mu\x
\ea\right)
=(\mu+\lambda)+(\mu-\lambda) \E(\x)\wedge \E'(\x).
\label{ratio:const}
\ee
For different parameters $\lambda, \mu$, as long as the ratio $\lambda:\mu$ is constant,
then (\ref{ratio:const}) gives the same spinor up to scale. In particular, as long as $\lambda\neq 0$, then
$
\left(\ba{c}
\x \\
\lambda\x
\ea\right)
\left(\ba{c}
\x \\
-\lambda\x
\ea\right)
$ induces the same rotation of angle $\pi$ about the axis
in direction $\x$ at the origin.





The following is a general criterion on the ``compressibility" of a spinor in factored form.

\bp
Let $\U$ be the Clifford product of 4 invertible vectors of $\RS$. Then $\U$ can be written as the
Clifford product of 2 invertible vectors of $\RS$ if and only if $\langle \U \rangle_4=0$ and
$\langle \U \rangle_2$ has an invertible vector factor.
\ep

{\it Proof.} Let $\U=\X_1\X_2\X_3\X_4$ where each vector $\X_i$ is invertible. 
When $\U=\Y_1\Y_2$ for some other two invertible vectors $\Y_1, \Y_2$, 
the two conditions of the proposition are obviously true. Conversely, when both conditions hold, 
then the $\X_i$ are in some 3-space of $\RS$, so $\langle \U \rangle_2$ must be decomposable into
the outer product of two vectors. Let 
$\a, \b$ be two invertible vector of $\langle \U \rangle_2$ such that 
$\langle \U \rangle_2=\a\wedge \b$. We claim that there are two invertible vectors 
$\Y_i=\lambda_i\a+\mu_i\b$ for $i=1,2$ such that $\U=\Y_1\Y_2$. 

Expanding $\U=\langle \U \rangle_0+\a\wedge \b=(\lambda_1\a+\mu_1\b)(\lambda_2\a+\mu_2\b)$,
we get
\be\ba{lll}
\lambda_1\mu_2-\lambda_2\mu_1 &=& 0, \\
\lambda_1\lambda_2\a^2+\mu_1\mu_2\b^2+\a\cdot \b(\lambda_1\mu_2+\lambda_2\mu_1) &=& \langle \U \rangle_0.
\ea\ee
When setting $\lambda_2=0$, then $\mu_2=\lambda_1^{-1}$ and 
$\mu_1=(\langle \U \rangle_0-\a\cdot \b)\lambda_1\b^{-2}$, so
$\U=(\a+(\langle \U \rangle_0-\a\cdot \b)\b^{-1})\b$. Since $\b$ and $\U$ are both invertible, so
is $\U\b^{-1}$. 
\endproof

In the rest of this section, we construct 4-tuples of reflections in $\RS$ inducing rigid-body motions in ${\cal E}^3$.
Affine transformation (\ref{affine:form}) is translation-free if and only if $\y_1\times \y_2=0$.
Consider the composition of two such affine transformations. Let
\be
\U_1=\left(\ba{c}
\x_1\\
\y_1
\ea\right),\ \ \
\V_1=\left(\ba{c}
\x_1\\
\z_1
\ea\right),\ \ \
\U_2=\left(\ba{c}
\x_2\\
\y_2
\ea\right),\ \ \
\V_2=\left(\ba{c}
\x_2\\
\z_2
\ea\right),
\label{uv:pair}
\ee
where $\y_i=\lambda_i\z_i$ for $i=1,2$ such that $\lambda_i\neq 1$.
For simplicity, assume $\x_i\cdot \z_i=1$ for $i=1,2$.

The $4\times 4$ matrix form of $Ad^*_{\U_1^{-T}\V_1\U_2^{-T}\V_2}$
is up to scale the following:
\be
\left(\ba{cc}
1 &\ \ 0 \\

0 &\ \ \lambda_1\I_{3\times 3}+(1-\lambda_1)\x_1{\z_1}^T 
\ea\right)
\left(\ba{cc}
1 &\ \ 0 \\

0 &\ \ \lambda_2\I_{3\times 3}+(1-\lambda_2)\x_2{\z_2}^T 
\ea\right)
=\left(\ba{cc}
1 &\ 0 \\
0 &\ \M
\ea\right),
\ee
where
\be
\M=\lambda_1\lambda_2\I_{3\times 3}+\lambda_1(1-\lambda_2)\x_2{\z_2}^T
+\lambda_2(1-\lambda_1)\x_1{\z_1}^T+(1-\lambda_1)(1-\lambda_2)(\x_2\cdot \z_1)\x_1{\z_2}^T.
\label{two:UV}
\ee

When $\M$ is the rotation matrix 
\be
\left(\ba{ccc}
\cos \theta &\ -\sin \theta &\ 0\\
\sin \theta &\ \phantom{-} \cos \theta &\ 0 \\
0 &\ 0 &\ 1
\ea\right), 
\label{mat:rot}
\ee
there is more than one way of constructing the four
reflection vectors (\ref{uv:pair}).

Option 1. Choose $\x_1=\e_1$ and $\x_2=\e_2$. Let $\z_1=\e_1+a_2\e_2+a_3\e_3$ and 
$\z_2=b_1\e_1+\e_2+b_3\e_3$. Comparing (\ref{two:UV}) and (\ref{mat:rot}), we get
\be
\lambda_1 = \cos \theta,\ \ \
\lambda_2 = \sec \theta, \ \ \
a_3=b_3=0, \ \ \
a_2=b_1=-\cot\frac{\theta}{2}.
\ee
So
\be\ba{lll}
\U_1 &=& \ds \left(\ba{c}
\e_1\\
\e_1\cos\theta-\e_2\cos\theta\cot\frac{\theta}{2}
\ea\right), \\ 

\Bigstrut
\V_1 &=& \ds \left(\ba{c}
\e_1\cos\frac{\theta}{2}\\
\e_1-\e_2\cot\frac{\theta}{2}
\ea\right), \\ 

\Bigstrut
\U_2 &=& \ds \left(\ba{c}
\e_2\\
\e_2\sec\theta-\e_1\sec\theta\cot\frac{\theta}{2}
\ea\right), \\ 

\Bigstrut
\V_2 &=& \ds \left(\ba{c}
\e_2\\
\e_2-\e_1\cot\frac{\theta}{2}
\ea\right).
\ea
\ee
This 4-tuple of reflection vectors has the strong defect that the following is required: $\theta\notin\{0, \pm \pi/2\}$.

Option 2. Set $\x_1=\e_1$ and $\x_2=c'\e_1+s'\e_2$, where ${\c'}^2+{\s'}^2=1$. Set $\z_1=\e_1+a\e_2$ and
$\z_2=p\e_1+q\e_2$ such that $\x_2\cdot \z_2=cp+sq=1$. Set $\y_1=\lambda \z_1$ and
$\y_2=\lambda^{-1} \z_2$ for some $\lambda\neq 1$. Again by 
comparing (\ref{two:UV}) and (\ref{mat:rot}) and solving the equations, we get
\be
\lambda=-1,\ \ \ \ a=0, \ \ \ \
s'=q=-\sin\frac{\theta}{2},\ \ \ \
c'=p=\cos\frac{\theta}{2}. 
\ee
So
\be\ba{lll}
\U_1 &=& \ds \left(\ba{c}
\phantom{-}\e_1\\
-\e_1
\ea\right), \\ 

\Bigstrut
\V_1 &=& \ds \left(\ba{c}
\e_1\\
\e_1
\ea\right), \\ 

\Bigstrut
\U_2 &=& \ds \left(\ba{c}
\phantom{-}\e_1\cos\frac{\theta}{2}-\e_2\sin\frac{\theta}{2}\\
-\e_1\cos\frac{\theta}{2}+\e_2\sin\frac{\theta}{2}
\ea\right), \\ 

\Bigstrut
\V_2 &=& \ds \left(\ba{c}
\e_1\cos\frac{\theta}{2}-\e_2\sin\frac{\theta}{2}\\
\e_1\cos\frac{\theta}{2}-\e_2\sin\frac{\theta}{2}
\ea\right).
\ea
\label{spinor:rotation}
\ee
This 4-tuple of reflection vectors is perfect. It can be reformulated as follows:

\bp
Let $\x_1, \x_2$ be unit vectors of ${\mathbb R}^3$. The rotation of ${\mathbb R}^3$ in the
$\x_1\x_2$ plane with angle $2\angle (\x_1, \x_2)$ is generated by four reflections in $\RS$
with respect to 
the following invertible vectors sequentially (from right to left):
\be
\left(\ba{r}
\x_2\\
-\x_2
\ea\right), \ \ \
\left(\ba{r}
\x_2\\
\x_2
\ea\right), \ \ \
\left(\ba{r}
\x_1\\
-\x_1
\ea\right), \ \ \
\left(\ba{r}
\x_1\\
\x_1
\ea\right).
\ee
\ep

{\it Remark.} The adjoint action of
$\left(\ba{r}
\x_1\\
-\x_1
\ea\right)
\left(\ba{r}
\x_1\\
\x_1
\ea\right)$ induces the rotation of angle $\pi$ in the plane $\x_1^\perp$. 
It is a geometric fact that this rotation followed by the rotation in the plane $\x_2^\perp$
leads to the rotation of angle $2\angle (\x_1, \x_2)$ in the $\x_1\x_2$ plane.

Next consider the translation by vector $\t$. In the 4 vectors of (\ref{uv:pair}), denote
$\w_i=\y_i-\z_i$ for $i=1,2$, and assume $\x_i\cdot \z_i=1$, then
the matrix form of $Ad^*_{\U_1^{\cal J}\V_1\U_2^{\cal J}\V_2}$
is up to scale the following:
\be\ba{l}
\hskip -.3cm
\ds \left(\ba{cc}
1 &\ \ 0 \\
-\w_1\times \z_1 &\ \ (\x_1\cdot \w_1+1)\I_{3\times 3}-\x_1{\w_1}^T
\ea
\hskip -.1cm\right)
\hskip -.1cm
\left(\ba{cc}
1 &\ \ 0 \\
-\w_2\times \z_2 &\ \ (\x_2\cdot \w_2+1)\I_{3\times 3}-\x_2{\w_2}^T
\ea
\hskip -.1cm\right)\\

=\Bigstrut
\ds
\left(\ba{cc}
1 &\ \ 0 \\
f_1 &\ \ f_2
\ea\right),
\ea
\label{matrix:trans:1}
\ee
where
\be\ba{lll}
f_1 &=& -\w_1\times \z_1-(\x_1\cdot \w_1)\w_2\times \z_2-\w_2\times \z_2+\x_1[\w_1\w_2\z_2], \\
f_2 &=& \bigstrut
(\x_1\cdot \w_1+1)(\x_2\cdot \w_2+1)\I_{3\times 3}
-(\x_1\cdot \w_1+1)\x_2 {\w_2}^T\\

&& \hfill
-(\x_2\cdot \w_2+1)\x_1 {\w_1}^T
+(\x_2\cdot \w_1)\x_1{\w_2}^T.
\ea
\label{translate:expr1}
\ee

Notice that in (\ref{translate:expr1}), vectors $\z_1, \z_2$ occur only in the translation part $f_1$.
That (\ref{matrix:trans:1}) is a translation matrix if and only
if $f_2=\I_{3\times 3}$. For simplicity, choose $\x_1=\x_2=\x$, and set $\x=\e_1$. 

Option 1. 
To simplify the expression of $f_2$, let $\x\cdot \w_i=0$ for $i=1,2$. Then equation $f_2=\I_{3\times 3}$
becomes
\be
-\x(\w_2+(\x\cdot \w_1)\w_2-\w_1)^T=0,
\label{translate:tensor}
\ee
hence $\w_1, \w_2$ are linearly dependent, and we get
$\w_2=-\w_1$. Choose $\w_1=\e_2$.

The translation vector is $\t=f_1=-\w_1\times (\z_1-\z_2)$. We can choose
$\t$ to be in the direction of $\e_1$, so that $\t=-\lambda \e_1$ for some scalar $\lambda$. 
Choose $\z_2=\e_1$, then $\z_1=\e_1+\lambda \e_3$.
The 4 reflection vectors inducing the translation by vector $-\lambda \e_1$ are 
the following from right to left:
\be
\left(\ba{c}
\e_1\\
\e_1+\lambda\e_3+\e_2
\ea\right), \ \ \
\left(\ba{c}
\e_1\\
\e_1+\lambda\e_3
\ea\right), \ \ \
\left(\ba{c}
\e_1\\
\e_1-\e_2
\ea\right), \ \ \
\left(\ba{c}
\e_1\\
\e_1
\ea\right).
\label{translate:opt1}
\ee
This 4-tuple can be formulated in a coordinate-free form as following:

\bp
Let $\x, \y$ be nonzero vectors of ${\mathbb R}^3$ satisfying
$\x\cdot \y=0$. The translation of ${\mathbb R}^3$ by vector $\x$
is generated by four reflections induced by the following invertible vectors
in $\RS$ sequentially (from right to left):
\be
\left(\ba{c}
\x\\
\ds \x-\y+\frac{\x\times \y}{\y^2}
\ea\right), \ \ \
\left(\ba{c}
\x\\
\ds \x+\frac{\x\times \y}{\y^2}
\ea\right), \ \ \
\left(\ba{c}
\x\\
\x+\y
\ea\right), \ \ \
\left(\ba{c}
\x\\
\x
\ea\right).
\ee
\ep

{\it Remark}. The first two reflections generate the following affine transformation:
\be
\left(\ba{cc}
1 & \ \ 0 \\
-\x\times \y+\x&\ \ 
\I_{3\times 3}+\x\y^T
\ea\right);
\label{translate:shear1}
\ee
the last two reflections generate the following affine transformation:
\be
\left(\ba{cc}
1 & \ \ 0 \\
\x\times \y&\ \ 
\I_{3\times 3}-\x\y^T
\ea\right).
\label{translate:shear2}
\ee
The two matrices are always commutative.

Geometrically, when setting $|\y|=1$, then
(\ref{translate:shear2}) represents a planar shear transformation followed by a
translation by vector $\x\times \y$, where the shearing occurs in plane
$(\x\times \y)^\perp$, with shear direction $\x$ and shear factor $-|\x|$: 
for any $\z\in \RT$, $\z\mapsto \z-(\z\cdot \y)\x+\x\times \y$. Similarly, 
(\ref{translate:shear1}) represents a planar shear transformation followed by a
translation by vector $\x-\x\times \y$, where the shearing occurs in direction $\x$ of plane $(\x\times \y)^\perp$,
while the shear factor is $|\x|$. 

Option 2. 
Choosing $\x\cdot \w_i=-2$ for $i=1,2$ simplifies the expression of $f_2$ just the same. Then equation $f_2=\I_{3\times 3}$
becomes
\be
-\x(\w_2+(\x\cdot \w_1)\w_2+\w_1)^T=0,
\label{translate:tensor:2}
\ee
hence $\w_1, \w_2$ are linearly dependent, and we get
$\w_2=\w_1$. Choose $\w_1=\w_2=-2\e_1$. 

The translation vector is $\t=f_1=-\w_1\times (\z_1-\z_2)$. We can choose
$\t$ to be in the direction of $\e_3$, so that $\t=-\lambda \e_3$ for some scalar $\lambda$. 
Choose $\z_2=\e_1$, then 
$\z_1=\e_1-\lambda \e_2/2$.
The 4 reflection vectors inducing the translation by vector $-\lambda \e_3$ are 
the following from right to left:
\be
\left(\ba{c}
\e_1\\
\ds -\e_1-\frac{\lambda}{2}\e_2
\ea\right), \ \ \
\left(\ba{c}
\e_1\\
\ds \e_1-\frac{\lambda}{2}\e_2
\ea\right), \ \ \
\left(\ba{c}
\e_1\\
-\e_1
\ea\right), \ \ \
\left(\ba{c}
\e_1\\
\e_1
\ea\right).
\label{translate:opt2}
\ee
This 4-tuple is just as perfect as the previous option. Its coordinate-free formulation is the following:

\bp
Let $\x, \y$ be nonzero vectors of ${\mathbb R}^3$ satisfying
$\x\cdot \y=0$. 
The translation of ${\mathbb R}^3$ by vector $\x$
is generated by four reflections induced by the following invertible vectors
in $\RS$ sequentially (from right to left):
\be
\left(\ba{c}
\y\\
\ds -\y+\frac{\x\times \y}{2}\Bigstrut
\ea\right), \ \ \
\left(\ba{c}
\y\\
\ds \y+\frac{\x\times \y}{2}\Bigstrut
\ea\right), \ \ \
\left(\ba{c}
\y\\
-\y
\ea\right), \ \ \
\left(\ba{c}
\y\\
\y
\ea\right).
\ee
\ep

{\it Remark.} The first two reflections generate the following affine transformation:
\be
\left(\ba{cc}
1 & \ \ 0 \\
\x&\ \ 
\ds -\I_{3\times 3}+2\frac{\y\y^T}{\y^2}
\ea\right);
\ee
it is the rotation of angle $\pi$ about the axis in direction $\y$ at point $\x/2\in \RT$.
The last two reflections generate the following affine transformation:
\be
\left(\ba{cc}
1 & \ \ 0 \\
0&\ \ 
\ds -\I_{3\times 3}+2\frac{\y\y^T}{\y^2}
\ea\right).
\ee
It is the rotation of angle $\pi$ about the axis in direction $\y$ at the origin.
The two matrices are not commutative.

We work out the spinor representation in factored form of
a general rigid-body motion $(\R, \t): 
\x\in {\mathbb R}^3\mapsto \R\x+\t$, where $\R$ is a rotation matrix, and $\t\in {\mathbb R}^3$.
It is a classical result that a rigid-body motion in space is either a pure translation or 
a {\it screw motion}: a rotation about an affine line called the {\it screw axis}, and then a translation along
the screw axis by a {\it screw driving distance}.

\bl
Let $\R$ be a rotation of angle $\theta$ about the axis in unit direction $\v$ at the origin, and let $\t\cdot \v=0$.
Then the rigid-body motion $\x\mapsto \R\x+\t$ of ${\mathbb R}^3$ is a pure rotation, and the axis of rotation
is the line in direction $\v$ and passing through point
\be
\c=\frac{\R_{\frac{\pi-\theta}{2}}\t}{2\sin\frac{\theta}{2}},
\label{expr:first:center}
\ee 
where $\R_\alpha$ denotes the rotation of angle $\alpha$ about unit direction $\v$ at the origin.
\el

{\it Proof.} We only need to prove that $\c$ is a fixed point. Set 
\be
\v=\e_3,\ \ \
\t=t_1\e_1+t_2\e_2,\ \ \
\c=c_1\e_1+c_2\e_2,\ \ \
\R=\left(\ba{ccc}
\phantom{-}\cos\theta &\ \ \sin\theta&\ \ 0 \\
-\sin\theta &\ \ \cos\theta&\ \ 0 \\
0 &\ \ 0&\ \ 1
\ea\right).  
\ee
Substituting them into $\c-\R\c=\t$, we get
\be
2s\left(\ba{cc}
\phantom{-}s &\ \ c \\
-c &\ \ s
\ea\right)
\left(\ba{c}
c_1\\
c_2
\ea\right)
=\left(\ba{c}
t_1\\
t_2
\ea\right),
\ee
where $c=\cos(\theta/2)$ and $s=\sin(\theta/2)$. The solution is 
\be
\left(\ba{c}
c_1\\
c_2
\ea\right)=\frac{1}{2s}\left(\ba{cc}
s &\ \ -c \\
c &\ \ \phantom{-}s
\ea\right)
\left(\ba{c}
t_1\\
t_2
\ea\right).
\ee

\bc
Let $\R$ be a rotation of angle $\theta$ about the axis in unit direction $\v$ at the origin.
Then the rigid-body motion $\x\mapsto \R\x+\t$ of ${\mathbb R}^3$ is a screw motion, whose screw axis is in direction $\v$ and
passing through point
\be
\c=\frac{\R_{\frac{\pi-\theta}{2}}(\t-(\t\cdot \v)\v)}{2\sin\frac{\theta}{2}}.
\label{point:centertrue}
\ee 
In fact, point $\c$ is the foot drawn from the origin to the screw axis, 
called the {\it original center} of the screw motion. The screw driving distance is 
$
d=\t\cdot \v.
$
\ec

{\it Proof.} Since $\t=(\t-(\t\cdot \v)\v)+(\t\cdot \v)\v$, point $\c$ of (\ref{point:centertrue}) 
is first mapped to $\R\c+(\t-(\t\cdot \v)\v) =\c$ by (\ref{expr:first:center}), then mapped to
$\c+(\t\cdot \v)\v$, so it is on the screw axis.
\endproof

Consider a rigid-body motion $(\R, \t)$ where
the rotation axis of $\R$ is $\e_3$. Let $\e_1$ be rotated to 
$\e_1\cos\theta+\e_2\sin\theta$, and let $\t=t_1\e_1+t_2\e_2+t_3\e_3$ 
be the translation vector. Motivated by the 4-reflection generators of pure rotation
and pure translation,
consider the 4 reflections generated by the following 
vectors:
\be
\U_1=\left(\ba{l}
\e_1 \\
\y_1
\ea\right), \ \ \,
\V_1=\left(\ba{l}
\e_1 \\
\y_2
\ea\right), \ \ \,
\U_2=\left(\ba{c}
c\e_1-s\e_2 \\
\q_1
\ea\right), \ \ \,
\V_2=\left(\ba{c}
c\e_1-s\e_2 \\
\q_2
\ea\right),
\ee
where
\be\ba{lll}
\y_1 &=& 
-\e_1+\alpha_3\e_2-\alpha_2\e_3, \\

\y_2 &=& 
\phantom{-} \e_1+\alpha_3\e_2-\alpha_2\e_3, \\

\q_1 &=& 
-(c\e_1-s\e_2)+\beta_3(s\e_1+c\e_2)-\beta_2\e_3, \\

\q_2 &=&
\phantom{-}(c\e_1-s\e_2)+\beta_3(s\e_1+c\e_2)-\beta_2\e_3.
\ea
\ee
The $4\times 4$ matrix of the affine transformation induced by the above 4 reflections is
\be\ba{l}
\ds
\left(
\ba{cc}
1 &\ \ \ 0 \\
-2(\alpha_3\e_3+\alpha_2\e_2) &\ \ \
-\I_{3\times 3}+2\e_1\e_1^T
\ea
\right)\\

\hskip 1.5cm\Bigstrut
\left(
\ba{cc}
1&\ \ \ 0 \\
-2(\beta_3\e_3+\beta_2(s\e_1+c\e_2)) &\ \ \
-\I_{3\times 3}+2(c\e_1-s\e_2)(c\e_1-s\e_2)^T
\ea
\right)\\

\\
= 
\left(\ba{cc}
1&\ \ \ 0 \\
\t' &\ \ \
\R'
\ea\right),
\ea
\ee
where
\be\ba{lll}
\t' &=& \e_1(-2s\beta_2)
+\e_2(-2\alpha_2+2c\beta_2)
+\e_3(-2\alpha_3+2\beta_3),\\

\R' &=& \Bigstrut
\I_{3\times 3}-2\e_1\e_1^T
+2(c\e_1+s\e_2)(c\e_1-s\e_2)^T
\ =\ 
\left(\ba{ccc}
\cos\theta &\ \ -\sin\theta&\ \ 0 \\
\sin\theta&\ \ \phantom{-} \cos\theta&\ \ 0 \\
0&\ \ 0&\ \ 1 
\ea\right).
\ea
\label{find:general}
\ee

By $\t'=\t$, we get
\be
\alpha_2=-\frac{ct_1+st_2}{2s},\ \ \
\beta_2=-\frac{t_1}{2s}, \ \ \
-\alpha_3+\beta_3=\frac{t_3}{2}.
\label{expr:tnew}
\ee

When $s\neq 0$, {\it i.e.}, $\theta\neq 0$,
we can set
\be
\beta_3=0,\ \ \
\alpha_3=-\frac{t_3}{2}.
\label{opt1:new1}
\ee
Then by (\ref{point:centertrue}) and $d=t_3$,
\be\ba{lll}
\alpha_3\e_2-\alpha_2\e_3 &=& \ds \e_1\times(\c+\frac{d}{2}\e_3), \\

\beta_3(s\e_1+c\e_2)-\beta_2\e_3 &=& \bigstrut
(c\e_1-s\e_2)\times \c.
\ea
\ee
Alternatively, we can set
\be
\alpha_3=0,\ \ \
\beta_3=\frac{t_3}{2}.
\label{opt2:new1}
\ee
Then 
\be\ba{lll}
\alpha_3\e_2-\alpha_2\e_3 &=& \ds \e_1\times \c, \\

\beta_3(s\e_1+c\e_2)-\beta_2\e_3 &=& \bigstrut\ds
(c\e_1-s\e_2)\times (\c-\frac{d}{2}\e_3).
\ea
\ee

When $s=0$, the rigid-body motion is a translation, and there is 
no rotation plane $\e_1\e_2$.
Still
as long as $t_1=0$, we get from (\ref{expr:tnew}) and (\ref{opt2:new1})
$\beta_2=\alpha_3=0$, and $\alpha_2=-t_2/2$, and
$\beta_3=t_3/2$. 
The 4 reflection vectors generating the translation are (from right to left):
\be
\left(\ba{c}
\e_1\\
\ds -\e_1+\frac{t_2}{2}\e_3
\ea
\right),\ \
\left(\ba{c}
\e_1\\
\ds \e_1+\frac{t_2}{2}\e_3
\ea
\right),\ \
\left(\ba{c}
\e_1\\
\ds -\e_1+\frac{t_3}{2}\e_2
\ea
\right),\ \
\left(\ba{c}
\e_1\\
\ds \e_1+\frac{t_3}{2}\e_2
\ea
\right).
\ee
It gives the third 4-tuple of reflection vectors generating pure translation by $\t=t_2\e_2+t_3\e_3$
other than (\ref{translate:opt1}) and (\ref{translate:opt2}).

\bp
Any rigid-body motion is induced by two pairs of reflections 
in $\RS$ such that each pair induces an affine
transformation. For a rotation in a plane spanned by two orthonormal directions $\v_1, \v_2$
with angle $\theta\neq 0$, followed by a translation along vector $\t$, 
the corresponding 4 reflections
are the following (from right to left):
\be
\left(\ba{c}
c\v_1+s\v_2\\
\ds -(c\v_1+s\v_2)+\m_1
\ea\right), \ \ 
\ds\left(\ba{c}
c\v_1+s\v_2\\
\ds c\v_1+s\v_2+\m_1
\ea\right),\ \ 
\ds\left(\ba{c}
\v_1\\
\ds -\v_1+\m_2
\ea\right),\ \ 
\ds\left(\ba{c}
\v_1\\
\ds \v_1+\m_2
\ea\right),
\label{general:euclid}
\ee
where $c=\cos(\theta/2)$ and $s=\sin(\theta/2)$,
and we can choose either
\be
\m_1 =(c\v_1+s\v_2)\times (\c+\frac{d}{2}\n), \ \ \ \
\m_2 = \v_1\times \c,
\label{them:u:fact1}
\ee
or
\be
\m_1 = (c\v_1+s\v_2)\times \c, \ \ \ \
\m_2 = \v_1\times (\c-\frac{d}{2}\n).
\label{them:u:fact2}
\ee
Here $\c$ is the original center (\ref{expr:first:center}) 
of the screw motion, and $d$ is the screw driving distance.
\ep

By (\ref{them:u:fact1}), 
the first two reflections 
generate the rotation of angle $\pi$ about the axis in direction
$c\v_1+s\v_2$ at point $\c+d\n/2$, and the last two reflections
generate the rotation of angle $\pi$ about the axis in direction $\v_1$ at 
point $\c$. By (\ref{them:u:fact2}), 
the first two reflections 
generate the rotation of angle $\pi$ about the axis in direction
$c\v_1+s\v_2$ at point $\c$, and the last two reflections
generate the rotation of angle $\pi$ about the axis in direction $\v_1$ at 
point $\c-d\n/2$. 
When $\theta=0$, for the choice of $\v_1$ such that $\t\cdot \v_1=0$, 
(\ref{general:euclid}) is still valid by taking the limit 
of $\theta\rightarrow 0$. 

\bc
Any rigid-body motion is the composition of two 3-D rotations of angle $\pi$.
\ec

Given a 3-D projective transformation or polarity, 
we see that the elements of $Pin^{sp}(3,3)$ generating it in factored form is far from being unique.
On the other hand, the expanded forms of the elements differ at most by a factor in $\{\pm 1, \pm \I_{3,3}\}$.
For a rigid-body motion, the exponential of a bivector of $\Lambda^2(\RS)$ provides a spinor that when expanded is equal to
the spinors in different factored forms that we find in this section. This is the topic of
the bivector representation of $se(3)$, the Lie algebra of rigid-body motions.

\section{Lie subalgebra of rigid-body motions and classical screw theory}
\setcounter{equation}{0}

First we recall some basic facts of the Lie algebra $se(3)$. 
By direct verification, we know that

\bp
Let $\R$ be the rotation matrix of angle $\theta$ about the axis in unit direction $\v$ at the origin. Then
\be
\left(\ba{cc}
1 &\ \ 0 \\
\t &\ \ \R
\ea\right)
=\exp\left(\ba{cc}
0 &\ \ 0 \\
\u &\ \ \theta\v\times \I_{3\times 3}
\ea\right),
\label{rigid:exp}
\ee
where 
\be
\u=(\t\cdot \v)\v
+\frac{\theta}{2\sin\frac{\theta}{2}}\R_{-\frac{\theta}{2}}(\t-(\t\cdot \v)\v)
=\theta\c\times \v+d\v,
\ee
and the original center $\c$ is given by 
(\ref{point:centertrue}), and $d=\t\cdot \v$.
\ep

The Lie algebra $se(3)$ is composed of $4\times 4$ matrices of the form
$
\left(\ba{cc}
0 &\ \ 0 \\
\t &\ \ \v\times \I_{3\times 3}
\ea\right),
$
where $\t, \v\in \RT$. Hence we can use a pair of vectors $\t, \v$ to represent the matrix.

\bdf
The map
\be
\left(\ba{cc}
0 &\ \ 0 \\
\t &\ \ \v\times \I_{3\times 3}
\ea\right)\in se(3) \mapsto 
\left(\ba{c}
\v \\
\t
\ea\right) \in \RS
\ee 
is called the {\it screw representation} of $se(3)$. The Lie bracket of $se(3)$ agrees with the following
{\it cross product} of screw forms:
\be
\left(\ba{c}
\v_1 \\
\t_1
\ea\right) \times
\left(\ba{c}
\v_2 \\
\t_2
\ea\right)
:=\left(\ba{c}
\v_1\times \v_2 \\
\v_1\times \t_2+\v_2\times \t_1
\ea\right).
\ee
\edf

By (\ref{rigid:exp}),
the screw form of the element of $se(3)$ 
generating the motion $\x\in \RT \mapsto \R\x+\t$ where $\R\neq \I_{3\times 3}$, is
\be
\theta\left(\ba{c}
\v \\
\c\times \v+\lambda \v
\ea\right),
\ee
where $\lambda=d/\theta=(\t\cdot \v)/\theta$ is called the {\it screw ratio}. It is the 
ratio of the screw driving distance $d$ by the angle of rotation $\theta$. The screw motion is positive if $\lambda>0$,
in which case the screw driving direction and the rotation orientation follow the right-hand rule. 
The screw motion is negative if $\lambda<0$, in which case the screw driving direction 
and the rotation orientation follow the left-hand rule. The screw motion is a pure rotation if $\lambda=0$.

Any nonzero vector $\X\in \RS$ is up to scale either of the form 
\be
\left(\ba{c}
0 \\
\t
\ea\right)=\E'(\t), 
\label{ex:screw1}
\ee
or of the form
\be
\X=\left(\ba{c}
\l\\
\c\times \l+\lambda\l
\ea\right),
\label{ex:screw2}
\ee
where $\l$ is a unit vector of $\RT$, $\c\in \RT$ and $\mu\in \mathbb R$. In terms of infinitesimal generators of screw motions,
(\ref{ex:screw1}) generates a pure rotation in direction $\t$. As to (\ref{ex:screw2}), 
since $\X^2=2\lambda$, the vector is
null if $\lambda=0$, and represents a line in direction $\l$ and passing through point $\c$, at the same time it is the generator of
a pure rotation about the axis represented by the 6-D vector itself. 
The vector is positive if $\lambda>0$, and
generates a positive screw motion about the axis $(\l,\c\times \l)^T$; 
the vector is negative if $\lambda>0$, and
represents a negative screw motion about the same axis. So the screw representation of $se(3)$ provides a 
{\it motion-generator interpretation} to all vectors of $\RS$. In contrast, in line geometry only null vectors of $\RS$ have geometric
interpretation.

Traditionally, a ``screw" refers to the screw form of a pure rotation or pure translation, {\it i.e.}, a null vector of $\RS$,
while a general vector of $\RS$ 
is called a {\it twist}, as it generates a screw motion. In this paper, we call the vectors unanimously a screw form.

The following proposition states that the cross product of the screw forms of two screw motions is
the screw form of a third screw motion, the latter being along the common perpendicular of the two screw axes of the former two motions. 

\bp
For $i=1,2$, let
\be
\X_i=\left(\ba{c}
\l_i \\
\c_i\times \l_i+\lambda_i\l_i
\ea\right)\in \RS,
\ee
such that (1) each $\l_i$ is a unit vector, (2) $\c_i\cdot \l_i=0$, and (3) $\l_1\times \l_2\neq 0$.
Let $\y\in \RT$ be an arbitrary point on the common perpendicular of the two lines $(\l_i,\c_i\times \l_i)$, 
and let $d$ be the signed distance from the first line to the second along the direction $\l_1\cdot \l_2$, then
\be
\X_1\times \X_2=\left(\ba{c}
\l_1\times \l_2 \\

\bigstrut
\y\times (\l_1\times \l_2)+\mu \l_1\times \l_2
\ea\right),
\label{geometric:cross}
\ee 
where 
\be
\mu=\lambda_1+\lambda_2+d\frac{\l_1\cdot \l_2}{|\l_1\times \l_2|}.
\ee

For example, the following point of $\RT$ is on the common perpendicular:
\be
\y=\frac{(\c_1\cdot \l_2)\l_2+(\c_2\cdot \l_1)\l_1}{(\l_1\times \l_2)^2},
\label{eval:y:perp}
\ee
and
\be
d=\frac{[(\c_2-\c_1)\l_1\l_2]}{|\l_1\times \l_2|}.
\label{eval:d:perp}
\ee
\ep

{\it Proof.}
Obviously, $\l_1, \l_2, \l_1\times \l_2$ form a basis of $\RT$. Let
\be\ba{lll}
\c_1 &=& p_2(\l_2-(\l_1\cdot \l_2)\l_1)+p_1\l_1\times \l_2, \\
\c_2 &=& q_2(\l_1-(\l_1\cdot \l_2)\l_2)+q_1\l_1\times \l_2.
\ea
\ee
Then (\ref{eval:y:perp}) and (\ref{eval:d:perp}) become
\be
\y=p_2\l_2+q_2\l_1,\ \ \
d=(q_1-p_1)|\l_1\times \l_2|.
\ee
It is easy to verify that line 
$\left(\ba{c}
\l_1\times \l_2 \\
\y\times (\l_1\times \l_2)
\ea\right)$ meets both lines $\left(\ba{c}
\l_i \\
\c_i\times \l_i
\ea\right)$ for $i=1,2$, as the inner products of the vector with the other two are both zero;
so the line is the common perpendicular of the other two lines, and $d$ is the signed distance between the other two lines.

By direction computation, we get that $\X_1\times \X_2$ equals
\[
\left(\ba{c}
\l_1\times \l_2 \\
\l_1\times \l_2(\lambda_1+\lambda_2+(q_1-p_1)\l_1\cdot \l_2)
+\l_1(p_2+q_2\l_1\cdot \l_2)
-\l_2(q_2+p_2\l_1\cdot \l_2)
\ea\right).
\]
By $\y\times (\l_1\times \l_2)=\l_1(p_2+q_2\l_1\cdot \l_2)
-\l_2(q_2+p_2\l_1\cdot \l_2)$, we get (\ref{geometric:cross}).
\endproof

\bp
\be
\left(\left(\ba{c}
\l_1\\
\q_1
\ea\right)\times
\left(\ba{c}
\l_2\\
\q_2
\ea\right)\right)\cdot \left(\ba{c}
\l_3\\
\q_3
\ea\right)
=[\l_1\l_2\q_3]+[\l_2\l_3\q_1]+[\l_3\l_1\q_2].
\ee
The result is invariant under a shift of the subscripts 1,2,3. In particular, when $\q_i=\c_i\times \l_i+\lambda_i\l_i$, the result is
\be
\l_1\cdot (\c_3-\c_2) (\l_2\cdot \l_3) +\l_2\cdot (\c_1-\c_3) (\l_3\cdot \l_1)+\l_3\cdot (\c_2-\c_1) (\l_1\cdot \l_2)
+(\lambda_1+\lambda_2+\lambda_3)[\l_1\l_2\l_3].
\ee
\ep


\bl
[Decomposition of vectors of $\RS$ with respect to the cross product] 
For any $\x, \y\in {\mathbb R}^3$, let $\z\in \RT$ be nonzero such that $\z\cdot \x=0$ and $\z\cdot \y=0$.
Then
\be
\left(\ba{c}
\x \\
\y
\ea\right)
=
\frac{1}{\z^2}
\left(\ba{c}
\z\times \x \\
\z\times \y
\ea\right)
\times
\left(\ba{c}
\z \\
0
\ea\right).
\ee
\el

\bp [Orthogonal decomposition of infinitesimal screw motions of $\RF$] Let $\v_1, \v_2, \v_3$ be an orthonormal basis of $\RT$. Then
\be
\left(\ba{c}
\v_1\\
\lambda \v_3\times \v_1+\mu\v_1
\ea\right) =
\left(\ba{c}
\v_2\\
\lambda \v_3\times \v_2
\ea\right)\times
\left(\ba{c}
\v_3\\
\mu\v_3
\ea\right).
\ee
It states that an infinitesimal screw motion about the axis in direction $\v_1$ and through point $\lambda\v_3$, 
is the cross product of an infinitesimal rotation about the axis in direction $\v_3\times \v_1$ and through point $\lambda\v_3$, 
with an infinitesimal screw motion about the axis in direction $\v_3$ and through the origin.
\ep





We come back to the definition of the cross product in $\RS$. It has some limited covariance in
$SO(3,3)$. Consider $\B\in SO(3,3)$. Let its matrix form with respect to Witt decomposition $\RS=\I_3\oplus \J_3$
be $\left(\ba{cc}
\B_{11} & \B_{12} \\ 
\B_{21} & \B_{22}
\ea\right)$.
Matrix $\B$ is called {\it affine} if its induced transformation in $\RF$ 
preserves the plane at infinity of ${\cal E}^3$, {\it i.e.}, 
$\B_{12}=0$. When $\B$ is affine,
by $\B^T\J\B=\J$, we have
\be
\B_{22}=\B_{11}^{-T},\hskip .7cm \B_{21}^T\B_{11}=-\B_{11}^T\B_{21}.
\ee
When $\B$ is affine, it is {\it Euclidean} if its induced transformation in $\RF$ 
preserves the inner product of $\RT=\langle \e_1, \e_2, \e_3\rangle$,
{\it i.e.}, matrix $\B_{11}$ is orthogonal, and $\B_{22}=\B_{11}$.
If $\B$ induces a translation, then 
$\B_{11}=\B_{22}=\pm \I_{3\times 3}$; if it induces a rotation about the origin, then 
$\B_{21}=0$.





\bp
The cross product is covariant under the Euclidean subgroup of $SO(3,3)$ inducing rigid-body motions of ${\cal E}^3$, {\it i.e.}, for any
\be
\M=\left(\ba{cc}
\R &\ 0 \\
\C\ & \R
\ea\right) \in SO(3,3),
\ee
where $\R$ is a rotation matrix, and $\C=-\A\C^T\A$, 
\be
\left(\M\left(\ba{c}
\v_1 \\
\t_1
\ea\right)\right)
\times
\left(\M\left(\ba{c}
\v_2 \\
\t_2
\ea\right)\right)
=\M\left(\left(\ba{c}
\v_1 \\
\t_1
\ea\right)\times
\left(\ba{c}
\v_2 \\
\t_2
\ea\right)
\right).
\ee
\ep

{\it Proof.}
By direct computation, we get that the cross product is covariant under the subgroup of $SO(3,3)$ inducing
translations, and is also covariant under the subgroup of $SO(3,3)$ inducing rotations
at the origin. So it is covariant under the Euclidean subgroup. \endproof

In classical screw theory, a force is represented by a line in space, and a torque is represented by a line at infinity.
For example, the force along line $(\e_0+\x)\l$ with magnitude $f$ is represented by $f(\e_0\l+\x\l)\in \Lambda^2(\RF)$,
the torque in the plane normal to unit vector $\n$ with scale $t$ is represented by $t\n^\perp\in \Lambda^2(\RF)$.
The linear space spanned by forces and torques is the space of {\it wrenches}. A general {\it wrench} is represented by a
vector of $\RS$, while a pure force or torque is represented by a null vector of $\RS$.

When an infinitesimal rigid-body motion is represented in screw form, {\it i.e.}, as a vector of $\RS$, then the
{\it virtual work} of a wrench and a screw form is defined as the inner product of the two corresponding vectors
in $\RS$. Let $(\f,\q)^T$ be a wrench,
where $\p$ is the force part and $\q$ is the torque part, and let 
$(\v,\u)^T$ be the screw form of a screw motion with axis direction $\v$. 
Then the virtual work between them is
\be
\left(\ba{c}\f\\ \q\ea\right)\cdot
\left(\ba{c}\v\\ \t\ea\right)
=\f\cdot \u+\q\cdot \v.
\ee 

We try to under the virtual work in the $4\times 4$ matrix form of $se(3)$. By (\ref{bilinear:2}), any vector $\X\in \RS$ corresponds to a
bilinear form ${\cal L}(\X)\in sl(4)$. When $\X=(\f,\q)^T$ represents a wrench, then the inner product of ${\cal L}(\X)$ with
the infinitesimal motion matrix $\left(\ba{cc}
0 &\ \ 0 \\
\u &\ \ \v\times \I_{3\times 3}
\ea\right)$ 
is 
\be
-\left(\ba{cc}
0 &\ \ -\f^T \\
\f &\ \ \q\times \I_{3\times 3}
\ea\right)\cdot
\left(\ba{cc}
0 &\ \ 0 \\
\u &\ \ \v\times \I_{3\times 3}
\ea\right)
=-\f\cdot \u-2(\q\cdot \v).
\label{virtual:w}
\ee
It is almost the virtual work except for the ratio of the two part of works by the force and the torque respectively.


Since $Spin(3,3)$ covers $SL(4)$, the Lie algebra $\Lambda^2(\RS)$ of $Spin(3,3)$ is isomorphic to $sl(4)$.
In particular, any element of $se(3)$ has a bivector form in $\Lambda^2(\RS)$. On the other hand, a wrench is only
a vector of $\RS$. How is the virtual work represented by some pairing between a bivector and a vector of $\Lambda(\RS)$?
Before embarking on the investigation of this problem, we make an analysis of the bivector representation of $se(3)$
by computing the expanded form of the spinors (\ref{spinor:rotation}), (\ref{translate:opt1}), (\ref{general:euclid})
inducing rigid-body motions and the corresponding exponential form.

\bp
Let $\x=\e_1$ and 
$\y=c\e_1+s\e_2$, where $c=\cos(\theta/2)$ and $s=\sin(\theta/2)$, then
\be
\ds \frac{1}{4}\left(\ba{c}
\y \\
-\y
\ea\right)
\left(\ba{c}
\y \\
\y
\ea\right)
\left(\ba{c}
\x \\
-\x
\ea\right)
\left(\ba{c}
\x \\
\x
\ea\right)
=
\exp(\frac{\theta}{2}(\E_{21'}-\E_{12'})).
\ee
\ep

{\it Proof.} On one hand,
\be\ba{ll}
& \ds \frac{1}{4}\left(\ba{c}
\y \\
-\y
\ea\right)
\left(\ba{c}
\y \\
\y
\ea\right)
\left(\ba{c}
\x \\
-\x
\ea\right)
\left(\ba{c}
\x \\
\x
\ea\right)
\\

=& \Bigstrut
((c\E_1+s\E_2)\wedge (c\E_1'+s\E_2'))(\E_1\wedge \E_1') \\

=& \bigstrut
c^2+cs(\E_{21'}-\E_{12'})-s^2\E_{121'2'}.
\ea
\ee
On the other hand, 
let $\B=(\E_{21'}-\E_{12'})/2$, let $\I=\E_{12'21'}$, and let
$\L=(1+\I)/2$. Then
\be\ba{lllll}
\I^2 &=& 1,\\
\B\I &=& \I\B &=& \B, \\
\B^2 &=& -\L, \\
\L^2 &=& \L.
\ea
\ee
So
\[\ba{lll}
\exp(\theta\B) &=& \ds 1+\B\theta+\frac{-\L\theta^2}{2!}+\frac{-\B\theta^3}{3!}
+\frac{\L\theta^4}{4!}+\frac{\B\theta^5}{5!}+\ldots \\

&=& \ds\Bigstrut
\B(\theta-\frac{\theta^3}{3!}+\frac{\theta^5}{5!}+\ldots)
+1-\L+\L(1-\frac{\theta^2}{2!}+\frac{\theta^4}{4!}+\ldots) \\

&=& \ds\Bigstrut
\B\sin\theta+\frac{1-\I}{2}+\frac{1+\I}{2}\cos\theta \\

&=& \ds\Bigstrut
\cos^2\frac{\theta}{2}-\I \sin^2\frac{\theta}{2}+2\B \cos\frac{\theta}{2} \sin\frac{\theta}{2}.
\ea
\]
\endproofs



\bp
The translation along vector $\x\in {\mathbb R}^3$ is induced by the following spinor of
$\CL(3,3)$:
\be
1-\frac{\E(\x)\cdot \J_3}{2}=\exp(-\frac{1}{2}\E(\x)\cdot \J_3).
\label{general:translate}
\ee
\ep

{\it Proof.} 
Expanding the Clifford product of the four vectors in (\ref{translate:opt1}) 
from left to right, we get $4+2\lambda \E_{2'3'}$. The exponential form of
spinor $1+\lambda \E_{2'3'}/2$ is obviously $\exp(\lambda \E_{2'3'}/2)$. 
\endproof

\bdf
Let $\x=(x_1, x_2, x_3)\in {\mathbb R}^3$. Define the following linear maps
from $\RT$ to $\Lambda^2(\RS)$:
\be\ba{llll}
\E\E: &\ \ \x=x_1\e_1+x_2\e_2+x_3\e_3 &\mapsto& x_1\E_{23}+x_2\E_{31}+x_3\E_{12}; \\

\E'\E': &\ \ \x=x_1\e_1+x_2\e_2+x_3\e_3 &\mapsto& x_1\E_{2'3'}+x_2\E_{3'1'}+x_3\E_{1'2'}; \bigstrut\\

\E_<\E': &\ \ \x=x_1\e_1+x_2\e_2+x_3\e_3 &\mapsto& x_1\E_{23'}+x_2\E_{31'}+x_3\E_{12'}; \bigstrut\\

\E_>\E': &\ \ \x=x_1\e_1+x_2\e_2+x_3\e_3 &\mapsto& x_1\E_{32'}+x_2\E_{13'}+x_3\E_{21'}.
\ea
\ee
\edf

Obviously, $\E\E(\x)=\E'(\x)\cdot \I_3$, and $\E'\E'(\x)=\E(\x)\cdot \J_3$. 
Under a special orthogonal transformation of $\RT$, the basis $\e_1, \e_2, \e_3$
is changed into another basis, and the basis vectors $\E_i, \E_j'$ of $\RS$ change accordingly.
With respect to the new basis, four new maps can be defined as above. The following invariance can be proved easily:

\bp
The maps $\E\E$, $\E'\E'$ and $\E_>\E'-\E_<\E'$ are all invariant under $SO(3)$.
Furthermore, for any $\x_1, \x_2\in \RT$,
\be
(\E_>\E'-\E_<\E')(\x_1)\cdot (\E_>\E'-\E_<\E')(\x_2)=-2(\x_1\cdot \x_2).
\ee
\ep

The elements of $\E'\E'(\x)$ for $\x\in \RT$ form a Lie subalgebra of $sl(4)$; it is the Lie subalgebra $\RT$ of translations.
The elements of $(\E_>\E'-\E_<\E')(\x)$ for $\x\in \RT$ form the Lie subalgebra $so(3)$. It is easy to prove that any nonzero
element of $so(3)$ in bivector form is not the outer product of two vectors of $\RS$. 

Consider the expanded form of the Clifford product of the 4 vectors in (\ref{general:euclid}):
\be\ba{ll}
& \ds 4\left(\E_{11'}-\frac{t_3}{2}\E_{1'2'}+\frac{c't_1-s't_2}{2s'}\E_{3'1'}\right)
\left(
{c'}^2\E_{11'}+s^2\E_{22'}+c's'(\E_{12'}+\E_{21'})
\phantom{\frac{t_1}{2}}
\right.\\

&\hfill \ds \left.+\frac{c't_1}{2s'}\E_{3'1'}-\frac{t_1}{2}\E_{2'3'}
\right)\\

=& \Bigstrut\ds
4\left(
{c'}^2-s^2\E_{12'21'}-\frac{s}{2}(st_1+ct_2)\E_{11'2'3'}-\frac{s}{2}(st_2-ct_1)\E_{21'2'3'}\right.\\

& \hfill\ds\left.
+cs(\E_{12'}-\E_{21'})-\frac{t_3}{2}\E_{1'2'}
+\frac{c}{2}(st_2-ct_1)\E_{2'3'}-\frac{c}{2}(st_1+ct_2)\E_{3'1'}
\right).
\ea
\label{general:exp:euclid}
\ee

After some computation, we get that
(\ref{general:exp:euclid}) equals $\exp(\B_2)$, where if $\theta\neq 0$, then 
\be\ba{lll}
\B_2= &=& \ds \frac{\theta}{2}(\E_{21'}-\E_{12'})-\frac{t_3}{2}\E_{1'2'}\\

&& \hfill\bigstrut\ds
-\frac{\theta}{4\sin\frac{\theta}{2}}\left((t_1\cos\frac{\theta}{2}+t_2\sin\frac{\theta}{2})\E_{2'3'}
+(-t_1\sin\frac{\theta}{2}+t_2\cos\frac{\theta}{2})\E_{3'1'}\right).
\ea
\ee 
If $\theta=0$, the above expression is still valid by (\ref{general:translate}), as the limit of the right side when $\theta\rightarrow 0$ is
$-(t_1\E_{2'3'}+t_2\E_{3'1'})/2$, and $\e_1, \e_2$ are selected so that 
the translation vector $\t$ is in the $\e_1\e_2$ plane.

\bp
The rigid-body motion $\x\mapsto \R\x+\t$, where $\R$ is the 
rotation matrix of angle $\theta$ about the axis in unit direction $\v$ at the origin, and $\t\in \RT$,
corresponds to the following spinor of $\RS$, where $\c$ is the original center defined by (\ref{point:centertrue}),
and $d=\t\cdot \v$ is the screw driving distance:
\be\ba{l}
\ds
\exp\frac{1}{2}\left\{(\E_>\E'-\E_<\E')(\theta\v)-\E'\E'(\theta\c\times \v+d\v)\right\}.
\ea
\ee
The Lie algebra $se(3)$ in bivector form is spanned by $\E'\E'(\x)$ and $(\E_>\E'-\E_<\E')(\y)$ for $\x, \y\in \RT$.
\ep

\bp
The cross product of $\RS$ is induced by the following trivector of $\Lambda^3(\RS)$:
\be
\C_3={\cal F}(\e_1)+{\cal F}(\e_2)+{\cal F}(\e_3)
=\E_{12'3'}+\E_{23'1'}+\E_{31'2'},
\label{def:C3}
\ee
such that for any $\X, \Y\in \RS$,
\be
\X\times \Y=-(\X\wedge \Y)\cdot \C_3.
\ee
Furthermore, for $\X=(
\x,
\y)^T$,
\be
\X\cdot \C_3=\E'\E'(\y)-(\E_>\E'-\E_<\E')(\x)\in se(3).
\ee
\ep

Let $\L_2=(\E_>\E'-\E_<\E')(\x)-\E'\E'(\y)$ be an element of $se(3)$.
The map
\be
\L_2\mapsto \left(\ba{c} \x \\ \y \ea\right)
\ee 
gives the {\it screw form} of bivector $\L_2\in se(3)$. It is also called the {\it $se(3)$-lift} of
vector $\x,\y)^T$.

The 
transformation $\X\mapsto \L_2\cdot \X$ for $\X\in \RS$ has the following $6\times 6$ matrix form:
\be
\M:=\left(\ba{cc} \x\times \I_{3\times 3} & \ \ 0 \\ 
\y\times \I_{3\times 3}
&\ \ \x\times \I_{3\times 3}
\ea\right).
\ee
For $\X=\left(\ba{c} \p \\ \q \ea\right)\in \RS$,
\be
\L_2\cdot \X=\M\X=
\left(\ba{c} \x\times \p \\ \x\times \q+\y\times \p \ea\right)
=\left(\ba{c} \x \\ \y \ea\right)\times \left(\ba{c} \p \\ \q \ea\right).
\ee

So the cross product of
two vectors of $\RS$ equals the inner product of the $se(3)$-lift of
the first vector with the second vector. When both vectors are lifted to bivectors by the
$se(3)$-lift, then the cross product of the two bivectors of $se(3)$ 
equals the cross product of the two original vectors of $\RS$.

\bp
By the following image of $\C_3$ under $\wedge^3 \cal J$:
\be
\D_3:={\cal F}''(\e_1^*)+{\cal F}''(\e_2^*)+{\cal F}''(\e_3^*)
=\E_{1'23}+\E_{2'31}+\E_{3'12},
\label{def:D3}
\ee
we have
\be
\X\cdot \D_3=\E\E(\x)-(\E_>\E'-\E_<\E')(\y)\in so(3,0,1).
\ee
Furthermore,
\be
\left(\left(\ba{c}
\v\\
\u
\ea\right)\cdot \C_3\right)
\cdot \left(\left(\ba{c}
\f\\
\q
\ea\right)\cdot \D_3\right)=
-\f\cdot \u-2(\q\cdot \v).
\label{virtual:w2}
\ee
\ep

{\it Proof.} That $\X\cdot \D_3\in so(3,0,1)$, the Lie algebra of the special orthogonal group of ${\mathbb R}^{3,0,1}$,
will be made clear in the next section. All others are by simple computation. \endproof

It is easy to prove that both $\C_3$ and $\D_3$ are invariant under any transformation of $SO(\I_3)$ and the induced transformation of
$SO(\J_3)$ with the same matrix form. In this sense, the two trivectors are basis-independent. 

By (\ref{virtual:w}) and (\ref{virtual:w2}), for two vectors $\X, \Y\in \RS$, 
if $\X$ is lifted by $\C_3$ to a bivector $\A_2$ of $se(3)$, and the other is lifted by $\D_3$ to a bivector $\B_2$ of $so(3,0,1)$, then
their inner product in $\Lambda^2(\RS)$ is the inner product of the $4\times 4$ matrix form of $\A_2\in sl(4)$ with the
bilinear form ${\cal L}(\Y)$. The lift by $\D_3$ realizes a pairing 
between a vector and a bivector of $se(3)$, and the result is almost the virtual work between a wrench and an infinitesimal
screw motion.

We have the following more general result on the pairing:

\bp
For $\C_3, \D_3$ given by (\ref{def:C3}) and (\ref{def:D3}), 
\be
\left(\left(\ba{c}
\x_1\\
\y_1
\ea\right)\cdot (\C_3+\lambda \I_3)\right)
\cdot \left(\left(\ba{c}
\x_2\\
\y_2
\ea\right)\cdot (\D_3+\lambda \J_3)\right)=
-2\x_1\cdot \y_2
-(1+\lambda^2)\y_1\cdot \x_2.
\ee
\ep


Having talked enough of $se(3)$, we consider the whole group of {\it Euclidean transformations}, {\it i.e.},
affine transformations of the form $\x\in \RT\mapsto \G\x+\t$, where $\G\in O(3,3)$. For matrix 
\be
\T:=\diag(-1,1,1,1),
\ee
we have
\be
\T\left(\ba{cc}
1 &\ \ \ 0 \\
\t &\ \ \ \G
\ea\right)
=-\left(\ba{cc}
1 &\ \ \ 0 \\
-\t &\ \ \ -\G
\ea\right),
\hskip .4cm
\left(\ba{cc}
1 &\ \ \ 0 \\
\t &\ \ \ \G
\ea\right)\T
=-\left(\ba{cc}
1 &\ \ \ 0 \\
\t &\ \ \ -\G
\ea\right).
\label{T:compose}
\ee
Both results are in $SO(3,3)$, so we can simply compute the spinor generator of each matrix, say $\U_l, \U_r$ respectively,
and then get the generator of $\x\mapsto \G\x+\t$ in ${\cal T}Spin(3,3)$: ${\cal T}\U_l=\U_r{\cal T}$.

For example, consider the 3-D reflections with respect to affine plane $(\n, d)$, the plane normal to unit vector 
$\n$ and whose distance from the origin along $\n$ is $-d$.
The foot drawn from the origin to the plane is $\f=-d\n$. Let $\y\in \RT$. Then the image of the reflection of $\y$
with respect to the plane is $\y-2((\y-\f)\cdot \n)\n=\y-2(\y\cdot \n)\n-2d\n$. The $4\times 4$ matrix form of this
transformation is
\be
\M=\left(\ba{cc}
1 &\ \ \ 0 \\
-2d\n &\ \ \ \I_{3\times 3}-2\n\n^T
\ea\right).
\label{reflection:matrix}
\ee

\bl
The screw motion of rotation angle $\pi$ and screw driving distance $d$ about the axis 
in normal direction $\n$ and passing though point $\c\in \RT$, where $\c\cdot \n=0$,
has the following $4\times 4$ matrix form:
\be
\left(\ba{cc}
1 &\ \ \ 0 \\
2\c+d\n &\ \ \ 2\n\n^T-\I_{3\times 3}
\ea\right).
\label{screw:pi}
\ee
\el

{\it Proof.} By (\ref{point:centertrue}), the translation vector is $\t=2\c+d\n$. 
The rotation matrix of angle $\pi$ about the axis in direction $\n$ at the origin is $2\n\n^T-\I_{3\times 3}$.
\endproof

By (\ref{T:compose}), 
(\ref{reflection:matrix}) and (\ref{screw:pi}), we get

\bc
The reflection with respect to plane $(\n, d)$ is the composition of $\T$ with the screw motion of angle $\pi$ and screw driving
distance $2d$ about the axis in direction $\n$ at the origin.
\ec

\bp
The screw form of $se(3)$ realizing the reflection with respect to plane $(\n, d)$ is $(\pi\n, 2d\n)^T\in \RS$. The corresponding
element of ${\cal T}Spin(3,3)$ is the following: 
\be
{\cal T}\left(\ba{c}
\v_2\\
-\v_2-d\v_1
\ea\right)
\left(\ba{c}
\v_2\\
\v_2-d\v_1
\ea\right)
\left(\ba{c}
\v_1\\
-\v_1
\ea\right)
\left(\ba{c}
\v_1\\
\v_1
\ea\right),
\ee
where $\v_1, \v_2, \n$ form an orthonormal basis of $\RT$.
When $d=0$, the above element can be simplified to 
\be
{\cal T}\left(\ba{c}
\n\\
-\n
\ea\right)
\left(\ba{c}
\n\\
\n
\ea\right).
\ee
\ep

\section{Other Lie subalgebras and corresponding screw forms}
\setcounter{equation}{0}

The 6-D Lie algebra $se(3)$ of rigid-body motions provides an interpretation of vectors of $\RS$ as screw forms of 
infinitesimal rigid-body motions. The inner product of $\RS$ can be transferred to $se(3)$, although incompatible with
the inner product of the bivector representation of $se(3)$; conversely, the cross product of $se(3)$ in bivector form
can be transferred to $\RS$, and the result is the cross product in classical screw theory for line geometry. The 
transformation-generator interpretation of vectors of $\RS$ is not necessarily restricted to the specific Lie subalgebra
$se(3)$. As long as there is a 6-D Lie subalgebra of $sl(4)$, a screw representation may be assigned to the Lie subalgebra, and
the Lie bracket of the Lie subalgebra may be translated to a new cross product of the screw forms.

First we investigate several typical 3-D Lie algebras of $sl(4)$.



1. Perspectivity group:

A {\it perspectivity} \cite{dorst15}
has the following matrix form:
\be
\left(\ba{cc}
1 &\ \ \x^T \\
0 & \ \ \I_{3\times 3}
\ea\right).
\label{special:projective}
\ee
It induces the following rational linear transformation on the line through the origin in unit direction $\l$:
\be
\lambda \l \in \RT \mapsto \frac{\lambda}{1+\lambda \x\cdot \l}\l\in \RT.
\ee
When $\lambda\rightarrow \infty$, the point at infinity of the line is mapped to point $\l/(\x\cdot \l)$.
When $\x\cdot \l=0$, every point on the line is invariant. The origin is fixed, so is every line through the origin.

Every point on the plane normal to $\x$ and through the origin is fixed.
Furthermore,
draw a family of parallel planes normal to $\x$, and let $\lambda$ be the signed distance of the plane from the origin
along $\x$. Then the plane with signed distance $\lambda$ is mapped to the plane with signed distance $\lambda/(1+\lambda d)$.
In particular, the plane with signed distance $-1/d$ is mapped to the plane at infinity, while the plane at infinity is mapped
to the plane with signed distance $1/d$.


\bp
The bivector generator of 
(\ref{special:projective}) in $\Lambda^2(\RS)$ is
$\E'\E'(\x)/2$. The Lie algebra of the group is isomorphic to $\RT$.
\ep

2. Anisotropic dilation group:

The transformation $(1, t_1, t_2, t_3)\in GL(4)$ where $t_1t_2t_3>0$, is called an 
{\it anisotropic dilation} (or {\it non-uniform scaling}). When all the $t_i>0$,
the dilation is said to be {\it positive}. All anisotropic dilations of ${\cal E}^3$ form an Lie group, denoted by $AD(3)$;
all positive anisotropic dilations of ${\cal E}^3$ form a Lie subgroup, denoted by $AD^+(3)$. We have
\be\ba{lll}
AD(3) &=& AD^+(3)\ \cup\  {\rm diag}(1,1,-1,-1)AD^+(3) \\
&&
\cup\ {\rm diag}(1,-1,1,-1)AD^+(3)\ \cup\ {\rm diag}(1,-1,-1,1)AD^+(3). 
\ea
\ee

\bl
Bivector $-\lambda \E_{11'}/2$ induces the following transformation of ${\cal E}^3$:
\be
{\rm diag}(1, 1, e^{\lambda}, e^{\lambda}).
\ee
The corresponding element of $sl(4)$ is $\diag(-\frac{\lambda}{2}, -\frac{\lambda}{2}, \frac{\lambda}{2}, \frac{\lambda}{2})$.
\el

{\it Proof.} Let $\U=\exp(\frac{\lambda}{2} \E_{11'})=\cosh \frac{\lambda}{2}+\E_{11'}\sinh \frac{\lambda}{2}$.
Then $Ad^*_{\U}\E_1=e^\lambda \E_1$ and $Ad^*_{\U}\E_{1'}=e^{-\lambda} \E_{1'}$, while all other basis vectors are fixed.
By direct computing, we get that $Ad^*_{\U}$ has the following $4\times 4$ matrix form:
${\rm diag}(e^{\frac{\lambda}{2}}, e^{\frac{\lambda}{2}}, e^{-\frac{\lambda}{2}}, e^{-\frac{\lambda}{2}})$.\endproof

\bdf
Denote
\be
\F_1=\frac{\E_{22'}+\E_{33'}-\E_{11'}}{2}, \ \ \
\F_2=\frac{\E_{33'}+\E_{11'}-\E_{22'}}{2},\ \ \
\F_3=\frac{\E_{11'}+\E_{22'}-\E_{33'}}{2}.
\label{def:fis}
\ee
The following are linear maps from $\RT$ to $\Lambda^2(\RS)$:
\be\ba{llll}



\E_=\E': &\ \ (x_1,x_2,x_3)^T &\mapsto& x_1(\E_{22'}+\E_{33'})+x_2(\E_{33'}+\E_{11'})+x_3(\E_{11'}+\E_{22'}); \\

\E_=\E'_3: &\ \ (x_1,x_2,x_3)^T &\mapsto& \ds x_1\F_1+x_2\F_2+x_3\F_3.\bigstrut


\ea
\ee
\edf

\bp
Any element of the Lie algebra of $AD^+(3)$ is of the form ${\rm diag}(-(t_1+t_2+t_3), t_1, t_2, t_3)$, and the corresponding
bivector of $\Lambda^2(\RS)$ is 
\be
-\frac{1}{2}\{
t_1(\E_{22'}+\E_{33'})+
t_2(\E_{33'}+\E_{11'})+
t_3(\E_{11'}+\E_{22'})\}.
\ee 
\ep

For example, the isotropic positive dilation ${\rm diag}(1, e^{t}, e^{t}, e^{t})$ is generated by
$-t(\E_{11'}+\E_{22'}+\E_{33'})$.

3. Upper triangular shear transformation group:

An upper triangular shear transformation 
of ${\cal E}^3$ by vector $\t=(t_1, t_2, t_3)^T$ has the following matrix form:
\be
\left(\ba{cccc}
1 &\ \ 0&\ \ 0& \ \ 0\\
0 &\ \ 1&\ \ t_3 &\ \ t_2\\
0 &\ \ 0&\ \ 1 &\ \ t_1\\
0 &\ \ 0&\ \ 0 &\ \ 1
\ea\right).
\label{scissor:def}
\ee

\bp
(\ref{scissor:def}) is generated by the following element of $\Lambda^2(\RS)$: 
\be
\frac{1}{2}(t_3\E_{12'}+t_1\E_{23'}+(t_2-\frac{t_1t_3}{2})\E_{13'}).
\ee
The corresponding element of $sl(4)$ is 
\be
\left(\ba{cccc}
0 &\ \ 0&\ \ 0& \ \ 0\\
0 &\ \ 0&\ \ t_3 &\ \ t_2-\frac{t_1t_3}{2}\\
0 &\ \ 0&\ \ 0 &\ \ t_1\\
0 &\ \ 0&\ \ 0 &\ \ 0
\ea\right).
\ee
\ep

{\it Proof.} Let $\H=x_3\E_{12'}+x_1\E_{23'}+x_2\E_{13'}$. By $\exp(\H)=1+\H-x_1x_3\E_{122'3'}$, we get
the $6\times 6$ matrix form of $Ad^*_{\exp(\H)}$ as following:
\be
\left(\ba{cccccc}
1 &\ \ 2x_3&\ \ 2(x_2+x_1x_3)& 0&\ \ 0&\ \ 0\\
0 &\ \ 1&\ \ 2x_1& 0&\ \ 0&\ \ 0\\
0 &\ \ 0&\ \ 1& 0&\ \ 0&\ \ 0\\
0&\ \ 0&\ \ 0& 1&\ \ 0&\ \ 0 \\
0&\ \ 0&\ \ 0& -2x_3&\ \ 1&\ \ 0 \\
0&\ \ 0&\ \ 0& -2(x_2-x_1x_3)&\ \ -2x_1&\ \ 1 
\ea\right).
\ee
The corresponding $4\times 4$ matrix of $SL(4)$ is
\[
\left(\ba{cccc}
1 &\ \ 0&\ \ 0& \ \ 0\\
0 &\ \ 1&\ \ 2x_3 &\ \ 2(x_2+x_1x_3)\\
0 &\ \ 0&\ \ 1 &\ \ 2x_1\\
0 &\ \ 0&\ \ 0 &\ \ 1
\ea\right).
\]
\endproofs

4. Other triangular shear transformation groups:

In the upper triangular shear transformation (\ref{scissor:def}), basis vector $\e_1$ is fixed, 
$\e_2$ is sheared towards $\e_1$, and $\e_3$ is sheared towards plane $\e_1\e_2$. The roles of $\e_1, \e_2, \e_3$
can be interchanged, resulting in five other branches of shear transformations. The six different branches
of shear transformations are generated by the following six 3-D Lie subalgebras, 
each being represented by its basis:
\be\ba{ll}
\E_{12'}, \E_{13'}, \E_{23'}; &\hskip 1.2cm \E_{12'}, \E_{13'}, \E_{32'}; \\
\E_{21'}, \E_{23'}, \E_{31'}; &\hskip 1.2cm \E_{21'}, \E_{23'}, \E_{13'}; \\
\E_{31'}, \E_{32'}, \E_{12'}; &\hskip 1.2cm \E_{31'}, \E_{32'}, \E_{21'}.
\ea
\ee

For example, the {\it lower triangular shear transformation}
$
\left(\ba{cccc}
1 &\ \ 0&\ \ 0& \ \ 0\\
0 &\ \ 1&\ \ 0 &\ \ 0\\
0 &\ \ t_3&\ \ 1 &\ \ 0\\
0 &\ \ t_2&\ \ t_1 &\ \ 1
\ea\right)
$ is generated by the following element of $\Lambda^2(\RS)$: 
\be
\frac{1}{2}(t_1\E_{32'}+t_3\E_{21'}+(t_2-\frac{t_1t_3}{2})\E_{31'}).
\ee

5. 3-D Lie subalgebra $sl(2)$:

In ${\cal E}^3$ let $\X\in \RF$ represent a point, and let $\Pi\in (\RF)^*$ represent a plane not incident to the point.
The subgroup of $SL(4)\cup SL^-(4)$ fixing point $\X$ and plane $\Pi$ is the group of 2-D projective transformations of the plane with respect to
the point, and is isomorphic to $SL(3)\cup SL^-(3)$. 

This embedding of $SL(3)\cup SL^-(3)$ into $SL(4)$ may be suitable for modeling a pin-hole camera. Let $\C, \X\in \RF$, and let 
$\Pi\in \Lambda^3(\RS)$ represent a plane, then the intersection of line $\C\X$ with the plane is represented by
the bivector $(\C\X)\cdot \Pi\in \Lambda^2(\RS)$. A null 2-space of $\Pi$ represents a unique point on the plane.
It can be taken as the local representation of the spatial point on the plane. The global representation of the point
is $(\C\X)\wedge ((\C\X)\cdot \Pi)\in \Lambda^3(\RS)$.

In the special case where the point is the origin, and the plane is the one at infinity,
The 8-D Lie subalgebra $sl(3)$ is spanned by the following elements of $\Lambda^2(\RS)$:
$\E(\x)\wedge \E'(\y)$ for $\x, \y\in \RT$.

When vectors $\e_0, \e_3$ are fixed, and the 2-space $\langle \e_1, \e_2\rangle$ is also fixed, the corresponding
Lie subgroup of $SL(3)$ is $SL(2)$. The Lie subalgebra $sl(2)$ has 
the following bivector basis:
$
\E_{12'},\
\E_{21'},\  
\E_{11'}-\E_{22'}. 
$


\bl
The Pl\"ucker transform induces the following Lie algebraic isomorphism from $sl(4)$ to $so(3,3)$:
\be
\M=\left(\ba{cc}
-\tr(\N) &\ \ \n_0^T \\
\m_0 &\ \ \N
\ea\right)
\in sl(4)\mapsto 
\left(\ba{cc}
-\tr(\N)\I_{3\times 3}+\N &\ \ -\n_0\times \I_{3\times 3} \\
\m_0\times \I_{3\times 3} &\ \ \tr(\N)\I_{3\times 3}-\N^T
\ea\right),
\label{expr:lift46}
\ee
where $\N=(\m_1\ \ \m_2\ \ \m_3)$, and $\m_i, \n_0\in \RT$ for $0\leq i\leq 3$.
\el

{\it Proof.}
As $\M$ is the derivative of $\exp(t\M)$ at $t=0$, for $\X, \Y\in \RF$, 
\be\ba{lll}
(\exp(t\M)\X)\wedge (\exp(t\M)\Y)
&=& (\X+t\M\X+o(t))\wedge (\X+t\M\Y+o(t))\\

&=& \X\wedge \Y+t(\M\X\wedge \Y+\X\wedge \M\Y)+o(t) \\

&=& \X\wedge \Y+t(\M\otimes \I_{4\times 4}+\I_{4\times 4}\otimes \M)(\X\wedge \Y)+o(t),
\ea
\ee
where 
\be
\M\otimes \I_{4\times 4}(\X\wedge \Y)=\frac{1}{2}\M\otimes \I_{4\times 4}(\X\otimes \Y-\Y\otimes \X)
=\frac{1}{2}(\M\X\otimes \Y-\M\Y\otimes \X).
\ee
So 
the lift of $\M\in sl(4)$ to $so(3,3)$ is the $6\times 6$ matrix representation of the action of
$\M\otimes \I_{4\times 4}+\I_{4\times 4}\otimes \M$ upon $\Lambda^2(\RF)=\RS$. For any positive permutation $ijk$ of 123, by
\[\ba{lll}
(\M\otimes \I_{4\times 4}+\I_{4\times 4}\otimes \M)\e_{0i}
&=& -\tr(\N)\e_{0i}+\m_0 \e_i+\e_0 \m_i, \\

(\M\otimes \I_{4\times 4}+\I_{4\times 4}\otimes \M)\e_{jk}
&=& -\n_0 \e_i+\e_j\m_k-\e_k \m_j,\bigstrut
\ea\]
we get (\ref{expr:lift46}).
\endproof

\bdf
The following are notations of some linear maps from $\RT$ to $3\times 3$ matrices:
\be\ba{llll}
{\rm diag}: &\ \ \x=x_1\e_1+x_2\e_2+x_3\e_3 &\mapsto& {\rm diag}(x_1, x_2, x_3); \\
\\

{\rm skew}: &\ \ \x=x_1\e_1+x_2\e_2+x_3\e_3 &\mapsto& \ds\left(\ba{ccc}
0&\ x_3&\ 0 \\
0&\ 0&\ x_1\\
x_2&\ 0&\ 0
\ea\right); \\
\\

{\rm skew}^T: &\ \ \x=x_1\e_1+x_2\e_2+x_3\e_3 &\mapsto& 
\ds\left(\ba{ccc}
0&\ 0&\ x_2 \\
x_3&\ 0&\ 0\\
0&\ x_1&\ 0
\ea\right)=({\rm skew}(\x))^T. 
\ea
\label{def:skewt}
\ee
\edf

\bp
The Pl\"ucker transform and the adjoint action of $Spin(3,3)$ induce the following
The Lie algebraic isomorphism from $sl(4)$ to $\Lambda^2(\RS)$: let 
\be
\N={\rm diag}(\n_1)+{\rm skew}(\n_2)+{\rm skew}^T(\n_3),
\ee
and let $\m_0, \n_i\in \RT$ for $0\leq i\leq 3$, then
\be\ba{l}
\M=\left(\ba{cc}
-\tr(\N) &\ \ \n_0^T \\
\m_0 &\ \ \N
\ea\right)
\in sl(4)\mapsto \ds \frac{1}{2}\{
-\E'\E'(\m_0)+\E\E(\n_0)-\E_=\E'(\n_1)\\

\hfill +\E_<\E'(\n_2)+\E_>\E'(\n_3)\}.
\label{map:46}
\ea
\ee
\ep

{\it Proof.} The right side of
(\ref{map:46}), denoted by $\B_2\in \Lambda^2(\RS)$, induces a linear transformation
$\X\in \RS\mapsto \B_2\cdot \X\in \RS$. It has the same matrix form (\ref{expr:lift46}).
\endproof

\bdf
Given an $n\times n$ symmetric real matrix $\K$, the set of matrices $\M\in SL(n)$ satisfying
$\M^T\K\M=\K$ is denoted by $SO(\K)$, called the {\it special orthogonal group with respect to $\K$}.
Its Lie algebra is denoted by $so(\K)$.
\edf

\bp
Let $\K$ be a $4\times 4$ symmetric real matrix of rank $\geq 3$. Then
$so(\K)$ is isomorphic to one of $so(4), so(3,1), so(3,0,1), so(2,2), so(2,1,1)$.
\ep

{\it Proof.} Let $\R\in O(4)$ be a matrix to diagonalize $\K$ by similarity transformation:
$\R^T\K\R=\L$, where $\L$ is a real diagonal matrix whose nonzero entries are the eigenvalues of $\K$.
Let $\T$ be a real diagonal matrix to change each nonzero entry of $\L$ to $\pm 1$ by 
matrix congruent transformation: $\T^T\L\T=\N$, where $\N$ is a real diagonal matrix whose nonzero entries are $\pm 1$.
Up to sign, $\N$ is one of
\[
\diag(0,1,1,1), \ \ \diag(0,-1,1,1), \ \ 
\I_{4\times 4},  \ \  \diag(-1,1,1,1), \ \  \diag(-1,-1,1,1),
\]
and 
\be
(\R\T)^T\K(\R\T)=\N. \label{compare:K}
\ee

Let $\M\in so(\K)$. Then $\exp(\M^T)\K\exp(\M)=\K$. It can be written in the following form by using
(\ref{compare:K}):
\be
\{(\R\T)^{-1}\exp(\M)(\R\T)\}^T\N\{(\R\T)^{-1}\exp(\M)(\R\T)\}=\N.
\ee
So $(\R\T)^{-1}\exp(\M)(\R\T)\in SO(\N)$, and $so(\K)$ is isomorphic to $so(\N)$.
\endproof

We start to discuss 6-D Lie subalgebras of type $so(\K)$.

6. $so(4)$: 

It is composed of
skew-symmetric matrices, {\it i.e.}, matrices of the form 
\be
\left(\ba{cc}
0 &\ \ \x^T \\
-\x&\ \ \y\times \I_{3\times 3}
\ea\right),
\ee
for $\x, \y\in \RT$. The following is a 
bivector basis of $so(4)$:
\be\ba{lll}
\E_{2'3'}+\E_{23}, &\ \ \ \E_{3'1'}+\E_{31}, &\ \ \ \E_{1'2'}+\E_{12}, \\
\E_{32'}-\E_{23'}, &\ \ \ \E_{13'}-\E_{31'}, &\ \ \ \E_{21'}-\E_{12'}.
\ea
\ee

For $\x=(x_1,x_2,x_3)^T$ and $\y=(y_1,y_2,y_3)^T$ of $\RT$, let
\be\ba{lll}
\Q(\x) &=& x_1(\E_{32'}-\E_{23'})+x_2(\E_{13'}-\E_{31'})+x_3(\E_{21'}-\E_{12'}),\\
\P(\y) &=& y_1(\E_{2'3'}+\E_{23})+y_2(\E_{3'1'}+\E_{31})+y_3(\E_{1'2'}+\E_{12}).
\ea
\ee
If we denote by $\left(\ba{c}\x \\ \y\ea\right)$ the bivector $\Q(\x)+\P(\y)$, then
\be
\left(\ba{c}\x_1 \\ \y_1\ea\right)\times_{so(4)}
\left(\ba{c}\x_2 \\ \y_2\ea\right)
=\left(\ba{c}
\x_1\times \x_2+\y_1\times \y_2\\

\x_1 \times \y_2+\y_1\times \x_2
\ea\right),
\ee
where ``$\times_{so(4)}$" denotes the Lie bracket of $so(4)$.
So $so(4)$ defines its own cross product of screw forms. 

If we set 
$\left(\ba{c}\x \\ \y\ea\right)=\E(\x)+\E'(\y)$, then
the corresponding trivector of $\Lambda^3(\RS)$ realizing the cross product of $so(4)$ is
\be
\E_{123}+\E_{12'3'}+\E_{23'1'}+\E_{31'2'}
={\cal F}(\e_0)+{\cal F}(\e_1)+{\cal F}(\e_2)+{\cal F}(\e_3).
\ee
It corresponds to the matrix $\diag(1,1,1,1)$ of the quadratic form preserved by group $SO(4)$. 

7. $so(3,1)$: 

It is composed of matrices of the form 
\be
\left(\ba{cc}
0 &\ \ \x^T \\
\x&\ \ \y\times \I_{3\times 3}
\ea\right),
\ee
for $\x, \y\in \RT$. It has the following
bivector basis in $\Lambda^3(\RS)$: 
\be\ba{lll}
\E_{2'3'}-\E_{23}, &  \ \ \  \E_{3'1'}-\E_{31},& \ \ \  \E_{1'2'}-\E_{12}, \\
\E_{32'}-\E_{23'}, & \ \ \  \E_{13'}-\E_{31'}, & \ \ \  \E_{21'}-\E_{12'}.
\ea
\ee

Let
\be
\Q(\x)=(\E'\E'-\E\E)(\x),\ \ \
\P(\y)=(\E_>\E'-\E_<\E')(\y).
\ee
If we denote by $\left(\ba{c}\x \\ \y\ea\right)$ the bivector $\Q(\x)+\P(\y)$, then
\be
\left(\ba{c}\x_1 \\ \y_1\ea\right)\times_{so(3,1)}
\left(\ba{c}\x_2 \\ \y_2\ea\right)
=\left(\ba{c}
\x_1\times \x_2-\y_1\times \y_2\\

\x_1 \times \y_2+\y_1\times \x_2
\ea\right),
\ee
where ``$\times_{so(3,1)}$" denotes the Lie bracket of $so(3,1)$.
So $so(3,1)$ defines another cross product of screw forms. 

If we set 
$\left(\ba{c}\x \\ \y\ea\right)=\E(\x)+\E'(\y)$, then
the corresponding trivector realizing the cross product of $so(3,1)$ is
\be
-\E_{123}+\E_{12'3'}+\E_{23'1'}+\E_{31'2'}
=-{\cal F}(\e_0)+{\cal F}(\e_1)+{\cal F}(\e_2)+{\cal F}(\e_3).
\ee
It corresponds to the matrix $\diag(-1,1,1,1)$ of the quadratic form preserved by group $SO(3,1)$. 

8. 
$so(3,0,1)$: 

It is composed of
matrices of the form 
\be
\left(\ba{cc}
0 &\ \ \x^T \\
0 &\ \ \y\times \I_{3\times 3}
\ea\right),
\ee
for $\x, \y\in \RT$. It has the following
bivector basis in $\Lambda^3(\RS)$: 
\be\ba{lll}
\E_{23}, &  \ \ \  \E_{31},& \ \ \  \E_{12}, \\
\E_{32'}-\E_{23'}, & \ \ \  \E_{13'}-\E_{31'}, & \ \ \  \E_{21'}-\E_{12'}.
\ea
\ee

Let
\be
\Q(\x)=\E\E(\x),\ \ \
\P(\y)=(\E_>\E'-\E_<\E')(\y).
\ee
If we denote by $\left(\ba{c}\x \\ \y\ea\right)$ the bivector $\Q(\x)+\P(\y)$, then
\be
\left(\ba{c}\x_1 \\ \y_1\ea\right)\times_{so(3,0,1)}
\left(\ba{c}\x_2 \\ \y_2\ea\right)
=\left(\ba{c}
\x_1\times \x_2\\

\x_1 \times \y_2+\y_1\times \x_2
\ea\right),
\ee
where ``$\times_{so(3,0,1)}$" denotes the Lie bracket of $so(3,0,1)$.
So $so(3,0,1)$ defines the same cross product of screw forms with that of $se(3)$, and the two algebras are isomorphic under 
the transpose of $4\times 4$ matrices.

The trivector realizing the cross product of $so(3,0,1)$ is
\be
\E_{12'3'}+\E_{23'1'}+\E_{31'2'}
={\cal F}(\e_1)+{\cal F}(\e_2)+{\cal F}(\e_3).
\ee
It corresponds to the matrix $\diag(0,1,1,1)$ of the quadratic form preserved by group $SO(3,0,1)$. 
The transpose of matrices from $se(3)$ to $so(3,0,1)$ is realized in 
$\Lambda^2(\RS)$ by $\wedge^2\cal J$.

9. 
$so(2,2)$:

It is composed of matrices of the form
\be
\left(\ba{cccc}
0&\ \ x_1&\ \ x_2&\ \ x_3 \\
-x_1&\ \ 0&\ \ y_3&\ \ -y_2 \\
x_2&\ \ y_3&\ \ 0&\ \ y_1 \\
x_3&\ \ -y_2&\ \ -y_1&\ \ 0
\ea\right).
\ee
It has the following
bivector basis in $\Lambda^3(\RS)$: 
\be\ba{lll}
\E_{2'3'}+\E_{23}, & \ \ \ \E_{3'1'}-\E_{31},& \ \ \  \E_{1'2'}-\E_{12}, \\
\E_{32'}-\E_{23'}, & \ \ \  \E_{13'}+\E_{31'}, & \ \ \  \E_{21'}+\E_{12'}.
\ea
\ee

\bdf
For $\x=(x_1,x_2,x_3)^T$ and $\y=(y_1,y_2,y_3)^T$ of $\RT$, 
\be
\x\times_1 \y:= \left(\ba{c}
x_2y_3-x_3y_2\\
-(x_3y_1-x_1y_3)\\
-(x_1y_2-x_2y_1)
\ea\right),\ \ \ 
\x\times_2 \y:= \left(\ba{c}
-(x_2y_3-x_3y_2)\\
-(x_3y_1+x_1y_3)\\
x_1y_2+x_2y_1
\ea\right).
\ee
\edf

In $so(2,2)$, 
for $\x=(x_1,x_2,x_3)^T$ and $\y=(y_1,y_2,y_3)^T$ of $\RT$, let
\be\ba{lll}
\Q(\x) &=& x_1(-\E_{32'}+\E_{23'})+x_2(-\E_{13'}-\E_{31'})+x_3(\E_{21'}+\E_{12'}),\\
\P(\y) &=& y_1(\E_{2'3'}+\E_{23})+y_2(\E_{3'1'}-\E_{31})+y_3(\E_{1'2'}-\E_{12}).
\ea
\ee
If we denote by $\left(\ba{c}\x \\ \y\ea\right)$ the bivector $\Q(\x)+\P(\y)$, then
\be
\left(\ba{c}\x_1 \\ \y_1\ea\right)\times_{so(2,2)}
\left(\ba{c}\x_2 \\ \y_2\ea\right)
=\left(\ba{c}
\x_1\times_1 \x_2+\y_1\times \y_2\\

\y_1\times_2 \x_2-\y_2 \times_2 \x_1
\ea\right),
\ee
where ``$\times_{so(2,2)}$" denotes the Lie bracket of $so(2,2)$.
So $so(2,2)$ defines its own cross product of screw forms. 

If we set 
$\left(\ba{c}\x \\ \y\ea\right)=\E(\x)+\E'(\y)$, then
the corresponding trivector realizing the cross product of $so(2,2)$ is
\be
-\E_{123}-\E_{12'3'}+\E_{23'1'}+\E_{31'2'}
=-{\cal F}(\e_0)-{\cal F}(\e_1)+{\cal F}(\e_2)+{\cal F}(\e_3).
\ee
It corresponds to the matrix $\diag(-1,-1,1,1)$ of the quadratic form preserved by group $SO(2,2)$.

10. 
$so(2,1,1)$:

It is composed of matrices of the form
\be
\left(\ba{cccc}
0&\ \ x_1&\ \ x_2&\ \ x_3 \\
0&\ \ 0&\ \ y_3&\ \ -y_2 \\
0&\ \ y_3&\ \ 0&\ \ y_1 \\
0&\ \ -y_2&\ \ -y_1&\ \ 0
\ea\right).
\ee
It has the following
bivector basis in $\Lambda^3(\RS)$:  
\be\ba{lll}
\E_{23}, & \ \ \ -\E_{31},& \ \ \  \E_{12}, \\
\E_{32'}-\E_{23'}, & \ \ \  \E_{13'}+\E_{31'}, & \ \ \  \E_{21'}+\E_{12'}.
\ea
\ee

In $so(2,1,1)$, 
for $\x=(x_1,x_2,x_3)^T$ and $\y=(y_1,y_2,y_3)^T$ of $\RT$, let
\be\ba{lll}
\Q(\x) &=& x_1(\E_{32'}-\E_{23'})+x_2(\E_{13'}+\E_{31'})+x_3(\E_{21'}+\E_{12'}),\\
\P(\y) &=& y_1\E_{23}-y_2\E_{31}+y_3\E_{12}.
\ea
\ee
If we denote by $\left(\ba{c}\x \\ \y\ea\right)$ the bivector $\Q(\x)+\P(\y)$, then
\be
\left(\ba{c}\x_1 \\ \y_1\ea\right)\times_{so(2,1,1)}
\left(\ba{c}\x_2 \\ \y_2\ea\right)
=\left(\ba{c}
\x_1\times_1 \x_2\\

\y_1\times_2 \x_2-\y_2 \times_2 \x_1
\ea\right),
\ee
where ``$\times_{so(2,1,1)}$" denotes the Lie bracket of $so(2,1,1)$.
So $so(2,1,1)$ defines its own cross product of screw forms. 

If we set 
$\left(\ba{c}\x \\ \y\ea\right)=\E(\x)+\E'(\y)$, then
the corresponding trivector realizing the cross product of $so(2,1,1)$ is
\be
-\E_{12'3'}+\E_{23'1'}+\E_{31'2'}
=-{\cal F}(\e_1)+{\cal F}(\e_2)+{\cal F}(\e_3).
\ee 
It corresponds to the matrix $\diag(0,-1,1,1)$ of the quadratic form preserved by group $SO(2,1,1)$.

\bdf
The 
{\it Hadamard product} of vectors in $\RT$ is the following multilinear, associative and commutative product:
\be
\left(\ba{c} x_1 \\
y_1\\
z_1
\ea\right)\odot
\left(\ba{c} x_2 \\
y_2\\
z_2
\ea\right)
:=\left(\ba{c} x_1x_2 \\
y_1y_2\\
z_1z_2
\ea\right).
\ee
The {\it $i$-th Hadamard power} of a vector $\x\in \RT$ is 
\be
\odot^i\x:=\x\odot\cdots\odot \x\ \, \hbox{ ($i$ times).}
\ee
The {\it exponential function} of $\x$ in the Hadamard product is denoted by $e^{\odot\x}$.
In particular,
\be
\odot^0\x=
e^{\odot 0}=\left(\ba{c} 1 \\
1\\
1
\ea\right).
\ee
\edf

11. 6-D Lie subalgebra of general anisotropic dilation group:

The general anisotropic dilation group is composed of matrices of the form
\be
\left(\ba{cccc}
1&\ \ 0&\ \ 0&\ \ 0 \\
t_1&\ \ \lambda_1&\ \ 0&\ \ 0 \\
t_2&\ \ 0&\ \ \lambda_2&\ \ 0 \\
t_3&\ \ 0&\ \ 0&\ \ \lambda_3
\ea\right).
\ee
The matrix represents either an anisotropic dilation centering at an affine point, or a pure translation.
The corresponding Lie subalgebra is denoted by $gad(3)$, and is 
composed of matrices of the form
\be
\left(\ba{cccc}
0&\ \ 0&\ \ 0&\ \ 0 \\
u_1&\ \ \mu_1&\ \ 0&\ \ 0 \\
u_2&\ \ 0&\ \ \mu_2&\ \ 0 \\
u_3&\ \ 0&\ \ 0&\ \ \mu_3
\ea\right).
\ee
Its bivector basis is spanned by elements of the form 
$\E'\E'(\x)$, $\E_=\E'(\y)$ for $\x, \y\in \RT$.

Let
\be\ba{lll}
\Q(\x) &=& x_1\E_{2'3'}+x_2\E_{3'1'}+x_3\E_{1'2'}, \\
\P(\y) &=& \ds y_1\F_1+y_2\F_2+y_3\F_3,
\ea
\ee
where the $\F_i$'s are defined by (\ref{def:fis}).
If we denote by $\left(\ba{c}\x \\ \y\ea\right)$ the bivector $\Q(\x)+\P(\y)$, then
\be
\left(\ba{c}\x_1 \\ \y_1\ea\right)\times_{gad}
\left(\ba{c}\x_2 \\ \y_2\ea\right)
=\left(\ba{c}
\x_1\odot \y_2-\y_1\odot \x_2\\

0
\ea\right),
\ee
where ``$\times_{gad}$" denotes the Lie bracket of this Lie algebra.

12. 7-D invariant group of ${\mathbb R}^{2,2}$ and its 6-D subgroups:

Let $\L_4$ be a fixed 4-space of signature ${\mathbb R}^{2,2}$, and let its orthogonal complement in $\RS$ be
$\L'_2$. Then $\L'_2$ has signature ${\mathbb R}^{1,1}$, so its null 1-spaces represent a pair of non-intersecting lines in space,
say $\l_1$ and $\l_2$.  
The null 1-spaces of $\L_4$ is denoted by $N(\L_4)$; it
represents all the lines in space incident to both $\l_1$ and $\l_2$. 

\bp
Any projective point on line $\l_1$ or $\l_2$ is on infinitely many lines of $N(\L_4)$, 
while any other projective point is on one and only one line of
$N(\L_4)$. Similarly, any plane passing through $\l_1$ or $\l_2$ also passes through infinitely many lines of $N(\L_4)$, 
while any other plane passes through one and only one line of
$N(\L_4)$.
\ep

{\it Proof.}
Let $\X_3$ be a null 3-space of $\Lambda(\RS)$. 
When $\X_3$ represents a point of $\l_1$, then it is incident to all the lines connecting $\X_3$ and $\l_2$; when
$\X_3$ represents a plane through $\l_1$, then it meets line $\l_2$ at a point $\Y_3$, and so is incident to all the lines
connecting $\Y_3$ and $\l_1$. 

When $\X_3$ represents a point on neither $\l_1$ nor $\l_2$, then in $\RS$, the intersection of the 3-space $\X_3$ with the 4-space
$\L_4$ is a null 1-space or 2-space. If point $\X_3$ is on more than one line of $N(\L_4)$, then $\l_1, \l_2$ must be coplanar,
violating the requirement that the two lines do not intersect. So point $\X_3$ is on on and only one line of $N(\L_4)$. 
Similarly, when $\X_3$ represents a plane supporting neither $\l_1$ nor $\l_2$, 
it must be incident to one and only one line of $N(\L_4)$. 
\endproof

Consider the subgroup of $SO(3,3)$ that leaves the 4-space $\L_4$ invariant. Denote the subgroup by $Inv(\L_4)$.
Then 
\be
Inv(\L_4)=(SO(1,1)\oplus SO(2,2))\cup (SO^-(1,1)\oplus SO^-(2,2)).
\ee
It is a 7D Lie subgroup. 

A typical example is that $\l_1$ is an affine line, say $\l_1=\e_{01}$, the line through the origin and in direction $\e_1$,
and $\l_2$ is a line at infinity, say $\l_2=\e_{23}$, the line at infinity normal to $\e_1$. Then $N(\L_4)$ is composed of all
lines meeting $\l_1$ and perpendicular to it. Then $\L_4=\E_{22'33'}$, and denote 
\be
Inv(2,2):=Inv(\E_{22'33'}).
\ee
The Lie algebra $inv(2,2)$ of $Inv(2,2)$ is spanned by $\E_{11'}, \E_{22'}, \E_{33'}, \E_{23}, \E_{2'3'}$, $\E_{23'}, \E_{32'}$. 
Its matrix form in $sl(4)$ is
\be
\left(\ba{cc}
\A &\ \ 0 \\
0 &\ \ \D
\ea\right),
\ee
where $\A, \D$ are two $2\times 2$ matrices such that $\tr(\A)+\tr(\D)=0$.
The group
$Inv(2,2)$ is composed of matrices of the form
\be
\left(\ba{cc}
\A &\ \ 0 \\
0 &\ \ \D
\ea\right),\ \ \hbox{ or }\ \
\left(\ba{cc}
0 &\ \ \B \\
\C &\ \ 0
\ea\right),
\label{couple:proj}
\ee
where $\det(\A\D)=1$, and $\det(\B\C)=-1$. 

For a general 4-space $\L_4$ of signature ${\mathbb R}^{2,2}$,
$Inv(\L_4)$ is isomorphic to $(GL(2)\times GL(2))/*$, where two elements are equivalent if and only if they differ by scale.
The space $\L_4$ is called the space of {\it coupled projective screws} with respect to the pair of axes
$\l_1, \l_2$. The {\it projective screw ratio} refers to $\det(\A)/\det(\D)$ in the first case of (\ref{couple:proj}), and
$\det(\B)/\det(\C)$ in the second case.

A typical 6-D Lie subalgebra of $inv(2,2)$ is spanned by 
$\E_{22'}, \E_{33'}, \E_{23}, \E_{2'3'}$, $\E_{23'}, \E_{32'}$. 
For $\x=(x_1,x_2,x_3)^T$ and $\y=(y_1,y_2,y_3)^T$,
let
\be\ba{lll}
\P(\x) &=& x_1(\E_{33'}+\E_{22'})+x_2\E_{23}+\x_3\E_{2'3'}, \\
\Q(\y) &=& y_1(\E_{33'}-\E_{22'})+y_2\E_{32'}+\y_3\E_{23'}.
\ea
\ee
Then
\be\ba{lll}
\P(\x)\times \Q(\y) &=& 0, \\

\P(\x)\times \P(\y) &=& -(x_2y_3-x_3y_2)(\E_{33'}+\E_{22'})+2(x_1y_2-x_2y_1)\E_{23}\bigstrut\\
&& \hfill +2(x_3y_1-x_1y_3)\E_{2'3'}, \\

\Q(\x)\times \Q(\y) &=& \phantom{-}(x_2y_3-x_3y_2)(\E_{33'}-\E_{22'})+2(x_1y_2-x_2y_1)\E_{32'}\bigstrut\\
&& \hfill +2(x_3y_1-x_1y_3)\E_{23'}. 
\ea
\ee

For the whole Lie algebra $sl(4)$, by (\ref{def:skewt}),
any element of it is naturally decomposed into the direct sum of five 3-D vectors:
\be
\M=\left(\ba{cc}
-\tr({\rm diag}(\c)) &\ \ \ \n^T \\
\m &\ \ \ {\rm diag}(\c)+{\rm skew}(\u)+{\rm skew}^T(\d)
\ea\right),
\label{M:decomp}
\ee
where $\m, \n, \c, \u, \d\in \RT$. They allow for defining a matrix of $sl(4)$ as a ``superscrew" in 
$\RT\times\RT\times\RT\times\RT\times\RT$, and defining the Lie bracket as a ``supercross product". 
Such a representation of $sl(4)$ is called the {\it superscrew representation}. For symmetry consideration,
instead of choosing the 3-D Lie subalgebra decomposition (\ref{M:decomp}), we make the following decomposition of $sl(4)$
into 3-spaces:

\bdf
When $sl(4)$ is represented by bivectors of $\Lambda^3(\RS)$, 
a {\it superscrew} of 3-D projective geometry is defined as follows:
for any $\x_i\in \RT$ where $1\leq i\leq 5$,
\be
\left(\ba{c}
\x_1\\
\x_2\\
\x_3\\
\x_4\\
\x_5
\ea\right)
:=
\ba{l}
\\

\E\E(\x_1)+\E'\E'(\x_2)+\E_=\E'_3(\x_3)\Bigstrut\\

\hskip .2cm +(\E_>\E'+\E_<\E')(\x_4)+(\E_>\E'-\E_<\E')(\x_5).\bigstrut
\ea
\ee
\edf

\bdf
The symmetric cross product of two vectors of $\RT$ is defined as follows:
\be
\left(\ba{c} x_1 \\
y_1\\
z_1
\ea\right)\star
\left(\ba{c} x_2 \\
y_2\\
z_2
\ea\right)
=\left(\ba{c} 
y_1z_2+z_1y_2 \\
z_1x_2+x_1z_2\\
x_1y_2+y_1x_2
\ea\right).
\ee
\edf

By direct computation, we get

\bp
The Lie bracket of $sl(4)$, when represented in the superscrew form, is the following:
for $\x_i, \y_j\in \RT$,
\be
\left(\ba{c}
\x_1\\
\x_2\\
\x_3\\
\x_4\\
\x_5
\ea\right)\times_{sl(4)}
\left(\ba{c}
\y_1\\
\y_2\\
\y_3\\
\y_4\\
\y_5
\ea\right)
=\left(\ba{c}
\z_1\\
\z_2\\
\z_3\\
\z_4\\
\z_5
\ea\right),
\ee
where
\be
\ba{lll}
\z_1 &=& -\x_1\odot \y_3+\x_3\odot \y_1+\x_1\star \y_4-\x_4\star \y_1+\x_1\times \y_5+\x_5\times \y_1, \\

\z_2 &=& -\x_2\star \y_4+\x_4\star \y_2+\x_2\odot \y_3-\x_3\odot \y_2+\x_2\times \y_5+\x_5\times \y_2, \bigstrut\\

\z_3 &=& (-\x_1\cdot \y_2+\x_2\cdot \y_1)\e^{\odot 0} -\x_1\odot \y_2+\x_2\odot \y_1
+ 4(\x_4\odot \y_5-\x_5\odot \y_4)\times \e^{\odot 0}, \bigstrut\\

\z_4 &=& \ds \frac{1}{2}(\x_1\star \y_2-\x_2\star \y_1)-(\x_3\times \e^{\odot 0})\odot \y_5
+\x_5\odot (\y_3\times \e^{\odot 0}) \Bigstrut\\

&& \hfill 
-\x_4\times \y_5-\x_5\times \y_4,\\

\z_5 &=& \ds \frac{1}{2}(\x_1\times \y_2+\x_2\times \y_1)+(\x_3\times \e^{\odot 0})\odot \y_4
-\x_4\odot (\y_3\times \e^{\odot 0}) \Bigstrut\\

&& \hfill
+\x_4\times \y_4+\x_5\times \y_5.
\ea
\ee
\ep

\section{Conclusion}
\setcounter{equation}{0}

In this paper, we establish a rigorous mathematical foundation for the line geometric model of 3-D projective geometry. 
We also extend screw theory from rigid-body motions to projective transformations, in hope that non-Newtonian mechanics
may find an algebraic language for the development of virtual work of projective motions. 

The connection between the $\CL(3,3)$ model of 3-D projective geometry and the $\CL(4,1)$ model of 3-D conformal geometry
is an interesting topic, and will be investigated in another paper. This connection makes it possible to investigate 
3-D non-Euclidean geometry together with its various realizations in Euclidean space via pin-hole cameras of ${\mathbb R}^{4,1}$.

\end{document}